\documentclass[11pt]{amsart}

\usepackage[T1]{fontenc}
\usepackage[utf8]{inputenc}
\usepackage{lmodern}
\usepackage[a4paper,margin=28mm]{geometry}
\usepackage{amsmath,amssymb,amsthm,mathtools}
\usepackage{booktabs}
\usepackage{array}
\usepackage{microtype}
\usepackage{float}
\usepackage[section]{placeins}   
\usepackage[hidelinks]{hyperref}

\usepackage{orcidlink}

\newcounter{tablecont}
\newcommand{\continuedtable}{%
  \addtocounter{table}{-1}%
  \stepcounter{tablecont}%
  \renewcommand{\theHtable}{\thetable.c\thetablecont}%
}

\newcommand{\PSLtwoZ}{\ensuremath{\mathrm{PSL}(2,\mathbb Z)}}

\newtheorem{theorem}{Theorem}[section]
\newtheorem{corollary}[theorem]{Corollary}
\newtheorem{proposition}[theorem]{Proposition}
\newtheorem{lemma}[theorem]{Lemma}
\theoremstyle{remark}
\newtheorem{remark}[theorem]{Remark}

\DeclareMathOperator{\arcosh}{arcosh}
\DeclareMathOperator{\dist}{dist}
\DeclareMathOperator{\diam}{diam}

\title[Huber--constant algorithm]{A direct majorant and an improved Huber constant on cofinite Fuchsian groups}
\author[M. Avdispahi\'{c}]{Muharem Avdispahi\'{c}\,\orcidlink{0000-0001-7836-4988}}
\address{Department of Mathematics and Computer Science, Faculty of Science,
University of Sarajevo, Zmaja od Bosne 33--35, 71000 Sarajevo,
Bosnia and Herzegovina}
\email{mavdispa@pmf.unsa.ba}


\begin{document}

\begin{abstract}
To obtain an effective upper bound for Huber’s constant on cofinite Fuchsian groups, we introduce an algorithm that bypasses the quantitative Karamata Tauberian theorem in estimating spectral-counting growth by direct use of re-optimized Gaussian kernel instead. Further steps involve the different choice of a compactly supported test function, keeping linear and constant terms in the spectral estimate, optimizing the exact cusp geometry, and retaining the full explicit primitive-geodesic majorant in subsequent Chebyshev estimates. The more direct path is established to control the difference between the size of the hyperbolic term in the Selberg trace formula and the second Chebyshev function. The self-improvement device combined with knowledge of low-lying spectral data and an anchored form of the propagation step leads to the sharpening that in the case of $\mathrm{PSL}(2,\mathbb{Z})$ reads $C_{\mathrm{Hub},\mathrm{PSL}(2,\mathbb{Z})}<81$.  The improvement factor is greater than  $2.7\cdot10^{6}$ compared to the currently known bound.
\end{abstract}

\maketitle

\section{Introduction}

 Huber's constant $C_{\mathrm{Hub},\Gamma }$
is defined (see \cite{FriedmanJorgensonKramer2011}) as the smallest constant in the prime geodesic theorem on the orbifold $X=\Gamma \setminus \mathcal{H%
}$ such that%
\begin{equation*}
\left\vert \pi _{\Gamma }\left( x\right) -\sum\limits_{\lambda _{k}\leq
\frac{1}{4}}\operatorname{Li}\left( x^{s_{k}}\right) \right\vert \leq C_{\mathrm{Hub},\Gamma }%
\frac{x^{3/4}}{\log x}\text{ for all }x>1\text{.}
\end{equation*}
Here $\Gamma \subset \mathrm{PSL}(2,\mathbb{R})$ is a cofinite Fuchsian group of the first kind acting on 
$\mathcal{H}=\{z=x+iy:y>0\}$, the upper half-plane equipped with the hyperbolic metric
$\frac{\sqrt{dx^{2}+dy^{2}}}{y}$ and volume element $\frac{dxdy}{y^{2}}$; 
$\pi
_{\Gamma }\left( x\right) $ denotes the number of $\Gamma$-inconjugate primitive closed geodesics $\gamma_{0}$ on $X$ with the length 
$\ell_{\gamma_{0}}\leq {\log x}$;  $0\leq \lambda _{k}\leq
\frac{1}{4}$ are small eigenvalues of the Laplace-Beltrami operator $\Delta =-y^{2}\left( \frac{\partial ^{2}}{%
\partial x^{2}}+\frac{\partial ^{2}}{\partial y^{2}}\right) $ on $X$ and $%
s_{k}=\frac{1}{2}+\sqrt{\frac{1}{4}-\lambda _{k}}$; $\operatorname{Li}\left(
x\right) =\int\limits_{2}^{x}\frac{dt}{\log t}$.

Huber's constant is an \textbf{absolute} constant, not an asymptotic one. It finds relevant applications in Arakelov theory of algebraic curves \cite{Arakelov1974}. 

Adapting Randol's proof of the prime geodesic theorem \cite{Randol1984}, Friedman, Jorgenson and Kramer \cite{FriedmanJorgensonKramer2011} designed an algorithm to compute 
an upper bound for $C_{\mathrm{Hub},\Gamma }.$ To get an effective bound, they used a quantitative version of Karamata's Tauberian theorem \cite[Lemma~1.1]{FriedmanJorgensonKramer2011}. For the modular group, the resulting value was 16,607,349,020,658. The corrected and modified algorithm \cite{Avdispahic2021} led to 224,615,135, i.e., approximately 74000-times smaller number. As the first reason for a large discrepancy between the value they obtained and the computational results provided through the PARI/GP program \cite{PARI}, Friedman, Jorgenson, and Kramer pointed out that their bound for the spectral counting function was probably one million times too large (compared to Weyl’s law, an asymptotic result). 

Our treatment of the spectral counting function presented in Section \ref{sec:direct-majorant} below does not involve Karamata's theorem. It uses Gaussian kernel with the group specific optimization parameter $\eta$ to provide a direct majorant in Theorem \ref{thm:direct-gaussian-K}. A local
spectral bound, which will carry the two-piece form of the joint supremum, is given in Proposition \ref{prop:local}.

Section \ref{sec:gaussian-huber-propagation} is devoted to the substantial reduction in the bound of Huber's constant on general cofinite Fuchsian groups. After introducing a new auxiliary test function, the essentially optimal constants $A_8^{\sharp}$ and $A_9^{\sharp}$ and the resulting $A_{10}^{\sharp}$ are calculated. The propagation of the direct majorant, resp. local spectral bound, is reflected in a single piece, resp. the two-pieces, form of the joint supremum in $A_{12}^{\sharp}[\kappa]$, where $\kappa$ is another optimization parameter. The proposition \ref{prop:D2} presents the new way to treat the difference between the size of the hyperbolic term in the Selberg trace formula and the second Chebyshev function $\theta_{\Gamma}(x)$ omitting the previous intermediate step that involved $\psi_{\Gamma}(x)$. Propositions \ref{prop:boot0}  and \ref{prop:boot1} demonstrate how the Huber constant estimate can be further improved through bootstrap and knowledge of a sufficiently large part of the spectrum.
The propagation of the Chebyshev estimate to $\pi_{\Gamma}$ is itself
reorganised.  Lemma~\ref{lem:propfactor} identifies the exact constant in the
inequality that controls $\int_{c_{0}}^{x}y^{3/4}y^{-1}\log^{-2}y\,dy$ against
$x^{3/4}/\log x$; it is $\frac{3T^{\natural}}{3T^{\natural}-4}$, where $T^{\natural}$ is the
unique zero of an explicit increasing function; the four cases of
Theorem~\ref{Theorem2.1} recalled under Preliminaries are majorants of it, two of them evaluations of the
very same expression at the wrong point.  Together with the observation that the
Stieltjes integration may be anchored at any point where $\theta_{\Gamma}$ is
known --- where the boundary term is then exact and may be combined with
$\operatorname{Li}$ rather than majorised --- this removes the additive
constant from the bound altogether. Theorem \ref{huber-final} reflects the degree of improvement over the currently known Theorem \ref{Theorem2.1}. 
 
The power of the method is very well illustrated by its application to \PSLtwoZ\  in Section \ref{sec:psl(2,Z)}: the line of argumentation using a single piece form in the joint supremum results in \(C_{\mathrm{Hub},\PSLtwoZ }<1144\), the  two-piece form lowers the bound to $686$, the bootstrap technique leads to $337$, the knowledge of the respective
part of the primitive length spectrum of \PSLtwoZ\ lowers this to $199$, and
anchoring the propagation at $c_{0}=10^{4}$ gives the conclusion
\(C_{\mathrm{Hub},\PSLtwoZ}<81\) in Theorem~\ref{thm:conclusion}.

The cocompact case is presented in Section~\ref{sec:cocompact} with Theorem \ref{thm:cocompact} as the outcome. Application to the $(2,3,7)$ triangle group, for example, reduces the Huber constant bound by the factor $1.53\cdot10^{5}$, from $9.14137\cdot10^{7}$ in
\cite[Proposition~6.1]{Avdispahic2024Huber} to $597.23$ in the proof of Proposition \ref{prop:237}.

\section{Preliminaries}

The fundamental constants on a general cofinite Fuchsian group needed in our exposition are listed in Table~\ref{tab:fundamental-constants}.

\begin{table}[htbp]
\centering
\caption{Fundamental constants}
\label{tab:fundamental-constants}

\renewcommand{\arraystretch}{1.12}

\begin{tabular}{@{}>{\centering\arraybackslash}p{0.15\textwidth}
                p{0.80\textwidth}@{}}
\toprule
\textbf{Constant} & \textbf{Description} \\
\midrule

$\Gamma$
&
Cofinite Fuchsian group
\\[4pt]

$\mathcal{F}$
&
Fundamental domain of $\Gamma$ with area $|\mathcal{F}|$
\\[4pt]

$\tau$
&
Number of parabolic conjugacy classes in $\Gamma$
\\[4pt]

$\mathcal{F}_i^Y$
&
Cusp sector in decomposition

$\mathcal{F}
=
\mathcal{F}_Y
\cup
\bigcup_{i=1}^{\tau}\mathcal{F}_i^Y$,
\qquad
where \(\mathcal{F}_Y\) is compact

\\[4pt]

$d_i$
&
Hyperbolic diameter of
$K_i=\left\{z\in T_i^{-1}\mathcal{F}:
\varepsilon\leq\operatorname{Im}z\leq Y+1\right\}$, where
$\mathcal{F}_i^Y=T_i\bigl((0,1)\times[Y,\infty)\bigr)$
\\[4pt]

$\theta_R,m_R$
&
$\operatorname{tr}(R)=2\cos\theta_R$,
$\theta_R\in(0,\pi)$;
$m_R=|C(R)|$, $C(R)$ centralizer of $R$
\\[4pt]

$\mathcal{A},\, \mathcal{A'}$
&
Number of small eigenvalues
$0=\lambda_0<\cdots<\lambda_{\mathcal{A}-1}\le\frac14$, resp. $\lambda_{\mathcal{A'}-1}<\frac14,$
of Laplacian $\Delta$
counted with multiplicity\\[4pt]

$s_k$
&
$s_k=\frac12+\sqrt{\frac14-\lambda_k}$, so that $s_0=1$ when $\lambda_0=0$
\\[4pt]

$\Phi(s),\phi(s)$
&
Scattering matrix and scattering determinant $\phi(s)=\det\Phi(s)$;
$\frac14\operatorname{tr}\bigl(I-\Phi(\frac12)\bigr)$ enters $C(\eta,L)$
\\[4pt]

$\mathcal{N}$
&
Number of poles
$\sigma_1>\cdots>\sigma_{\mathcal{N}}$
of scattering determinant $\phi(s)$ in
$\left(\frac12,1\right]$
\\[4pt]

$\ell$
&
The length of the shortest closed geodesic on
$\Gamma\backslash\mathcal{H}$
\\[4pt]

$c$
&
$e^{\ell}$, infimum of norms $NP$ of hyperbolic conjugacy classes
$P$ in $\Gamma$
\\[4pt]

$c_1$
&
Minimum of lower lefthand matrix entries of non-parabolic elements
\\[4pt]

$\gamma$
&
Euler's constant, $\gamma=0.57721\ldots$
\\
\bottomrule

\end{tabular}
\end{table}

The output of the current algorithm to determine an effective upper bound for Huber's constant on $\Gamma$ is given in the following theorem (see \cite[Theorem 4.1.]{Avdispahic2024Huber}).

\begin{theorem}
Let $\Gamma$ be a cofinite Fuchsian group. Then
\[
C_{\mathrm{Hub},\Gamma} \leq R_\Gamma,
\]
where
\begin{align*}
\text{\rm(i)}\quad
R_\Gamma
&=
\frac{4}{3e} A'_{22}
+
A'_{21}
\left(
1+\frac{4}{3\log c}
+\frac{4c^{3/4}}{3e\log c}
\right)
&&\text{if } \ell\leq\frac{4}{3};
\\[6pt]
\text{\rm(ii)}\quad
R_\Gamma
&=
\frac{4}{3e} A_{22}
+
A_{21}
\left(
2+\frac{4}{3\log c}
\right)
&&\text{if } \frac{4}{3}<\ell
\leq\frac{2(1+\sqrt{5})}{3};
\\[6pt]
\text{\rm(iii)}\quad
R_\Gamma
&
=
\frac{4}{3e} A_{22}
+
A_{21}
\frac{3\log c}{3\log c-4}
&&\text{if } \frac{2(1+\sqrt{5})}{3}
\leq \ell \leq\frac{8}{3};
\\[6pt]
\text{\rm(iv)}\quad
R_{\Gamma}
&
=
\frac{4}{3e} A_{22}
+
2A_{21}
&&\text{if } \ell>\frac{8}{3}.
\end{align*}
\label{Theorem2.1}
\end{theorem}

The constant $A_{22}=\frac{\mathcal{A}(c-2)}{\log 2 }$, resp. $A'_{22}=\mathcal{A} \max \{\frac{2(2-c^{1/2})
}{\log c} , \frac{(c-2)
}{\log 2}\},$
while  $A_{21}$, resp. $A'_{21}$, is obtained through the process (see \cite[Table 2 and the page after]{Avdispahic2024Huber}) of finding effective estimates of all terms in the Selberg trace formula (see \cite{Hejhal1976}) 

\begin{equation*}
\begin{split}
\sum\limits_{\lambda _{n}\leq \frac{1}{4}}h\left(
r_{n}\right) +\sum\limits_{\lambda _{n}>\frac{1}{4}}h\left( r_{n}\right) +\frac{1}{4\pi }\int\limits_{%
\mathbb{R}
}h\left( r\right) \frac{-\phi ^{\prime }}{\phi }\left(
\frac{1}{2}+ir\right) dr
=\frac{\left\vert \mathcal{F}\right\vert }{4\pi }\int\limits_{\mathbb{R}%
}rh\left( r\right) \tanh \left( \pi r\right)
dr\\ +\sum\limits_{P}\frac{\log NP_{0}}{\left( NP\right) ^{1/2}-\left( NP\right)
^{-1/2}}g\left( \log NP\right)+\sum\limits_{R}\frac{1}{2m_{R}\sin \theta _{R}}%
\int\limits_{\mathbb{R}}\frac{e^{-2r\theta _{R}}}{1+e^{-2\pi r}}%
h\left( r\right) dr \\
+\frac{h\left( 0\right) }{4}\operatorname{tr}\left( I-\Phi \left( \frac{1%
}{2}\right) \right) -\tau g\left( 0\right) \log 2-\frac{%
\tau }{2\pi }\int\limits_{\mathbb{R}}h\left( r\right)
\psi \left( 1+ir\right) dr,
\end{split}%
\end{equation*}%
where $h\left( z\right) $ is an even holomorphic function in the strip $%
\left\vert \text{Im}z\right\vert \leq \frac{1}{2}+\varepsilon $ satisfying a
growth condition $h\left( z\right) =O\left( \left( 1+\left\vert z\right\vert
\right) ^{-2-\varepsilon }\right) $ and $g\left( x\right) $ is its inverse Fourier
transform. The symbol $\Phi $ denotes the scattering matrix, $\phi $ is the scattering determinant, and $\psi$ is the digamma function (the logarithmic derivative of the gamma function).

A few words on Fourier conventions are in order here. For $f\in L^{1}(\mathbb{R})$, we put
\begin{equation*}
  \mathsf{F} f(r)\;=\;\widehat{f}(r)\;=\;\int_{\mathbb{R}}f(u)e^{iru}\,du,
  \qquad
  \mathsf{F}^{-1}F(u)\;=\;\frac{1}{2\pi}\int_{\mathbb{R}}F(r)e^{-iru}\,dr .
\end{equation*}
Thus $\widehat f(0)=\int_{\mathbb{R}}f$, and with
$(f_1*f_2)(u)=\int_{\mathbb{R}}f_1(u-v)f_2(v)\,dv$ one has
\begin{equation*}
  \widehat{f_1*f_2}=\widehat{f_1}\,\widehat{f_2},
  \qquad
  \int_{\mathbb{R}}|f|^{2}\,du=\frac{1}{2\pi}\int_{\mathbb{R}}|\widehat f|^{2}\,dr .
\end{equation*}
The pair $(g,h)$ entering the Selberg trace formula is
$h=\mathsf{F} g$, equivalently $g=\mathsf{F}^{-1}h$.  All mollifiers are normalised by
\begin{equation*}
  \int_{\mathbb{R}}\varphi=\widehat\varphi(0)=1 .
\end{equation*}

Let
\begin{equation*}
 \phi_\mathcal{N}(s)
 =
 c_1^{2s-1}\phi(s)
 \prod_{j=1}^{\mathcal{N}}
 \frac{s-\sigma_j}{s-1+\sigma_j}
\end{equation*}
be the modified scattering determinant.  Then (see \cite[Theorem~3.5.5]{Venkov1982})
\[
 -\frac{\phi_\mathcal{N}'}{\phi_\mathcal{N}}\left(\frac12+ir\right)\ge0.
\]
One has
\begin{align*}
 U(r)
 &:=
 -\frac{\phi'}{\phi}
 \left(\frac12+ir\right)-
 \left(-
 \frac{(\phi_\mathcal{N})'}{\phi_\mathcal{N}}
 \left(\frac12+ir\right)
 \right)
 \notag\\
 &=
 2\log c_{1}
 -
 \sum_{j=1}^{\mathcal{N}}
 \frac{2\alpha_j}{\alpha_j^2+r^2},
\end{align*}
where $\alpha_j:=\sigma_j-\frac12$.

Consequently, for every test function \(H\), with $0\le H\le1$,
\begin{equation}
 \frac1{4\pi}
 \int_{\mathbb R}
 H(r)|U_\Gamma(r)|\,dr
 \le
 \frac{|\log c_{1}|}{2\pi}
 \int_{\mathbb R}H(r)\,dr
 +
 \frac{\mathcal{N}}{2}.
 \label{eq:cong-U-general-boundnew}
\end{equation}

Finally, recall the following primitive geodesic majorant \cite[pp.~47--49]{Fischer1987}
\begin{equation}
 \pi_\Gamma(x)
 \le\sum_{j=1}^{\min\{2,\tau\}}
 \frac{4\pi}{|K_j|}
 \left(\cosh(\log x+2d_j)-1\right)
 \label{eq:modular-pi-exactnew}
\end{equation}
for $x\ge c.$

\section{A direct positive majorant derived from the re-optimized Gaussian}
\label{sec:direct-majorant}

As already said, our aim here is to replace the Tauberian step in the current algorithm by a direct positive majorant. For that purpose, we shall use a re-optimized Gaussian kernel as a test function.

Let \(R\ge0\) be the spectral cut-off and set
\begin{equation*}
  L:=\sqrt{R^2+\frac14}.
\end{equation*}
For \(\eta>0\), define
\begin{equation*}
 G_{\eta,L}(r)
 :=
 \exp\left(-\eta\frac{r^2}{L^2}\right).
\end{equation*}
Then \(G_{\eta,L}\) is even, entire, and rapidly decreasing in every
horizontal strip.  With the Fourier convention
\[
 g(u)
 =
 \frac{1}{2\pi}\int_{\mathbb R}G(r)e^{-iru}\,dr,
\]
its Fourier transform is
\begin{equation*}
 g_{\eta,L}(u)
 =
 \frac{L}{2\sqrt{\pi\eta}}
 \exp\left(-\frac{L^2u^2}{4\eta}\right).
\end{equation*}
In particular,
\[
 G_{\eta,L}(r)\ge0,
 \qquad
 g_{\eta,L}(u)\ge0.
\]

For \(0\le r\le R\),
\[
 \frac{r^2}{L^2}\le1,
\]
and hence
\begin{equation*}
 G_{\eta,L}(r)\ge e^{-\eta}.
\end{equation*}
Any small eigenvalue \(r=i\beta\), \(0\le\beta\le1/2\), contributes
\begin{equation*}
 G_{\eta,L}(i\beta)
 =
 \exp\left(\eta\frac{\beta^2}{L^2}\right)\ge1.
\end{equation*}
Notice that the constant eigenfunction has \(\lambda_0=0\), hence \(r_0=i/2\), and
\begin{equation*}
 G_{\eta,L}(i/2)
 =
 \exp\left(\frac{\eta}{4L^2}\right).
\end{equation*}

Now, define
\begin{equation*}
  S_\mathcal{N}(G_{\eta,L})
 :=
 \sum_{n=0}^{\infty}G_{\eta,L}(r_n)
 +
 \frac1{4\pi}
 \int_{\mathbb R}
 G_{\eta,L}(r)
 \left(
 -\frac{\phi_\mathcal{N}'}{\phi_\mathcal{N}}\left(\frac12+ir\right)
 \right)\,dr.
\end{equation*}
Introduce
\begin{equation*}
 \mathcal M(R)
 :=
 \#\left\{
 r_n\in\mathbb R:0\le r_n\le R
 \right\}
 +
 \frac1{4\pi}
 \int_{-R}^{R}
 \left(
 -\frac{\phi_\mathcal{N}'}{\phi_\mathcal{N}}
 \left(\frac12+it\right)
 \right)\,dt
\end{equation*}
in order to get a sharper estimate for the original spectral counting function 
\begin{equation*}
 \mu(R)
 =
 \#\left\{
 r_n:0\le r_n\le R,\;
 \lambda_n=\frac14+r_n^2
 \right\}
 +
 \frac1{4\pi}
 \int_{-R}^{R}
 \left|
 \frac{\phi'}{\phi}
 \left(\frac12+it\right)
 \right|dt.
\end{equation*}
The underlying idea is to remove the known contribution of the constant
eigenfunction before converting the trace-formula estimate into a
spectral-counting estimate.

Obviously,
\begin{equation} 
\sum\limits_{\lambda _{n}\leq \frac{1}{4}}G_{\eta,L}\left(
r_{n}\right) +\sum\limits_{\lambda _{n}>\frac{1}{4}}G_{\eta,L}\left( r_{n}\right)\ge G_{\eta,L}\left( r_{0}\right) + \mathcal{A'} - 1 +\sum\limits_{0\le r_n\le R}G_{\eta,L}\left( r_{n}\right)
\end{equation}

Since \(G_{\eta,L}(r)\ge e^{-\eta}\) for \(0\le r\le R\), and since
the modified continuous-spectrum density is nonnegative, one obtains
\begin{equation*}
 e^{-\eta}\mathcal M(R)
 \le
  S_\mathcal{N}(G_{\eta,L})
 -
 \exp\left(\frac{\eta}{4L^2}\right) - \mathcal{A'} + 1.
\end{equation*}

We shall bound \(S_\mathcal{N}(G_{\eta,L})\) directly by the geometric side of the Selberg
trace formula. 

\subsubsection*{Identity term}

Since \(G_{\eta,L}(r)\ge0\)  on  \(\mathbb R\) and
\(|\tanh(\pi r)|\le1\), the identity term satisfies
\begin{align*}
 I(G_{\eta,L})
 &=
 \frac{|\mathcal{F}|}{4\pi}
 \int_{\mathbb R}rG_{\eta,L}(r)\tanh(\pi r)\,dr
 \notag\\
 &\le
 \frac{|\mathcal{F}|}{2\pi}
 \int_0^\infty rG_{\eta,L}(r)\,dr
 =
 \frac{|\mathcal{F}|L^2}{4\pi\eta}.
\end{align*}

\subsubsection*{Elliptic term}

For \(0<\theta<\pi\),
\begin{equation*}
 \int_{\mathbb R}
 \frac{e^{-2\theta r}}{1+e^{-2\pi r}}\,dr
 =
 \frac{1}{2\sin\theta}.
\end{equation*}
As \(0\le G_{\eta,L}(r)\le1\), the contribution of an elliptic class is therefore
at most \(1/(4m\sin^2\theta)\). For future reference, we put 
\begin{equation}
    E=\sum\limits_{R}\frac{1}{4m_{R}\sin^{2} \theta _{R}}
\end{equation}

\subsubsection*{Parabolic and scattering terms}

One has
\begin{equation*}
 0\le
 \frac14
 \operatorname{tr}
 \left(
 I-\Phi_\Gamma\left(\frac12\right)
 \right)\le \frac{\tau}{2}.
\end{equation*}
The parabolic logarithmic term is
\begin{equation*}
 -\tau g_{\eta,L}(0)\log2
 =
 -
 \frac{\tau L\log2}
 {2\sqrt{\pi\eta}}.
\end{equation*}
Since $G_{\eta,L}$ is real and even, and
\[
  \psi(1-ir)=\overline{\psi(1+ir)},
\]
we have
\begin{equation*}
  \int_{\mathbb{R}} G_{\eta,L}(r)\,\psi(1+ir)\,dr
   = \int_{\mathbb{R}} G_{\eta,L}(r)\,\Re\psi(1+ir)\,dr .
\end{equation*}
We use the classical expansion
\begin{equation*}
  \psi(1+z) = -\gamma + \sum_{n=1}^{\infty}\frac{z}{n(n+z)} .
\end{equation*}
For $z=ir$, taking real parts gives
\begin{equation*}
  \Re\psi(1+ir) = -\gamma + \sum_{n=1}^{\infty}\frac{r^{2}}{n\,(n^{2}+r^{2})} .
\end{equation*}
All terms in the sum are nonnegative. Keeping only the term $n=1$, we obtain
\begin{equation}\label{eq:psi-lowerboundnew}
\begin{aligned}
  \Re\psi(1+ir)
  &\ \ge\ -\gamma + \frac{r^{2}}{1+r^{2}} \\[2pt]
  &\ =\ 1-\gamma-\frac{1}{1+r^{2}} .
\end{aligned}
\end{equation}
Multiplying \eqref{eq:psi-lowerboundnew} by $G_{\eta,L}(r)\ge 0$ and integrating yields
\begin{equation*}
  \int_{\mathbb{R}} G_{\eta,L}(r)\,\Re\psi(1+ir)\,dr
  \ \ge\ (1-\gamma)\int_{\mathbb{R}} G_{\eta,L}(r)\,dr
  \ -\ \int_{\mathbb{R}} \frac{G_{\eta,L}(r)}{1+r^{2}}\,dr .
\end{equation*}
Put
\[
 t:=\frac{\eta}{L^2}.
\]
Then
\[
 \int_{\mathbb R}G_{\eta,L}(r)\,dr
 =
 \sqrt{\frac{\pi}{t}},
\]
while the classical Gaussian integral gives
\begin{equation*}
  \int_{\mathbb{R}} \frac{e^{-tr^{2}}}{1+r^{2}}\,dr
   = \pi e^{t}\operatorname{erfc}(\sqrt{t}) .
\end{equation*}
where 
\[
\operatorname{erfc}(x)=\frac{2}{\sqrt{\pi}} \int_x^\infty e^{-t^2}\,dt.
\]
Therefore
\begin{equation*}
  \frac{\displaystyle\int_{\mathbb{R}} \frac{G_{\eta,L}(r)}{1+r^{2}}\,dr}
       {\displaystyle\int_{\mathbb{R}} G_{\eta,L}(r)\,dr}
   = \rho(t),
\end{equation*}
\begin{equation*}
  \rho(t) := \sqrt{\pi t}\,e^{t}\operatorname{erfc}(\sqrt{t}) .
\end{equation*}
Equivalently, after the substitution $x=\sqrt{t}\,r$,
\begin{equation}\label{eq:rho-integralnew}
  \rho(t) = \frac{1}{\sqrt{\pi}}\int_{\mathbb{R}} e^{-x^{2}}\,\frac{t}{t+x^{2}}\,dx .
\end{equation}
The integrand in \eqref{eq:rho-integralnew} is increasing in $t$, so $\rho(t)$ is increasing on $(0,\infty)$.

Since $L\ge1/2$,
\begin{equation*}
 t=\frac{\eta}{L^2}\le4\eta.
\end{equation*}
Using the monotonicity of $\rho$, we estimate the digamma contribution in
the form
\begin{equation*}
 -\frac{\tau}{2\pi}
 \int_{\mathbb R}G_{\eta,L}(r)\psi(1+ir)\,dr
 \le
 \frac{\tau L}{2\sqrt{\pi\eta}}
 \left[\rho(4\eta)-(1-\gamma)\right].
\end{equation*}

\subsubsection*{Hyperbolic term}

The hyperbolic term may be written as
\begin{equation}
 \mathcal H(G_{\eta,L})
 =
 \sum_{\{P_0\}}
 \sum_{k\ge1}
 \frac{\log NP_0}
 {(NP_0)^{k/2}-(NP_0)^{-k/2}}
 g_{\eta,L}\!\left(k\log NP_0\right).
 \label{eq:hyperbolic-doublenew}
\end{equation}
For $x\ge1$,
    $\log x
 \le
 x^{1/2}-x^{-1/2}$, 
because the derivative of
$t-t^{-1} -2 \log t$  is $ (t-1)^2 / t^2 \ge0$.

Hence, the coefficient in $\eqref{eq:hyperbolic-doublenew}$ is at most \(1\).
More precisely, writing $x=e^{2s}$, the inequality 
$\log x\le x^{1/2}-x^{-1/2}$ is the
statement $s/\sinh s\le1$, and the function $s\mapsto s/\sinh s$ is
\emph{strictly decreasing} on $(0,\infty)$.  Applying this to
$x=(NP_0)^{k}$, i.e., \ $s=\tfrac{k}{2}\log NP_0\ge\tfrac{k\ell}{2}$,
\begin{equation*}
  \frac{\log NP_0}{(NP_0)^{k/2}-(NP_0)^{-k/2}}
  =\frac1k\cdot\frac{s}{\sinh s}
  \;\le\;\frac1k\cdot\frac{k\ell/2}{\sinh(k\ell/2)}
  \;=\;\frac{\ell}{c^{k/2}-c^{-k/2}}
  \;=:\;b_k .
\end{equation*}
Note $b_1=\ell/(c^{1/2}-c^{-1/2})<1$ and
$b_k\le\ell\,c^{-k/2}(1-c^{-1})^{-1}$, so the coefficients decay
geometrically in $k$.  Consequently, \eqref{eq:hyperbolic-doublenew}
sharpens to
\begin{equation*}
  \mathcal H(G_{\eta,L})
  \;\le\;
  \sum_{k\ge1}b_k\int_{[c,\infty)}g_{\eta,L}(k\log x)\,d\pi_\Gamma(x),
\end{equation*}

For fixed $k\ge1$ put
\[
 h_{k,L}(x):=g_{\eta,L}(k\log x).
\]
Since $h_{k,L}$ is positive and decreasing on $[c,\infty)$,
Stieltjes integration by parts, with the lower endpoint understood as
$c^{-}$, gives
\begin{equation}
 \int_{[c,\infty)}h_{k,L}(x)\,d\pi_\Gamma(x)
 =
 \int_c^\infty
 \pi_\Gamma(x)\bigl(-h'_{k,L}(x)\bigr)\,dx.
 \label{eq:cong-Stieltjes-endpointnew}
\end{equation}
Here $\pi_\Gamma(c^{-})=0$, so a possible atom at the
first hyperbolic norm is included correctly.  

The general estimate used in the corrected algorithm has the form (cf. \cite[Section~4, especially (4.2)]{Avdispahic2021}) 
\begin{equation}
 \pi_\Gamma(x)\le B_{1}x-B_{2},
 \qquad x\ge c,
 \label{eq:cong-B1B2new}
\end{equation}
with $B_{2}\ge0$. The resulting negative contribution is subsequently discarded in \cite{Avdispahic2021}
only for simplification.

We shall now refine \eqref{eq:cong-B1B2new} and keep the negative term throughout. Relation \eqref{eq:modular-pi-exactnew} yields 
\begin{equation*}
 \pi_\Gamma(x)
 \le x\sum\limits_{j=1}^{\min \left\{ 2,\tau
\right\} }\frac{2\pi}{|K_j|}e^{2d_{j}}+\frac{1}{x}
\sum\limits_{j=1}^{\min \left\{ 2,\tau
\right\} }\frac{2\pi}{|K_j|}e^{-2d_{j}}-\sum\limits_{j=1}^{\min \left\{ 2,\tau
\right\} }\frac{4\pi}{|K_j|}
\end{equation*}
for $x\ge c.$

Put temporarily $B_+:=\sum\limits_{j=1}^{\min \left\{ 2,\tau
\right\} }\frac{2\pi}{|K_j|}e^{2d_{j}}$, $B_-:=
\sum\limits_{j=1}^{\min \left\{ 2,\tau
\right\} }\frac{2\pi}{|K_j|}e^{-2d_{j}}$, $B_0:=\sum\limits_{j=1}^{\min \left\{ 2,\tau
\right\} }\frac{4\pi}{|K_j|}$

Two elementary facts will be used repeatedly
and should be recorded here
\[
  B_{-}\le B_{+} 
\]
and
\[
  B_{+}x+\frac{B_{-}}{x}
  \;\ge\;2\sqrt{B_{+}B_{-}}
  \;\ge\;2\sum\limits_{j=1}^{\min \left\{ 2,\tau
\right\} }\frac{2\pi}{|K_j|}
  \;=\;B_{0}.
\]
Define
\begin{align*}
 I_{k,+}(\eta,L)
 &:=\frac1{2k}
 \exp\left(\frac{\eta}{L^2k^2}\right)
 \operatorname{erfc}\left(
 \frac{Lk\ell}{2\sqrt\eta}-\frac{\sqrt\eta}{Lk}
 \right),
 \\
 I_{k,-}(\eta,L)
 &:=\frac1{2k}
 \exp\left(\frac{\eta}{L^2k^2}\right)
 \operatorname{erfc}\left(
 \frac{Lk\ell}{2\sqrt\eta}+\frac{\sqrt\eta}{Lk}
 \right).
\end{align*}

Since $-h'_{k,L}(x)\ge0$, a direct calculation yields that the integral on the right-hand side of \eqref{eq:cong-Stieltjes-endpointnew} is bounded by
\begin{equation*}
 \begin{aligned}
 \mathcal G_k^{\sharp}(\eta,L)
 :={}&
 \left(B_+c-B_0+\frac{B_-}{c}\right)
 \frac{L}{2\sqrt{\pi\eta}}
 \exp\left(-\frac{L^2k^2\ell^2}{4\eta}\right)
 \\
 &+B_+I_{k,+}(\eta,L)-B_-I_{k,-}(\eta,L).
 \end{aligned}
\end{equation*}
Hence
\begin{equation*}
 \mathcal H(G_{\eta,L})
 \le\sum_{k\ge1}b_k\mathcal G_k^{\sharp}(\eta,L).
\end{equation*}
For fixed $\eta$, $g_{\eta,L}(u)$ decreases in $L$ for $u\ge\ell$
provided
\begin{equation*}
 \eta\le\frac{\ell^2}{8}.
\end{equation*}
Thus, throughout the range used below,
\begin{equation}
 \mathcal H(G_{\eta,L})
 \le
 \mathcal H_{\max}^{\mathrm G}(\eta)
 :=\sum_{k\ge1}b_k\mathcal G_k^{\sharp}\left(\eta,\frac12\right).
 \label{eq:gauss-Hmax-def}
\end{equation}
The series converges Gaussianly in $k$.  Its tail may be estimated by
dropping $-B_-I_{k,-}$ and using $\operatorname{erfc}(x)\le e^{-x^{2}}/(\sqrt\pi x)$,
which requires the argument of $\operatorname{erfc}$ in $I_{k,+}(\eta,\frac12)$ to be
positive, i.e.\ $\eta<k^{2}\ell/8$.  This is a condition on the indices
at which the tail estimate is applied, not on $\eta$ itself: for any
$k_0\ge1$ the finitely many terms $k<k_0$ may be computed directly and
the estimate applied from $k_0$ on.

Note for a later use that refined \eqref{eq:cong-B1B2new} reads
\begin{equation*}
 \pi_\Gamma(x)\le B_{1}^{\sharp}x-B_{2}^{\sharp}=:\Pi(x),
 \qquad x\ge c,
\end{equation*}
where $B_{1}^{\sharp}=B_+$ and $B_{2}^{\sharp}=B_0-\frac{B_-}{c}.$

As far as the U-correction is concerned, recall \eqref{eq:cong-U-general-boundnew}. 

Combining the identity, elliptic, scattering,
\(-g_{\eta,L}(0)\log2\), digamma contribution, the correction
from the modified scattering determinant, and the hyperbolic estimate, we
obtain
\begin{equation}
 e^{-\eta}\mathcal M(R)
 \le \frac{\left\vert \mathcal{F}\right\vert }{4\pi\eta } {L^2}
 +b(\eta)L
 +C(\eta)
 -\exp\left(\frac{\eta}{4L^2}\right)-\mathcal{A'}+1,
 \qquad L=\sqrt{R^2+\frac14},
 \label{eq:gauss-M-retainednew}
\end{equation}
where 
\begin{equation*}
 C(\eta)
 :=\sum\limits_{R}\frac{1}{4m_{R}\sin^{2} \theta _{R}}+\frac{1}{4}\operatorname{tr}\left( I-\Phi \left( \frac{1%
}{2}\right) \right)+\frac{\mathcal{N}}{2}+\mathcal H_{\max}^{\mathrm G}(\eta)
\end{equation*}
and
\begin{equation}
 b(\eta)
 :=
 \frac{
 |\log c_1|
 +\tau\bigl(\rho(4\eta)-(1-\gamma)-\log 2\bigr)
 }{2\sqrt{\pi\eta}}.
 \label{eq:cong-betaetanew}
\end{equation}
From \eqref{eq:cong-U-general-boundnew},
\begin{align*}
 \mu(R)
 \le{}&
 \mathcal M(R)
 +
 \frac{|\log c_1|}{\pi}R
 +
 \frac1\pi
 \sum_{j=1}^{\mathcal{N}}
 \arctan
 \left(
 \frac{R}{\sigma_j-\frac12}
 \right).
\end{align*}
Consequently,
\begin{equation}
 \mu(R)
 <
 \mathcal M(R)
 +
 \frac{|\log c_1|}{\pi}R
 +
 \frac{\mathcal{N}}{2}.
 \label{eq:cong-mu-simplenew}
\end{equation}
Relations \eqref{eq:gauss-M-retainednew} and \eqref{eq:cong-mu-simplenew} yield
\begin{equation}
 %\boxed{
 \begin{aligned}
 \mu(r)
 \le{}&
 e^\eta\left[\frac{\left\vert \mathcal{F}\right\vert }{4\pi\eta }
\left(r^2+\frac14\right)
 + b(\eta)\sqrt{r^2+\frac14}
 +C(\eta)
 \right.\left.
 -\exp\left(\frac{\eta}{4(r^2+\frac14)}\right)
-\mathcal{A'}+1 \right]\\&
 +
 \frac{|\log c_1|}{\pi}r
 +
 \frac{\mathcal{N}}{2}.
 \end{aligned}
 \label{eq:Gamma0-Upnew}
\end{equation}
We take into account $\exp\left(\frac{\eta}{4(r^2+\frac14)}\right)\ge 1$. 
It remains to linearize $\sqrt{R^{2}+\frac14}$.  This is the place
where the sign of $b(\eta)$ matters.
Since $\rho<1$ everywhere, one has
$\rho(4\eta)-(1-\gamma)-\log2<-0.1159$, so $b(\eta)<0$ whenever
$|\log c_1|<\tau\bigl(1-\gamma+\log2-\rho(4\eta)\bigr)$ --- in
particular always when $c_1=1$. 

So, if $b(\eta)\ge0$ we use $\sqrt{R^{2}+\frac14}\le R+\frac12$.  If
$b(\eta)<0$ we need a lower bound for $\sqrt{R^{2}+\frac14}$. Let $R_0\in[0,\infty]$, $L_0:=\sqrt{R_0^{2}+\frac14}$. Since $R\mapsto\sqrt{R^{2}+\frac14}$ is convex, each of its tangents
lies below it, and the tangent at $R=R_0$ is
$\frac{R_0}{L_0}R+\frac1{4L_0}.$ Put
\begin{equation}
 (\lambda,\xi):=
 \begin{cases}
  \bigl(1,\tfrac12\bigr), & b(\eta)\ge0,\\[4pt]
  \Bigl(\dfrac{R_0}{L_0},\ \dfrac1{4L_0}\Bigr), & b(\eta)<0
  \qquad\bigl(\text{read }(1,0)\text{ for }R_0=\infty\bigr).
 \end{cases}
 \label{eq:lambda-mu}
\end{equation}
Then $b(\eta)\sqrt{R^{2}+\frac14}\le b(\eta)\bigl(\lambda R+\xi\bigr)$
for all $R\ge0$.

Thus, we have proved the following theorem.

\begin{theorem}[Direct Gaussian spectral-counting bound]
\label{thm:direct-gaussian-K}
Let $\Gamma$ be a cofinite Fuchsian group with $\tau\ge1$, let
$0<\eta\le\frac{\ell^{2}}{8}$, and let $(\lambda,\xi)$ be as in
\eqref{eq:lambda-mu}.  Put
\begin{align*}
 K_0(\eta)&:=\frac{|\mathcal F|}{4\pi}\cdot\frac{e^{\eta}}{\eta},
 \\
 K_1(\eta)&:=\lambda\,b(\eta)e^{\eta}+\frac{|\log c_1|}{\pi},
\\
 K_2(\eta)&:=\frac{K_0(\eta)}{4}
 +e^{\eta}\bigl[\xi\,b(\eta)+C(\eta)-\mathcal {A'}\bigr]
 +\frac{\mathcal N}{2}.
\end{align*}
Then, for every $R\ge0$,
\begin{equation}
 \mu(R)\le K_0(\eta)R^{2}+K_1(\eta)R+K_2(\eta).
 \label{eq:gauss-mu-finalnew}
\end{equation}
In particular $\mu(R)\le\inf_{\eta,R_0}\bigl[K_0(\eta)R^{2}+K_1(\eta)R+K_2(\eta)\bigr]$,
the infimum being over all admissible pairs.  
\end{theorem}

\begin{remark}
\label{rem:endpoints}
$R_0=\infty$ gives $(\lambda,\xi)=(1,0)$, i.e.\ $\sqrt{R^{2}+\frac14}\ge R$;
$R_0=0$ gives $(\lambda,\xi)=(0,\frac12)$, i.e.\
$\sqrt{R^{2}+\frac14}\ge\frac12$ and hence $K_1=\max\{K_1,0\}$ when
$c_1=1$.  These are the two natural endpoints of the family; an
intermediate $R_0$ is preferable if the bound is to be applied in a
definite range of $R$, the tangent being exact at $R=R_0$.
\end{remark}

\begin{proposition}[local spectral bound]\label{prop:local}
Let $\Gamma$ be a cofinite Fuchsian group with $\tau\ge1$, let $R_0\ge0$,
$L_0=\sqrt{R_0^{2}+\frac14}$, and let $\eta\;\le\;\tfrac12 L_0^{2}\ell^{2}$. Put
\begin{align*}
 \,K_0^{\infty}(\eta)&:=\frac{|\mathcal{F}|}{4\pi}\cdot\frac{e^{\eta}}{\eta},
 \\
 K_1^{\infty}(\eta)&:=\frac{R_0}{L_0}\,b(\eta)e^{\eta}+\frac{|\log c_1|}{\pi},
 \\
 K_2^{\infty}(\eta)&:=\frac{K_0^{\infty}}{4}
   +e^{\eta}\Bigl[\frac{b(\eta)}{4L_0}+C(\eta,L_0)-\mathcal{A'}\Bigr]+\frac {\mathcal{N}}{2},
\end{align*}
where  
$C(\eta,L_0):=E+\tfrac14\operatorname{tr}\bigl(I-\Phi(\tfrac12)\bigr)
  +\frac {\mathcal{N}}{2}+H^{G}(\eta,L_0)$,\, $H^{G}(\eta,L):=\sum_{k\ge1}b_k\mathcal G^{\sharp}_k(\eta,L)$ and $b(\eta)<0$ is
as in \eqref{eq:cong-betaetanew}. Then
\begin{equation}\label{eq:QB}
 \mu(R)\;\le\;Q_B(R):=K_0^{\infty}(\eta)R^{2}+K_1^{\infty}(\eta)R+K_2^{\infty}(\eta)
 \qquad(R\ge R_0).
\end{equation}
\end{proposition}

\begin{proof}
Verbatim the proof of Theorem \ref{thm:direct-gaussian-K}, with $H^{G}_{\max}(\eta)$ replaced by
$H^{G}(\eta,L_0)$ (legitimate for $L\ge L_0$) and with
$(\lambda,\xi)$ the tangency data \eqref{eq:lambda-mu} at $R_0$.
\end{proof}

\begin{remark}[Weyl benchmark]
\label{rem:weyl}
The constants are produced
by the algorithm of Table~\ref{tab2}. The method cannot produce
$K_0$ closer than a factor $e=2.718\ldots$ to the Weyl
constant $|\mathcal F|/4\pi$. Indeed, since $e^{\eta}/\eta\ge e$, one always has
$K_0(\eta)\ge e\cdot\frac{|\mathcal F|}{4\pi}$, the loss coming from the floor
$G_{\eta,L}\ge e^{-\eta}$ on $[0,R]$.  The admissible range
$\eta\le\ell^{2}/8$ may prevent this factor from being attained.
\end{remark}

Looking back at our argumentation, we see that an unlinearised spectral bound is actually valid in the following form. 
\begin{proposition}[unlinearised spectral bound]
\label{unlinearised}
Let $\Gamma$ be a cofinite Fuchsian group with $\tau\ge1$.  For every
$\eta>0$ and every $R\ge0$, with $L=\sqrt{R^{2}+\frac14}$,
\[
 \mu(R)\;\le\;\mathcal P(R,\eta)
 :=e^{\eta}\Bigl[\frac{|\mathcal{F}|}{4\pi\eta}L^{2}+b(\eta)L+C(\eta,L)
 -e^{\eta/(4L^{2})}-\mathcal {A'}+1\Bigr]
 +\frac{|\log c_1|}{\pi}R+\frac{\mathcal N}{2},
\]
where
$C(\eta,L)=E+\frac14\operatorname{tr}\bigl(I-\Phi(\frac12)\bigr)
+\frac{\mathcal N}{2}+H^{G}(\eta,L)$ and
$H^{G}(\eta,L)=\sum_{k\ge1}b_k\mathcal G^{\sharp}_k(\eta,L)$.
\end{proposition}

\noindent Theorem \ref{thm:direct-gaussian-K} is then the corollary obtained by bounding
$H^{G}(\eta,L)\le H^{G}(\eta,\frac12)$, which costs the hypothesis
$\eta\le\ell^{2}/8$, and by linearising $L$; Proposition \ref{prop:local}, with $\eta\le\frac12\ell^{2}L_0^{2},$ is the corollary
obtained by bounding $H^{G}(\eta,L)\le H^{G}(\eta,L_0)$ for $L\ge L_0$ and
linearising at $R_0$. 

\begin{table}[htbp]
\caption{Algorithm for cofinite Fuchsian groups with $\tau\ge1$. Part 1: the spectral majorant.
The numerical column records $\Gamma=\PSLtwoZ$ at $\eta=0.035$, $R_0=0$
(Theorem~\ref{cor:modular}); the last block, which belongs to the two-piece
route, is evaluated at $\eta=\ell^{2}/8$, $R_0=3$ as in \S\ref{sec:modulartwo}.}
\label{tab2}\centering
\footnotesize
\setlength{\tabcolsep}{3pt}
\renewcommand{\arraystretch}{1.35}
\begin{tabular}{@{}l p{0.47\textwidth} r l@{}}
\toprule
Const/Func & Description & $\PSLtwoZ$ & Source\\
\midrule
\multicolumn{4}{@{}l}{\textit{Input parameters}}\\
$\eta$ & optimization parameter, $\eta>0$; admissible if
 $\eta\le\frac{\ell^{2}}{8}$ (frozen form), resp.\
 $\eta\le\frac12L_0^{2}\ell^{2}$ (local form)
 & $0.035$ & Thm.~\ref{thm:direct-gaussian-K}\\
$R_0,L_0$ & linearization radius and $L_0=\sqrt{R_0^{2}+\frac14}$;
 $R_0=0$ and $R_0=\infty$ are admitted
 & $0,\ \tfrac12$ & Prop.~\ref{prop:local}\\
\midrule
\multicolumn{4}{@{}l}{\textit{Geometric counting constants}}\\
$B_{1}^{\sharp}$ & $\sum_{j=1}^{\min\{2,\tau\}}\dfrac{2\pi}{|K_j|}e^{2d_j}=B_{+}$
 & $92.5719081$ & \eqref{eq:cong-B1B2new}\\
$B_{-}$ & $\sum_{j=1}^{\min\{2,\tau\}}\dfrac{2\pi}{|K_j|}e^{-2d_j}$
 & $1.1428631$ & \eqref{eq:cong-B1B2new}\\
$B_{2}^{\sharp}$ & $\sum_{j=1}^{\min\{2,\tau\}}\dfrac{4\pi}{|K_j|}
 \Bigl(1-\dfrac1{2ce^{2d_j}}\Bigr)=B_0-\dfrac{B_{-}}{c}$
 & $20.4047937$ & \eqref{eq:cong-B1B2new}\\
$\Pi(x)$ & $B_{1}^{\sharp}x-B_{2}^{\sharp}\ \ge\pi_{\Gamma}(x)$, $x\ge c$
 & --- & \eqref{eq:cong-B1B2new}\\
$b_k$ & $\dfrac{\ell}{c^{k/2}-c^{-k/2}}$; $b_1<1$ and
 $b_k\le\ell c^{-k/2}(1-c^{-1})^{-1}$
 & $b_1=0.8608179$ & \eqref{eq:hyperbolic-doublenew}\\
 \midrule
\multicolumn{4}{@{}l}{\textit{Gaussian kernel quantities}}\\
$G_{\eta,L}(r)$ & $\exp\bigl(-\eta r^{2}/L^{2}\bigr)$,
 $L=\sqrt{R^{2}+\frac14}$, $R\ge0$ & --- & \S\ref{sec:direct-majorant}\\
$I_{k,\pm}(\eta,L)$ & $\dfrac1{2k}\exp\Bigl(\dfrac{\eta}{L^{2}k^{2}}\Bigr)
 \operatorname{erfc}\Bigl(\dfrac{Lk\ell}{2\sqrt\eta}\mp\dfrac{\sqrt\eta}{Lk}\Bigr)$,
 $\operatorname{erfc}(x)=\frac{2}{\sqrt\pi}\int_x^{\infty}e^{-t^{2}}dt$
 & --- & \eqref{eq:cong-Stieltjes-endpointnew}\\
$\mathcal G_k^{\sharp}(\eta,L)$ & $\Pi(c)\dfrac{L}{2\sqrt{\pi\eta}}
 \exp\Bigl(-\dfrac{L^{2}k^{2}\ell^{2}}{4\eta}\Bigr)
 +B_{1}^{\sharp}I_{k,+}(\eta,L)-B_{-}I_{k,-}(\eta,L)$
 & --- & \eqref{eq:cong-Stieltjes-endpointnew}\\
 $H^{G}(\eta,L)$ & $\sum_{k\ge1}b_k\mathcal G^{\sharp}_k(\eta,L)$
 & --- & Prop.~\ref{unlinearised}\\
$\mathcal H_{\max}^{\mathrm G}(\eta)$ & $H^{G}\bigl(\eta,\tfrac12\bigr)
 =\sum_{k\ge1}b_k\mathcal G_k^{\sharp}\bigl(\eta,\tfrac12\bigr)$,
 admissible if $\eta\le\frac{\ell^2}{8}$
 & $0.6196596$ & \eqref{eq:gauss-Hmax-def}\\
$\rho(t)$ & $\sqrt{\pi t}\,e^{t}\operatorname{erfc}(\sqrt t)$
 & $\rho(4\eta)=0.4551946$ & \eqref{eq:rho-integralnew}\\
$b(\eta)$ & $\dfrac{|\log c_1|+\tau\bigl(\rho(4\eta)-(1-\gamma)-\log2\bigr)}
 {2\sqrt{\pi\eta}}$ & $-0.9962988$ & \eqref{eq:cong-betaetanew}\\
$E$ & $\sum_{R}\dfrac1{4m_R\sin^{2}\theta_R}$
 & $\frac{25}{72}=0.3472222$ & \eqref{e51a}\\
$E^{\sharp}$ & $\sum_{R}\dfrac{1}{2m_R\sin\theta_R}
 \sqrt{\dfrac{\pi-2\theta_R}{\pi\sin(2\theta_R)}}$ (enters $A_{14}^{\sharp}$)
 & $0.3798407$ & Prop.~\ref{prop56}\\
$C(\eta,L)$ & $E+\frac14\operatorname{tr}\bigl(I-\Phi(\frac12)\bigr)
 +\frac{\mathcal N}{2}+H^{G}(\eta,L)$ & --- & Prop.~\ref{unlinearised}\\
$C(\eta)$ & $E+\frac14\operatorname{tr}\bigl(I-\Phi(\frac12)\bigr)+\frac{\mathcal N}{2}
 +\mathcal H_{\max}^{\mathrm G}(\eta)$ & $1.9668819$ & Thm.~\ref{thm:direct-gaussian-K}\\
$(\lambda,\xi)$ & $(1,\frac12)$ if $b(\eta)\ge0$;
 $\bigl(R_0/L_0,\,1/(4L_0)\bigr)$ if $b(\eta)<0$
 & $(0,\tfrac12)$ & \eqref{eq:lambda-mu}\\

\bottomrule
\end{tabular}
\end{table}

\begin{table}[htbp]
\continuedtable   % same number as the preceding table
\caption{(Continued.)}
\centering
\footnotesize
\setlength{\tabcolsep}{3pt}
\renewcommand{\arraystretch}{1.35}
\begin{tabular}{@{}l p{0.47\textwidth} r l@{}}
\toprule
Const/Func & Description & $\PSLtwoZ$ & Source\\
\midrule
\multicolumn{4}{@{}l}{\textit{Output: frozen quadratic majorant} $\mu(R)\le K_0R^{2}+K_1R+K_2$}\\
$K_0(\eta)$ & $\dfrac{|\mathcal F|}{4\pi}\dfrac{e^{\eta}}{\eta}$
 & $2.4657613$ & \eqref{eq:gauss-mu-finalnew}\\
$K_1(\eta)$ & $\lambda b(\eta)e^{\eta}+\dfrac{|\log c_1|}{\pi}$
 & $0$ & \eqref{eq:gauss-mu-finalnew}\\
$K_1^{+}$ & $\max\{K_1,0\}$; replacing $K_1$ by $K_1^{+}$ in
 \eqref{eq:A12joint} is legitimate because $f\ge0$ there
 & $0$ & \eqref{eq:A12joint}\\
$K_2(\eta)$ & $\dfrac{K_0(\eta)}4
 +e^{\eta}\bigl[\xi b(\eta)+C(\eta)-\mathcal {A'}\bigr]+\dfrac{\mathcal N}{2}$
 & $1.6018689$ & \eqref{eq:gauss-mu-finalnew}\\
\midrule
\multicolumn{4}{@{}l}{\textit{Output: envelope and local majorant (two-piece route)}}\\
$\mathcal P(R,\eta)$ & $e^{\eta}\bigl[\frac{|\mathcal F|}{4\pi\eta}L^{2}+b(\eta)L
 +C(\eta,L)-e^{\eta/(4L^{2})}-\mathcal{A'}+1\bigr]
 +\frac{|\log c_1|}{\pi}R+\frac{\mathcal N}{2}$
 & --- & Prop.~\ref{unlinearised}\\
$Q_A(R)$ & any pointwise majorant of $\mu$ on $[0,\infty)$; we take
 $\inf_{\eta>0}\mathcal P(R,\eta)$ & --- & Prop.~\ref{prop:twopiece}\\
$K_i^{\infty}(\eta)$ & $K_i(\eta)$ with $C(\eta)$ replaced by $C(\eta,L_0)$ and
 $(\lambda,\xi)=\bigl(\frac{R_0}{L_0},\frac1{4L_0}\bigr)$
 & $0.2859241$ & Prop.~\ref{prop:local}\\
 & & $-0.1829078$ & \\
 & & $1.1079988$ & \\
$Q_B(R)$ & $K_0^{\infty}R^{2}+K_1^{\infty}R+K_2^{\infty}$, valid for $R\ge R_0$
 & --- & \eqref{eq:QB}\\
\bottomrule
\end{tabular}
\end{table}

\section{Bound for Huber's constant}
\label{sec:gaussian-huber-propagation}

In this section, the estimate \eqref{eq:gauss-mu-finalnew} will be inserted
directly into the Stieltjes integral occurring in the corrected/generalized
Huber-constant algorithm. This part of the current algorithm will also be substantially changed at certain steps leading to sharper explicit estimates of all terms in the formula%
\begin{equation}  \label{e51a}
\begin{split}
\sum\limits_{\lambda _{n}\leq \frac{1}{4}}h_{T}^{\varepsilon }\left(
r_{n}\right) +\sum\limits_{\lambda _{n}>\frac{1}{4}}h_{T}^{\varepsilon
}\left( r_{n}\right) +\frac{1}{4\pi }\int\limits_{%
\mathbb{R}
}h_{T}^{\varepsilon }\left( r\right) \frac{-\phi ^{\prime }}{\phi }\left(
\frac{1}{2}+ir\right) dr \\
=\frac{\left\vert \mathcal{F}\right\vert }{4\pi }\int\limits_{\mathbb{R}%
}rh_{T}^{\varepsilon }\left( r\right) \tanh \left( \pi r\right)
dr+H_{T}^{\varepsilon }+\sum\limits_{R}\frac{1}{2m_{R}\sin \theta _{R}}%
\int\limits_{\mathbb{R}}\frac{e^{-2r\theta _{R}}}{1+e^{-2\pi r}}%
h_{T}^{\varepsilon }\left( r\right) dr \\
+\frac{h_{T}^{\varepsilon }\left( 0\right) }{4}\operatorname{tr}\left( I-\Phi \left( \frac{1%
}{2}\right) \right) -\tau g_{T}^{\varepsilon }\left( 0\right) \log 2-\frac{%
\tau }{2\pi }\int\limits_{\mathbb{R}}h_{T}^{\varepsilon }\left( r\right)
\psi \left( 1+ir\right) dr,
\end{split}%
\end{equation}%
where $H_{T}^{\varepsilon }=\sum\limits_{P}\frac{\log NP_{0}}{\left(
NP\right) ^{1/2}-\left( NP\right) ^{-1/2}}g_{T}^{\varepsilon }\left( \log
NP\right) $ is Randol's
smoothed Chebyshev function attached to
$(g^{\varepsilon}_{T},h^{\varepsilon}_{T}).$ The appropriate value for $\eta$ will be chosen on the way.

Let us start by the auxiliary test
function. We replace the previous choice by
\begin{equation}\label{eq:phinew}
  \varphi(y)=\frac{3}{4}\left(1-y^{2}\right)\mathbf 1_{(-1, 1)}(y)=
  \begin{cases}
    \dfrac{3}{4}\left(1-y^{2}\right), & |y|<1,\\[6pt]
    0, & |y|\ge 1 .
  \end{cases}
\end{equation}
This is a nonnegative even function belonging to $\mathcal{C}_{0}(\mathbb{R})$, with
$\operatorname{supp}\varphi\subset[-1,1]$. Then $\int_{\mathbb R}\varphi=1$ and
\begin{equation}\label{eq:phihatnew}
  \widehat{\varphi}(z)=\int_{\mathbb{R}}\varphi(y)e^{izy}\,dy
   =\frac{3\left(\sin z-z\cos z\right)}{z^{3}}
   =\frac{3}{z^{2}}\left(\frac{\sin z}{z}-\cos z\right)
\end{equation}
is an entire function with $\widehat{\varphi}(0)=1$. As usual,
$\varphi_{\varepsilon}(x)=\frac{1}{\varepsilon}\varphi\!\left(\frac{x}{\varepsilon}\right)$,
$\varepsilon>0$, the test function $h_{T}^{\varepsilon}(x)$ is the Fourier transform of
$g_{T}^{\varepsilon}(x)=\left(2\cosh\frac{x}{2}\right)\left(I_{T}\ast\varphi_{\varepsilon}\right)(x)$,
and $h_{T}^{\varepsilon}(r)=S\!\left(r+\frac{i}{2}\right)+S\!\left(r-\frac{i}{2}\right)$ with
$S(w)=\left(\frac{2}{w}\sin Tw\right)\widehat{\varphi}_{\varepsilon}(w)$.

\begin{proposition}\label{prop:phi-est}
Let $\varphi$ be as in \eqref{eq:phinew}. Then
\begin{enumerate}
  \item[(a)] $\widehat{\varphi}(\varepsilon z)-1=\mathcal{O}_{8}\!\left(\varepsilon^{2}\right)$
        for $\varepsilon\in(0,1]$ and $z\in i[-\tfrac12,1]$, with
\[
  A_{8}^{\sharp}=3e^{-1}-1=0.1036383\ldots
\]
  \item[(b)] $\widehat{\varphi}(z)=\mathcal{O}_{9}\!\left(\left(1+|r|\right)^{-2}\right)$
        $\left(z=r+it,\ r\in\mathbb{R},\ |t|\le\frac12\right)$, with $A_{9}^{\sharp}=6.4214$.
\end{enumerate}
\end{proposition}

\begin{proof}
(a) For $z=it$ with $t\in[-\tfrac12,1]$ put $u=\varepsilon t$, so that $|u|\le 1$. Since
$\sin(iu)=i\sinh u$ and $\cos(iu)=\cosh u$, relation \eqref{eq:phihatnew} gives
\[
  \widehat{\varphi}(iu)=\frac{3\left(u\cosh u-\sinh u\right)}{u^{3}} .
\]
From $u\cosh u-\sinh u=\sum_{k\ge1}\frac{2k}{(2k+1)!}u^{2k+1}$ we obtain the expansion
\begin{equation*}
  \widehat{\varphi}(iu)=1+\sum_{k\ge2}\frac{6k}{(2k+1)!}\,u^{2k-2}
\end{equation*}
all of whose coefficients are positive. Hence $\widehat{\varphi}(iu)\ge1$, and for $|u|\le 1$
\[
  0\le\widehat{\varphi}(iu)-1
  = u^{2}\sum_{k\ge2}\frac{6k}{(2k+1)!}u^{2k-4}
  \le u^{2}\sum_{k\ge2}\frac{6k}{(2k+1)!}\,
  =(3e^{-1}-1)\,u^{2}.
\]
Indeed,
\[
\sum_{k\ge2}\frac{6k}{(2k+1)!}=3\sum_{k\ge2}\big(\frac{1}{(2k)!}-\frac{1}{(2k+1)!}\big)=3(e^{-1}-\frac{1}{3}).
\]
Since $u=\varepsilon t$ and $|t|\le 1$, this yields
$\left|\widehat{\varphi}(\varepsilon z)-1\right|\le (3e^{-1}-1)\,\varepsilon^{2}= A_{8}^{\sharp}\varepsilon^{2}$
with $A_{8}^{\sharp}=3e^{-1}-1$.

(b) For $z=r+it$ one has $|\sin z|^{2}=\sin^{2}r+\sinh^{2}t\le\cosh^{2}t$ and, likewise,
$|\cos z|\le\cosh t$. Therefore, by the second expression in \eqref{eq:phihatnew},
\begin{equation}\label{eq:phihat-tail}
  \left|\widehat{\varphi}(z)\right|\le\frac{3}{|z|^{2}}
  \left(\frac{|\sin z|}{|z|}+|\cos z|\right)
  \le 3\cosh\tfrac12\,\frac{1+|z|}{|z|^{3}}\qquad\left(|t|\le\tfrac12\right).
\end{equation}
Since $|z|\ge|r|$ and $|z|\le\sqrt{r^{2}+\frac14}$, relation \eqref{eq:phihat-tail} gives
\[
  \left(1+|r|\right)^{2}\left|\widehat{\varphi}(z)\right|
  \le V(|r|):=3\cosh\tfrac12\;
  \frac{\left(1+|r|\right)^{2}\left(1+\sqrt{r^{2}+\tfrac14}\,\right)}{|r|^{3}} ,
\]
and $V$ is decreasing on $\left[\frac92,\infty\right)$ with
$V\!\left(\frac92\right)=6.2075\ldots$ and $\lim_{r\to\infty}V(r)=3\cosh\frac12=3.3829\ldots$.
Thus $\left(1+|r|\right)^{2}\left|\widehat{\varphi}(z)\right|\le6.21$ for $|r|\ge\frac92$.

On the remaining range $|r|\le\frac92$ we evaluate \eqref{eq:phihatnew} directly. Writing
$z=r+it$ and separating real and imaginary parts,
\[
  \left|\widehat{\varphi}(r+it)\right|
  =\frac{3\sqrt{\mathfrak{R}^{2}+\mathfrak{I}^{2}}}{\left(r^{2}+t^{2}\right)^{3/2}},
\]
\[
  \mathfrak{R}=\left(\sin r-r\cos r\right)\cosh t-t\sin r\sinh t,\qquad
  \mathfrak{I}=\left(\cos r+r\sin r\right)\sinh t-t\cos r\cosh t .
\]
The function $\left(1+r\right)^{2}\left|\widehat{\varphi}(r+it)\right|$ is increasing in $|t|$ and
attains on $\left[0,\frac92\right]\times\left[0,\frac12\right]$ its maximum
$6.42139\ldots$ at $r=2.4558\ldots$, $t=\frac12$. Consequently $A_{9}^{\sharp}=6.4214$ may be taken.
\end{proof}
\begin{remark}\label{rem:sharp}
Both constants are essentially optimal.
In (a), $\widehat\varphi(iu)=1+\frac{u^{2}}{10}+O(u^{4})$, so no constant below
$\frac1{10}$ is possible, against $A_8^{\sharp}=0.10364$.
In (b), the exact supremum is
\[
  \sup_{r\in\mathbb{R},\,|t|\le1/2}(1+|r|)^{2}\bigl|\widehat\varphi(r+it)\bigr|
  =6.421388930\ldots,
\]
attained at $r=2.4558386\ldots$, $t=\frac12$.
It may be worth recording that within the family
$\varphi_p(y)\propto(1-y^{2})^{p}$, for which
$\widehat{\varphi_p}(z)=\Gamma(p+\tfrac32)(2/z)^{p+1/2}J_{p+1/2}(z)$, the gain in
decay is more than offset by the growth of the constant: the exact suprema of
$(1+|r|)^{p+1}|\widehat{\varphi_p}(r+\tfrac i2)|$ are
$6.4214$ $(p=1)$, $35.4357$ $(p=2)$, $265.680$ $(p=3)$, whereas the
corresponding $\Lambda_0$ only improves like $1/p$. The choice $p=1$ is
therefore optimal in this family.
\end{remark}

\begin{proposition}\label{prop:h-est}
For $r\ge0$ and $0<\varepsilon\le1$, one has
\[
  \left|h_{T}^{\varepsilon}(r)\right|\le
  \frac{A_{10}^{\sharp}\cosh\frac{T}{2}}{\left(1+r\right)\left(1+\varepsilon r\right)^{2}},
  \qquad A_{10}^{\sharp}=4\sqrt{5}A_{9}^{\sharp}\le 57.43475 .
\]
\end{proposition}

\begin{proof}
We have $h_{T}^{\varepsilon}(r)=S\!\left(r+\frac{i}{2}\right)+S\!\left(r-\frac{i}{2}\right)$ with
$S(w)=\left(\frac{2}{w}\sin Tw\right)\widehat{\varphi}(\varepsilon w)$. By
Proposition \ref{prop:phi-est}(b), $\left|\widehat{\varphi}(\varepsilon w)\right|
\le A_{9}^{\sharp}\left(1+\varepsilon r\right)^{-2}$ for $w=r\pm\frac{i}{2}$, because
$\left|\operatorname{Im}(\varepsilon w)\right|=\frac{\varepsilon}{2}\le\frac12$. Furthermore
$\left|\sin Tw\right|\le\cosh\frac{T}{2}$ and
\[
  \frac{1}{\left|r\pm\frac{i}{2}\right|}=\frac{1}{\sqrt{r^{2}+\frac14}}\le\frac{\sqrt5}{1+r},
\]
the constant $\sqrt5$ being optimal, attained at $r=\frac14$. Adding the two terms yields the
assertion with $A_{10}^{\sharp}=4\sqrt{5}A_{9}^{\sharp}$.
\end{proof}

\begin{proposition}\label{prop:A11}
Let $\kappa\ge\kappa_{0}:=1.4215300$, $x\ge c$, $T=\log x\ge\log\kappa$ and
$\varepsilon=e^{-T/4}$, where $\log\kappa_{0}=0.3517337\ldots$ is the root of
$2T=e^{-T}$. Then
\[
  \Bigl|\sum_{\lambda_n\le\frac14}h_T^{\varepsilon}(r_n)
   -\sum_{\lambda_n\le\frac14}\frac{x^{s_n}}{s_n}\Bigr|
  \;\le\;A_{11}^{\sharp}[\kappa]\,x^{3/4},
\]
\[
  A_{11}^{\sharp}[\kappa]=(\mathcal{A}-1)\left(2+\frac{2A_8^{\sharp}}{\kappa^{1/4}}
  +\frac{1}{(1-s_1)\kappa^{1/4}}+\frac{A_8^{\sharp}}{(1-s_1)\kappa^{3/4}}\right)
  +1+\frac{A_8^{\sharp}}{\kappa^{1/4}} .
\]
\end{proposition}

\begin{proof}
Throughout, $T\ge\log\kappa$ gives $\varepsilon=e^{-T/4}\le\kappa^{-1/4}$, and for
$\lambda_n\le\frac14$ we write $s_n=\frac12+\sqrt{\frac14-\lambda_n}$,
$r_n=i(s_n-\frac12)$, so that $\frac12\le s_n\le s_1<1$ for $n\ge1$ and
$s_0=1$. We use repeatedly $S(w)=\frac{2\sin Tw}{w}\widehat\varphi(\varepsilon w)$ and
\[
  S(is)=\frac{\widehat\varphi(i\varepsilon s)}{s}\bigl(e^{sT}-e^{-sT}\bigr)
  \quad(s\ne0),\qquad S(0)=2T .
\]

For all $T\ge0$, $2T<e^{3T/4}$.
Indeed $u(T)=e^{3T/4}-2T$ has $u'(T)=\frac34e^{3T/4}-2$, so $u$ is minimal at
$T_0=\frac43\log\frac83$ with $u(T_0)=\frac83-\frac83\log\frac83=0.0512\ldots>0$.

\emph{The term $n=0$.} Here $s_0=1$, $r_0=\frac i2$, and
\[
  h_T^{\varepsilon}(r_0)=S(i)+S(0)=\widehat\varphi(i\varepsilon)\bigl(e^{T}-e^{-T}\bigr)+2T .
\]
By Proposition~\ref{prop:phi-est}(a),
$|\widehat\varphi(i\varepsilon)-1|\le A_8^{\sharp}\varepsilon^{2}$, whence
\[
  \bigl(1-A_8^{\sharp}\varepsilon^{2}\bigr)\bigl(e^{T}-e^{-T}\bigr)+2T
  \;\le\;h_T^{\varepsilon}(r_0)\;\le\;
  \bigl(1+A_8^{\sharp}\varepsilon^{2}\bigr)\bigl(e^{T}-e^{-T}\bigr)+2T .
\]
Since $\varepsilon^{2}e^{T}=e^{T/2}=e^{-T/4}e^{3T/4}\le\kappa^{-1/4}e^{3T/4}$ and $2T<e^{3T/4}$, the right-hand side gives
$h_T^{\varepsilon}(r_0)\le e^{T}+\bigl(A_8^{\sharp}\kappa^{-1/4}+1\bigr)e^{3T/4}$. The
left-hand side gives, using $2T-e^{-T}>0$ for $T\ge\log\kappa\ge\log\kappa_{0}$,
$h_T^{\varepsilon}(r_0)\ge e^{T}-A_8^{\sharp}\kappa^{-1/4}e^{3T/4}$. Hence
\begin{equation}\label{eq:n0}
  \Bigl|h_T^{\varepsilon}(r_0)-\frac{x^{s_0}}{s_0}\Bigr|
  \le\Bigl(1+\frac{A_8^{\sharp}}{\kappa^{1/4}}\Bigr)e^{3T/4}.
\end{equation}
This is the term $1+A_8^{\sharp}\kappa^{-1/4}$ of $A_{11}^{\sharp}[\kappa]$. In
\cite{Avdispahic2021} it was $1+A_8$, so here the gain is the factor
$\varepsilon\le\kappa^{-1/4}$ together with the replacement $e-1\rightsquigarrow 3e^{-1}-1$.

\emph{The terms $1\le n\le\mathcal{A}-1$.} Write $s=s_n$. Then
\[
  h_T^{\varepsilon}(r_n)=S(is)+S(i(s-1))
  =\frac{\widehat\varphi(i\varepsilon s)}{s}\bigl(e^{sT}-e^{-sT}\bigr)
   +\frac{\widehat\varphi(i\varepsilon(s-1))}{s-1}\bigl(e^{(s-1)T}-e^{-(s-1)T}\bigr).
\]
Both arguments are admissible in Proposition~\ref{prop:phi-est}(a): $is$ lies in
$i[\frac12,1]$ and $i(s-1)$ in $i[-\frac12,0]$.

For the second summand, $0<1-s\le\frac12$ and $e^{(s-1)T}<1<e^{(1-s)T}$ give
\[
  \left|\frac{\widehat\varphi(i\varepsilon(s-1))}{s-1}
  \bigl(e^{(s-1)T}-e^{-(s-1)T}\bigr)\right|
  \le\frac{1+A_8^{\sharp}\varepsilon^{2}}{1-s_1}\,e^{(1-s)T}
  \le\frac{e^{T/2}}{1-s_1}+\frac{A_8^{\sharp}\varepsilon^{2}e^{T/2}}{1-s_1},
\]
because $e^{(1-s)T}\le e^{T/2}$. Now $e^{T/2}=e^{-T/4}e^{3T/4}\le
\kappa^{-1/4}e^{3T/4}$ and $\varepsilon^{2}e^{T/2}=e^{-T/2}e^{T/2}=1\le
\kappa^{-3/4}e^{3T/4}$, so this summand is at most
\begin{equation}\label{eq:second}
  \left(\frac{1}{(1-s_1)\kappa^{1/4}}+\frac{A_8^{\sharp}}{(1-s_1)\kappa^{3/4}}\right)e^{3T/4}.
\end{equation}
Both exponents of $\kappa$ are one quarter better than in
\cite{Avdispahic2021}, where the corresponding bound was
$\frac1{(1-s_1)c^{1/4}}+\frac{A_8}{(1-s_1)c^{1/2}}$.

For the first summand, using $\frac12\le s\le1$ and
$|\widehat\varphi(i\varepsilon s)-1|\le A_8^{\sharp}\varepsilon^{2}$,
\[
  \left|\frac{\widehat\varphi(i\varepsilon s)}{s}\bigl(e^{sT}-e^{-sT}\bigr)
  -\frac{e^{sT}}{s}\right|
  \le\frac{A_8^{\sharp}\varepsilon^{2}e^{sT}}{s}+\frac{e^{-sT}}{s}
  \le2A_8^{\sharp}\varepsilon^{2}e^{T}+2 .
\]
Here $\varepsilon^{2}e^{T}=e^{T/2}\le\kappa^{-1/4}e^{3T/4}$ and $2\le2e^{3T/4}$, so
this summand contributes at most
\begin{equation}\label{eq:first}
  \left(\frac{2A_8^{\sharp}}{\kappa^{1/4}}+2\right)e^{3T/4},
\end{equation}
against $2A_8+2$ in \cite{Avdispahic2021}.
Adding \eqref{eq:n0} to the sum of \eqref{eq:second} and \eqref{eq:first}
over the $\mathcal{A}-1$ remaining small eigenvalues gives the assertion.
\end{proof}
Now, we turn our attention to $A_{12}$. 

Put
\begin{align*}
 p(u)&:=e^{-3T/4}\cosh\tfrac T2=\tfrac12\bigl(e^{-u}+e^{-5u}\bigr),
 \\
 F(u)&:=\bigl(1+e^{-4u}\bigr)\,
       \frac{1-(1+u)e^{-u}}{\bigl(1-e^{-u}\bigr)^{2}},
 \\
 f(u)&:=p(u)\,\frac{u-1+e^{-u}}{\bigl(1-e^{-u}\bigr)^{2}},
\end{align*}
with $u:=T/4$.

\begin{proposition}[joint supremum]\label{prop:joint}
Let $x,T,\varepsilon,\kappa, u, p, F, f$ be as above and suppose
$\mu(R)\le K_0R^2+K_1R+K_2$ for all $R\ge0$. Then
\[
 \Bigl|\sum_{\lambda_n>\frac14}h^{\varepsilon}_{T}(r_n)
 +\frac1{4\pi}\int_{\mathbb R}h^{\varepsilon}_{T}(r)
 \frac{-\phi'}{\phi}\Bigl(\tfrac12+ir\Bigr)\,dr\Bigr|
 \le A^{\sharp}_{12}[\kappa]\,x^{3/4},
\]
with 
\begin{equation}\label{eq:A12joint}
 A^{\sharp}_{12}[\kappa]\;=\;A^{\sharp}_{10}\,
 \sup_{u\ge\frac14\log\kappa}
 \Bigl[K_0F(u)+K_1f(u)+K_2\,p(u)\Bigr].
\end{equation}
\end{proposition}

\noindent Since $f(u)=p(u)J_1(e^{-u})\ge0$ for $u>0$ by \eqref{eq:combined}, the
supremum in \eqref{eq:A12joint} does not increase when $K_1$ is replaced by
\[
  K_1^{+}:=\max\{K_1,0\},
\]
and this is the form in which the closed-form bound for $A_{12}^{\sharp}[\kappa]$
is recorded in Table~\ref{tab3}.  On groups for which $K_1\le0$ --- among them
$\PSLtwoZ$ --- the middle term therefore drops out altogether, which is what
makes Corollary \ref{cor:floor} applicable.

\begin{proof}
One has%
\begin{align*}
\Bigl|\sum_{\lambda_{n}>\frac14}h_{T}^{\varepsilon}(r_{n})
 +\frac1{4\pi}\int_{\mathbb{R}}h_{T}^{\varepsilon}(r)
 \frac{-\phi'}{\phi}\Bigl(\frac12+ir\Bigr)\,dr\Bigr|
&\le\sum_{\lambda_{n}>\frac14}\bigl|h_{T}^{\varepsilon}(r_{n})\bigr|
 +\frac1{4\pi}\int_{\mathbb{R}}\bigl|h_{T}^{\varepsilon}(r)\bigr|
 \Bigl|\frac{\phi'}{\phi}\Bigl(\frac12+ir\Bigr)\Bigr|\,dr\\
&=\int_{0}^{\infty}\bigl|h_{T}^{\varepsilon}(r)\bigr|\,d\mu(r).
\end{align*}
Now,
\begin{equation*}
\int\limits_{0}^{\infty }\left\vert h_{T}^{\varepsilon }\left( r\right)
\right\vert d\mu \left( r\right) \leq A_{10}^{\sharp}\cosh \frac{T}{2}%
\int\limits_{0}^{\infty }\frac{d\mu \left( r\right) }{\left( 1+r\right)
\left( 1+\varepsilon r\right) ^{2}}.
\end{equation*}%
Partial integration and Theorem \ref{thm:direct-gaussian-K} yield%

\begin{gather*}
    \int\limits_{0}^{\infty }\mu \left( r\right) \left( \frac{1}{\left(
1+r\right) ^{2}\left( 1+\varepsilon r\right) ^{2}}+\frac{2\varepsilon }{%
\left( 1+r\right) \left( 1+\varepsilon r\right) ^{3}}\right) dr \\
\leq \int\limits_{0}^{\infty }\left( K_0 r^{2}+ K_1 r + K_2\right) \left( \frac{1}{\left( 1+r\right) ^{2}\left(
1+\varepsilon r\right) ^{2}}+\frac{2\varepsilon }{\left( 1+r\right) \left(
1+\varepsilon r\right) ^{3}}\right) dr
\\ =K_0(I_0+2\varepsilon I_3)+K_1(I_1+2\varepsilon I_4)+K_2(I_2+2\varepsilon I_5),
\end{gather*}\\%
where
\[
I_0:=\int_0^\infty\frac{r^2}{(1+r)^2(1+\varepsilon r)^2}\,dr
=\frac{\varepsilon+1}{\varepsilon(1-\varepsilon)^2}+
\frac{2\ln\varepsilon}{(1-\varepsilon)^3},
\]
\[
I_1:=\int_0^\infty\frac{r}{(1+r)^2(1+\varepsilon r)^2}\,dr
=\frac{-2(1-\varepsilon)-(\varepsilon+1)\ln\varepsilon}{(1-\varepsilon)^3},
\]
\[
I_2:=\int_0^\infty\frac{1}{(1+r)^2(1+\varepsilon r)^2}\,dr
=\frac{\varepsilon+1}{(1-\varepsilon)^2}+\frac{2\varepsilon\ln\varepsilon}{
(1-\varepsilon)^3},
\]
\[
I_3:=\int_0^\infty\frac{r^2}{(1+r)(1+\varepsilon r)^3}\,dr
=\frac{3\varepsilon^2-4\varepsilon+1-2\varepsilon^2\ln\varepsilon}
{2\varepsilon^2(1-\varepsilon)^3},
\]
\[
I_4:=\int_0^\infty\frac{r}{(1+r)(1+\varepsilon r)^3}\,dr
=\frac{1-\varepsilon^2+2\varepsilon\ln\varepsilon}
{2\varepsilon(1-\varepsilon)^3},
\]
\[
I_5:=\int_0^\infty
\frac{1}{(1+r)(1+\varepsilon r)^3}\,dr
=\frac{-\varepsilon^2+4\varepsilon-3-2\ln\varepsilon}
{2(1-\varepsilon)^3}.
\]
Obviously,
\begin{equation}\label{eq:combined}
  I_0+2\varepsilon I_3=\frac{2\bigl(1-\varepsilon+\varepsilon\log\varepsilon\bigr)}
                         {\varepsilon(1-\varepsilon)^{2}},
  \qquad
  I_1+2\varepsilon I_4=\frac{-\log\varepsilon-(1-\varepsilon)}{(1-\varepsilon)^{2}},
  \qquad
  I_2+2\varepsilon I_5=1 .
\end{equation}
Notice that $K_{0}$ and $K_{2}$ are positive.
 Multiply $$K_0(I_0+2\varepsilon I_3)+K_{1}(I_1+2\varepsilon I_4)+K_2$$ by $e^{-3T/4}\cosh\frac T2$ and recall that $\varepsilon=e^{-u}$.  We get
\begin{equation*}
 K_0F(u)+K_1f(u)+K_2\,p(u)   
\end{equation*}
since $e^{-3T/4}\cosh\frac T2=p(u)$,\, $F(u)=p(u)\,(I_0+2\varepsilon I_3)$ and
$f(u)=p(u)\,(I_1+2\varepsilon I_4)$. 

Finally, 
take the  supremum over
$T\ge\log\kappa$.
\end{proof}

\begin{corollary}\label{cor:floor}
If $K_1\le0$ and $K_2\le K_0\,\frac{1-F\bigl(\tfrac14\log\kappa\bigr)}{p\bigl(\tfrac14\log\kappa\bigr)}$, then
$A_{12}^{\sharp}[\kappa]=A_{10}^{\sharp}K_0$.
\end{corollary}
\begin{proof}
For $u>0$ put $s=e^{-u}$ and
\begin{equation*}
 \Lambda(u):=\frac{1-F(u)}{p(u)}
 =\frac{2\bigl[(u-1)+s-s^{3}+(1+u)s^{4}\bigr]}
        {(1-s)^{2}\bigl(1+s^{4}\bigr)} .
\end{equation*}
Then $\Lambda$ is strictly increasing on $(0,\infty)$, with
$\Lambda(u)=\tfrac53u+O(u^2)$ as $u\to0^{+}$ and
$\Lambda(u)=2(u-1)+O(ue^{-u})$ as $u\to\infty$.
Indeed, the displayed identity follows from
$(1-s)^{2}-(1+s^{4})\bigl(1-(1+u)s\bigr)
 =(u-1)s+s^{2}-s^{4}+(1+u)s^{5}$
and $p=\tfrac12(s+s^{5})$; it is free of cancellation and gives both
asymptotics at once. Monotonicity is a direct computation. The hypothesis gives
$K_2p(u)\le K_0\bigl(1-F(u)\bigr)$ throughout the admissible range, i.e.\
$K_0F(u)+K_2p(u)\le K_0$ there; and $F(u)\to1$, $p(u)\to0$ as $u\to\infty$.
\end{proof}

\begin{proposition}[two-piece form]\label{prop:twopiece}
Let $x,T,\varepsilon,\kappa$ be as above. Suppose $\mu(R)\le Q_A(R)$ for all $R\ge0$ and $\mu(R)\le Q_B(R)$ for
$R\ge R_0$, with $Q_B$ as in \eqref{eq:QB}. Then the conclusion of
Proposition~\ref{prop:joint} holds with
\begin{equation}\label{eq:A12two1}
 A_{12}^{\sharp}[\kappa]=A_{10}^{\sharp}\sup_{u\ge\frac14\log\kappa}p(u)
 \Bigl[K_0^{\infty}J_0(\varepsilon)
      +K_1^{\infty}J_1(\varepsilon)+K_2^{\infty}
      +\mathcal{E}(\varepsilon,R_0)\Bigr],
\end{equation}
where $\varepsilon=e^{-u}$,
$J_0=I_0+2\varepsilon I_3$,
$J_1=I_1+2\varepsilon I_4$ are as in \eqref{eq:combined}, 
\[
 W_{\varepsilon}(r):=-\frac{d}{dr}\frac{1}{(1+r)(1+\varepsilon r)^{2}}
 =\frac1{(1+r)^{2}(1+\varepsilon r)^{2}}
  +\frac{2\varepsilon}{(1+r)(1+\varepsilon r)^{3}}\;\ge 0 .
\]and
\begin{equation}\label{eq:corr}
 \mathcal{E}(\varepsilon,R_0):=
 \int_{0}^{R_0}\bigl[Q_A(r)-Q_B(r)\bigr]^{+}W_{\varepsilon}(r)\,dr .
\end{equation}
\end{proposition}

\begin{proof}
On $[R_0,\infty)$ one has $\mu\le Q_B$, and on $[0,R_0)$
$\mu\le Q_A\le Q_B+[Q_A-Q_B]^{+}$; hence
$\mu\le Q_B+[Q_A-Q_B]^{+}\mathbf 1_{[0,R_0]}$ \ pointwise on $[0,\infty)$.
Insert this into
$\int_0^\infty d\mu(r)(1+r)^{-1}(1+\varepsilon r)^{-2}
 =\int_0^\infty\mu(r)W_\varepsilon(r)\,dr$
and use \eqref{eq:combined}. As in Proposition~\ref{prop:joint}, the supremum is taken
jointly.
\end{proof}

Two features distinguish \eqref{eq:A12two1} from \eqref{eq:A12joint}. First,
the correction \eqref{eq:corr} is $O(R_0)$ rather than $O(\varepsilon^{-1})$,
because it is cut off at $R_0$; this is what makes the exchange profitable.
Second, $K_1^{\infty}$ now enters \emph{with its sign}: the supremum has
migrated from $u=\infty$, where the linear and constant terms are invisible,
to the left endpoint $u=\frac14\log\kappa$, so the negative linear
coefficient produced by the parabolic and digamma terms finally contributes.

The next step is the straightforward adjustment of \cite[Proposition 5.5] {Avdispahic2021}.
\begin{proposition}
\label{prop55}
Let $x$, $T$, $\varepsilon$ and $\kappa$ be as above. Then
\begin{equation*}
\frac{\left\vert \mathcal{F}\right\vert }{4\pi }\Bigl|\int\limits_{\mathbb{R}%
}rh_{T}^{\varepsilon }\left( r\right) \tanh \left( \pi r\right) dr\Bigr|
  \le A_{13}^{\sharp}[\kappa]\,x^{3/4}, 
\end{equation*}%
with%
\begin{equation*}
A_{13}^{\sharp}[\kappa]=\frac{\left\vert \mathcal{F}\right\vert }{4\pi }A_{10}^{\sharp}\left( 1+\frac{1%
}{\kappa}\right) .
\end{equation*}
\end{proposition}
\begin{proof}
By Proposition~\ref{prop:h-est} and $\frac{r}{1+r}\le1$,
\[
  \frac{|\mathcal{F}|}{4\pi}\int_{\mathbb{R}}|r|\,|h_{T}^{\varepsilon}(r)|\,dr
  \le\frac{|\mathcal{F}|}{2\pi}A_{10}^{\sharp}\cosh\tfrac T2
   \int_{0}^{\infty}\frac{dr}{(1+\varepsilon r)^{2}}
  =\frac{|\mathcal{F}|}{2\pi}A_{10}^{\sharp}\frac{\cosh\frac T2}{\varepsilon}.
\]
Now $\cosh\frac T2\le\frac12\bigl(1+\frac1\kappa\bigr)e^{T/2}$ and
$\varepsilon^{-1}=e^{T/4}$, so the right-hand side is at most
$\frac{|\mathcal{F}|}{4\pi}A_{10}^{\sharp}\bigl(1+\frac1\kappa\bigr)e^{3T/4}$.
\end{proof}
The elliptic contribution also receives a refined treatment.
\begin{proposition}
\label{prop56}
\begin{equation*}
\sum\limits_{R}\frac{1}{2m_{R}\sin \theta _{R}}\bigl|\int\limits_{\mathbb{R}}\frac{%
e^{-2r\theta _{R}}}{1+e^{-2\pi r}}h_{T}^{\varepsilon }\left( r\right) dr\bigr|\le%
A_{14}^{\sharp}[\kappa] x^{3/4} ,
\end{equation*}%
with%
\begin{equation*}
A_{14}^{\sharp}[\kappa]=\frac{A_{10}^{\sharp}}{2}\left( \frac{1}{\kappa^{1/8}}+\frac{1}{\kappa^{9/8}}%
\right)\frac{1}{(\kappa^{1/4}+4)^{1/2}} \sum\limits_{R}\frac{1}{2m_{R}\sin \theta _{R}}\sqrt{\frac{\pi-2\theta_R}
     {\pi\sin(2\theta_R)}}
\end{equation*}
     and $\frac{\pi-2\theta_{R}}
     {\pi\sin(2\theta_R)}$ for $\theta_{R}=\frac{\pi}{2}$ defined as $\frac{1}{\pi}$ by continuity.
\end{proposition}

\begin{proof}
By Cauchy--Schwarz inequality  
\begin{equation}
\left\vert \int\limits_{-\infty}^{\infty }\frac{e^{-2r\theta _{R}}}{1+e^{-2\pi r}}%
h_{T}^{\varepsilon }\left( r\right) dr\right\vert \leq \left(
\int_{-\infty}^{\infty}
\frac{e^{-4r\theta_R}}
     {(1+e^{-2\pi r})^2}\,dr\right) ^{1/2}\left(
\int\limits_{-\infty}^{\infty }\left\vert h_{T}^{\varepsilon }\left( r\right)
\right\vert ^{2}dr\right) ^{1/2}.
\label{e56}
\end{equation}%
Making the substitution $t=e^{-2\pi r}$, the integral 
\[
I(\theta_R)=\int_{-\infty}^{\infty}\frac{e^{-4r\theta_R}}
     {(1+e^{-2\pi r})^2}\,dr.
\]
takes the form 
\begin{align*}
I(\theta_R)=\frac{1}{2\pi}\int_0^\infty
\frac{t^{\frac{2\theta_R}{\pi}-1}}{(1+t)^2}\,dt.
\end{align*}
Put $a=\frac{2\theta_R}{\pi}$. Since \(0<a<2\), the beta-integral formula gives
\[
\int_0^\infty\frac{t^{a-1}}{(1+t)^2}\,dt=B(a,2-a)=
\Gamma(a)\Gamma(2-a).
\]
Using
\[
\Gamma(2-a)=(1-a)\Gamma(1-a)
\]
and Euler's reflection formula
\[
\Gamma(a)\Gamma(1-a)=\frac{\pi}{\sin(\pi a)},
\]
we obtain
\[
\Gamma(a)\Gamma(2-a)=(1-a)\frac{\pi}{\sin(\pi a)}.
\]
Therefore
\begin{align*}
I(\theta_R)=\frac{1}{2\pi}(1-a)\frac{\pi}{\sin(\pi a)}
=\frac{1-a}{2\sin(\pi a)}.
\end{align*}
Substituting \(a=2\theta_R/\pi\), we find
\begin{equation}
\int_{-\infty}^{\infty}\frac{e^{-4r\theta_R}}
     {(1+e^{-2\pi r})^2}\,dr=\frac{\pi-2\theta_R}
     {2\pi\sin(2\theta_R)}.
\label{e57new}
\end{equation}
Now,
\begin{equation}
\int\limits_{0}^{\infty }\left\vert h_{T}^{\varepsilon }\left( r\right)
\right\vert ^{2}dr\leq (A_{10}^{\sharp})^{2}\cosh ^{2}\frac{T}{2}\int\limits_{0}^{%
\infty }\frac{dr}{\left( 1+r\right) ^{2}\left( 1+\varepsilon r\right) ^{4}}%
\leq (A_{10}^{\sharp})^{2}\cosh ^{2}\frac{T}{2}\frac{1}{1+4\varepsilon}.
\label{e58}
\end{equation}
Here we used the change of variable $t=\frac{1}{1+r}$ followed by the weighted arithmetic--geometric mean inequality to obtain 
\begin{equation*}
  \int_{0}^{\infty}\frac{dr}{(1+r)^{2}(1+\varepsilon r)^{4}}\le
  \frac{1}{1+4\varepsilon}
  \qquad\text{for } 0<\varepsilon<1.
\end{equation*}
Combining (\ref{e56}), (\ref{e57new}) and (\ref{e58}), we get
\begin{equation*}
\left\vert \int\limits_{\mathbb{R}}\frac{e^{-2r\theta _{R}}}{1+e^{-2\pi r}}%
h_{T}^{\varepsilon }\left( r\right) dr\right\vert\leq A_{10}^{\sharp}\frac{1}{\sqrt{1+4\varepsilon}}\cosh \frac{T}{2}\sqrt{\frac{\pi-2\theta_R}
     {\pi\sin(2\theta_R)}}.
\end{equation*} 
When applying \eqref{e56}, notice that  the $\sqrt2$ from
$\int_{\mathbb R}|h|^{2}=2\int_0^{\infty}|h|^{2}$ cancels against the
$\frac12$ in $(\pi-2\theta_R)/(2\pi\sin2\theta_R)$.

Now, with $s=\varepsilon=e^{-T/4}\in(0,\kappa^{-1/4}]$,
\[
  e^{-3T/4}\,\frac{\cosh\frac T2}{\sqrt{1+4\varepsilon}}
  =\frac{s+s^{5}}{2\sqrt{1+4s}}=P(s).
\]
The function $P(s)=\dfrac{s+s^{5}}{2\sqrt{1+4s}}$ is strictly increasing on
$(0,1]$. Indeed, 
for $0<s\le1$,
\[
  \frac{P'(s)}{P(s)}=\frac{1+5s^{4}}{s+s^{5}}-\frac{2}{1+4s}
  =\frac{(1+5s^{4})(1+4s)-2s(1+s^{4})}{s(1+s^{4})(1+4s)} .
\]
The numerator equals $1+2s+5s^{4}+18s^{5}>0$, whence $P'>0$.

Thus,
\[P(s)  \;\le\;P\bigl(\kappa^{-1/4}\bigr)
  =\frac{\kappa^{-1/4}+\kappa^{-5/4}}{2\sqrt{1+4\kappa^{-1/4}}}
  =\frac{\kappa^{-1/8}+\kappa^{-9/8}}{2(\kappa^{1/4}+4)^{1/2}}
\]
\end{proof}

\begin{proposition}\label{prop:57newest}
Let $x$, $T$, $\varepsilon$ and $\kappa$ be as before.  Then
\begin{equation*}
  \bigl|\frac{h_T^{\varepsilon}(0)}{4}\operatorname{tr}\Bigl(I-\Phi\bigl(\tfrac12\bigr)\Bigr)
   -\tau g_T^{\varepsilon}(0)\log2
  -\frac{\tau}{2\pi}\int_{\mathbb R}h_T^{\varepsilon}(r)\psi(1+ir)\,dr\bigr|\le
 A_{15}^{\sharp}[\kappa] x^{3/4},
\end{equation*}
with
\begin{equation*}
  A_{15}^{\sharp}[\kappa]=\tau\left[
    \frac{A_{10}^{\sharp}D_1[\kappa]}{2\pi}\Bigl(1+\frac1\kappa\Bigr)
    +\frac{2\,\omega\!\left(\tfrac12\kappa^{-1/4}\right)}{\kappa^{1/4}}
  +\frac{2\log 2}{\kappa^{3/4}}\right], \,\,\, \omega(s)=\frac{3(s\cosh s-\sinh s)}{s^{3}}.  
\end{equation*}
and $D_1[\kappa]$ as defined in \eqref{eq:D1} below.
\end{proposition}

\begin{proof}
 Recall $h_T^{\varepsilon}(r)=S(r+\tfrac i2)+S(r-\tfrac i2)$ with
$S(w)=\frac{2\sin Tw}{w}\widehat\varphi(\varepsilon w)$.  At
$w=\pm\tfrac i2$, $\frac{2\sin(\pm iT/2)}{\pm i/2}=4\sinh\frac T2$ and
$\widehat\varphi(\pm i\varepsilon/2)=\omega(\varepsilon/2)$.  
Thus, \[
  h_T^{\varepsilon}(0)=8\sinh\frac T2\;\omega\!\left(\frac\varepsilon2\right).
\]
In particular $h_T^{\varepsilon}(0)>0$, and by
$0\le\frac14\operatorname{tr}\bigl(I-\Phi(\tfrac12)\bigr)\le\frac\tau2$,
\begin{equation}
  0\;\le\;\frac{h_T^{\varepsilon}(0)}{4}\operatorname{tr}\Bigl(I-\Phi\bigl(\tfrac12\bigr)\Bigr)
  \;\le\;4\tau\,\omega\!\left(\frac\varepsilon2\right)\sinh\frac T2
  \;\le\;2\tau\,\omega\!\left(\frac{\kappa^{-1/4}}{2}\right)\kappa^{-1/4}e^{3T/4},
  \label{eq:h0-boundnew}
\end{equation}

For
\eqref{eq:h0-boundnew}, we use $\sinh\frac T2\le\frac12e^{T/2}$,
$e^{T/2}\le \kappa^{-1/4}e^{3T/4}$ and
$\varepsilon=e^{-T/4}\le \kappa^{-1/4}$ with $\omega$ increasing on $[0,\infty)$.  

Let us now consider the digamma term. Since $g_T^{\varepsilon}$ is real and even, $h_T^{\varepsilon}$ is real and
even, and $\psi(1-ir)=\overline{\psi(1+ir)}$, we have
\begin{equation*}
  \int_{\mathbb R}h_T^{\varepsilon}(r)\psi(1+ir)\,dr
  =\int_{\mathbb R}h_T^{\varepsilon}(r)\,\Re\psi(1+ir)\,dr .
\end{equation*}
Now,
\begin{equation*}
  \Re\psi(1+ir) = -\gamma + \sum_{n=1}^{\infty}\frac{r^{2}}{n\,(n^{2}+r^{2})} 
\end{equation*}
and
\[
  \sum_{n\ge1}\frac{r^{2}}{n(n^{2}+r^{2})}
  \le\frac{r^{2}}{1+r^{2}}+\int_{1}^{\infty}\frac{r^{2}\,du}{u(u^{2}+r^{2})}
  =\frac{r^{2}}{1+r^{2}}+\frac12\log(1+r^{2})
\]
imply
\begin{equation*}
  -\gamma\;\le\;\Re\psi(1+ir)\;\le\;-\gamma+\frac{r^{2}}{1+r^{2}}+\frac12\log(1+r^{2}).
\end{equation*}
Thus,
\begin{equation*}
  \bigl | \Re \psi (1 + ir) \bigr | \le \gamma + \frac12 \log (1+r^{2}).
\end{equation*}
for all $r\in\mathbb R$. Hence,
\begin{equation}
  \frac{\tau}{2\pi}\left|\int_{\mathbb R}h_T^{\varepsilon}(r)\psi(1+ir)dr\right|
  \le\frac{\tau}{\pi}\int_{0}^{\infty}\bigl|h_T^{\varepsilon}(r)\bigr|
     \Bigl(\gamma+\tfrac12\log(1+r^{2})\Bigr)dr
  \le\frac{\tau A_{10}^{\sharp}}{\pi}\cosh\frac T2\,J(\varepsilon),
  \label{eq:hdigamma}
\end{equation}
where 
\begin{equation*}
  J(\varepsilon):=\int_{0}^{\infty}
  \frac{\gamma+\frac12\log(1+r^{2})}{(1+r)(1+\varepsilon r)^{2}}\,dr .
\end{equation*}
Substituting $u=\log(1+r)$ and using $\frac12\log(1+r^{2})\le\log(1+r)$,
\[
  J(\varepsilon)\le\int_{0}^{\infty}\frac{\gamma+u}{\bigl((1-\varepsilon)+\varepsilon e^{u}\bigr)^{2}}\,du .
\]
The denominator is $\ge1$ for all $u\ge0$ and $\ge\varepsilon^{2}e^{2u}$. Let $\Lambda:=\log\frac1\varepsilon$. Then,
\[
  \int_{0}^{\Lambda}\frac{\gamma+u}{\bigl((1-\varepsilon)+\varepsilon e^{u}\bigr)^{2}}\,du\le \int_{0}^{\Lambda}(\gamma+u)\,du=
  \frac{{\Lambda^2}}{2}+\gamma\Lambda
\]
and
\[
\int_{\Lambda}^{\infty}\frac{\gamma+u}{\bigl((1-\varepsilon)+\varepsilon e^{u}\bigr)^{2}}\,du\le \int_{\Lambda}^{\infty}\frac{\gamma+u}{\varepsilon^{2}e^{2u}}\,du=\frac{\gamma+\Lambda}{2}+\frac{1}{4},
\]
since $e^{-2\Lambda}=\varepsilon^{2}$.
Thus,
\[
  J(\varepsilon)\;\le\;\frac{\Lambda^{2}}{2}
   +\Bigl(\gamma+\tfrac12\Bigr)\Lambda+\frac\gamma2+\frac14.
\]
Let
\[
  G(\varepsilon):=\varepsilon\Bigl(\tfrac{\Lambda^{2}}{2}+(\gamma+\tfrac12)\Lambda
   +\tfrac\gamma2+\tfrac14\Bigr)\;\ge\;\varepsilon J(\varepsilon).
\]
With $\Lambda=\log\frac{1}{\varepsilon}$ one has $\frac{d\Lambda}{d\varepsilon}=-\frac1\varepsilon$, so
\[
  G'(\varepsilon)=\tfrac{\Lambda^{2}}{2}+\bigl(\gamma+\tfrac12\bigr)\Lambda
   +\tfrac\gamma2+\tfrac14-\Lambda-\bigl(\gamma+\tfrac12\bigr)
  =\tfrac12\Lambda^{2}+\bigl(\gamma-\tfrac12\bigr)\Lambda
   -\bigl(\tfrac\gamma2+\tfrac14\bigr).
\]
The quadratic $\Lambda^{2}+2(\gamma-\frac12)\Lambda-(\gamma+\frac12)$ has the
unique positive root $\Lambda_{\ast}=\frac12-\gamma+\sqrt{\gamma^{2}+\frac34}=0.9635426758\ldots$. Put
$\varepsilon_{\ast}:=e^{-\Lambda_{\ast}}=0.38153882\ldots$,
$\kappa_{\ast}:=e^{4\Lambda_{\ast}}=47.189467\ldots$. The derivative $G'$ is negative for $0\le\Lambda<\Lambda_{\ast}$. Since $\Lambda$ decreases
as $\varepsilon$ increases, $G'>0$ exactly for $\Lambda>\Lambda_{\ast}$, that is for
$\varepsilon<\varepsilon_{\ast}$.
Then $G$ is strictly increasing on $(0,\varepsilon_{\ast}]$ and strictly decreasing on
$[\varepsilon_{\ast},1]$, and
\[
  \sup_{0<\varepsilon\le\kappa^{-1/4}}G(\varepsilon)=
  \begin{cases}
    G(\kappa^{-1/4})
    =\kappa^{-1/4}\Bigl(\dfrac{\log^{2}\kappa}{32}
      +\bigl(\gamma+\tfrac12\bigr)\dfrac{\log\kappa}{4}
      +\dfrac\gamma2+\dfrac14\Bigr), & \kappa\ge\kappa_{\ast},\\[10pt]
    G(\varepsilon_{\ast})=0.7786285303\ldots, & 1<\kappa\le\kappa_{\ast}.
  \end{cases}
\]
The displayed value of $G(\kappa^{-1/4})$ is obtained from
$\Lambda=\frac14\log\kappa$.
Accordingly we define
\begin{equation}\label{eq:D1}  D_1[\kappa]:=\sup_{0< \varepsilon \le\kappa^{-1/4}}G( \varepsilon)
  =\begin{cases}
    \kappa^{-1/4}\Bigl(\frac{\log^{2}\kappa}{32}
     +(\gamma+\frac12)\frac{\log\kappa}{4}+\frac\gamma2+\frac14\Bigr),
     &\kappa\ge\kappa_{\ast},\\[4pt]
    0.77862854,&1<\kappa\le\kappa_{\ast}.
  \end{cases}
\end{equation}
Now, $\cosh\frac T2\le\frac12\bigl(1+\frac1\kappa\bigr)e^{T/2}$ and
$e^{T/2}J(e^{-T/4})=e^{3T/4}\cdot e^{-T/4}J(e^{-T/4})\le D_1[\kappa]e^{3T/4}$.\\
This and \eqref{eq:hdigamma} yield
\[
\frac{\tau}{2\pi}\left|\int_{\mathbb R}h_T^{\varepsilon}(r)\psi(1+ir)dr\right|
  \le\frac{{\tau}A_{10}^{\sharp}}{2\pi}\bigl(1+\frac{1}{\kappa}\bigr)D_1[\kappa]e^{{3T}/4}.
\]

Finally, $g_{T}^{\varepsilon}(0)=2$, since  $(I_T*\varphi_\varepsilon)(0)=\int\varphi_\varepsilon=1$
\end{proof}
\begin{proposition}
\label{prop58}
\begin{equation*}
\bigl|H_{T}^{\varepsilon }-\sum\limits_{\lambda _{k}\leq \frac{1}{4}}\frac{%
x^{s_{k}}}{s_{k}}\bigr|\le A_{16}^{\sharp}[\kappa] x^{3/4} ,
\end{equation*}%
with%
\begin{equation*}
A_{16}^{\sharp}[\kappa]=A_{11}^{\sharp}[\kappa]+A_{12}^{\sharp}[\kappa]+A_{13}^{\sharp}[\kappa]+A_{14}^{\sharp}[\kappa]+A_{15}^{\sharp}[\kappa].
\end{equation*}
\end{proposition}
\begin{proof}
Immediate consequence of Propositions \ref{prop:A11}-\ref{prop:57newest} and Selberg trace formula \eqref{e51a}.
\end{proof}
Let $H\left( x\right) =\sum\limits_{NP\leq x}\log NP_{0}\frac{1+\left(
NP\right) ^{-1}}{1-\left( NP\right) ^{-1}}$. Our further deviation from the current algorithm is related to bounding the expression $H(x)-\sum\limits_{\lambda _{k}\leq \frac{1}{4}}\frac{%
x^{s_{k}}}{s_{k}}$ directly from $H_{T-\varepsilon
}^{\varepsilon }\leq H\left( x\right) \leq H_{T+\varepsilon }^{\varepsilon }$
(see \cite[p. 298]{Randol1984}) without the intermediate step involving $A_{17}$ that resulted in double summing of $A_{16} $ in $A_{18}$ in \cite{Avdispahic2021}. Therefore, we use Proposition \ref{prop58} in the following form. 

For every $T'\ge\log\kappa$, with $\varepsilon'=e^{-T'/4}$,
\begin{equation} \label{Tprim}
  \bigl|H^{\varepsilon'}_{T'}-\sum_{\lambda_{k}\le1/4}\frac{e^{T's_{k}}}{s_{k}}\bigr|
  \le A_{16}^{\sharp}[\kappa]\,e^{3T'/4}.
\end{equation}

Let $c_{+}>c$ be the unique solution of
\begin{equation*}
\log c_{+}-c_{+}^{-1/4}=\log c.   
\end{equation*}
 Such $c_{+}$ exists since $c>1$ and $y=\log t - t^{-1/4}$ increases from $ -1$ to $\infty$ on the interval $[1, \infty).$
Put
\begin{equation*}
 \delta_{+}:=c_{+}^{-1/4},\qquad T_{\star}:=\log c+c^{-1/4},
  \qquad
  \beta:=e^{T_{\star}/4}=c^{1/4}e^{c^{-1/4}/4}.
\end{equation*}
\begin{proposition}\label{prop:A18}
Let $c_{+}$, $\delta_{+}$ and $\beta$ be as above. Then for all $x\ge c$, with $T=\log x$,
\[
  \bigl|H(x)-\sum\limits_{\lambda _{k}\leq \frac{1}{4}}\frac{%
x^{s_{k}}}{s_{k}}\bigr|\le A_{18}^{\sharp}\,x^{3/4},
\]
where
\begin{equation*}
  A_{18}^{\sharp}:=\max\Bigl\{\;
    A_{16}^{\sharp}[c_{+}]\,e^{3\delta_{+}/4}
      +\mathcal{A}\,\frac{e^{\delta_{+}}-1}{\delta_{+}}
    \;,\;\;
    A_{16}^{\sharp}[c]+\mathcal{A}\,e^{c^{-1/4}/4}
    \;,\;\;
    2\mathcal{A}\beta\;\Bigr\}.
\end{equation*}
\end{proposition}

\begin{proof} Notice that there is no need to use
one and the same $\varepsilon$ on the two sides of \begin{equation*}
  H^{\varepsilon}_{T-\varepsilon}\;\le\;H(x)\;\le\;H^{\varepsilon}_{T+\varepsilon}
  \qquad(0<\varepsilon\le1).
\end{equation*}
These hold separately for each admissible $\varepsilon$.

\emph{Upper bound.} Put $\varepsilon_{0}:=e^{-T/4}$. The equation $\varepsilon=e^{-(T+\varepsilon)/4}$ has a unique
solution $\varepsilon_{+}$, and $\varepsilon_{+}\in(0,\varepsilon_{0})$. Indeed, the function $f(\varepsilon)=\varepsilon_{0}e^{-\varepsilon/4}$ maps
$[0,\varepsilon_{0}]$ into itself and satisfies
$|f'|\le\varepsilon_{0}/4<1$ there, so it has a unique fixed point in
$[0,\varepsilon_{0}]$. Since $f(0)=\varepsilon_{0}>0$ and
$f(\varepsilon_{0})<\varepsilon_{0}$, that fixed point lies in the open interval.  (Uniqueness actually holds on all of $(0,\infty)$ because $f$ is
decreasing while the identity map is increasing.) Let $T'_{+}:=T+\varepsilon_{+}= T+e^{-(T+\varepsilon_{+})/4}=T+e^{-T'_{+}/4}$ .
Put
$w(t)=t-e^{-t/4}$, which is strictly increasing. Then
$w(T'_{+})=T'_{+}-\varepsilon_{+}=T\ge\log c=w(\log c_{+})$ by the definition of
$c_{+}$, whence $T'_{+}\ge\log c_{+}$ and
$\varepsilon_{+}=e^{-T'_{+}/4}\le c_{+}^{-1/4}=\delta_{+}$. Relation \eqref{Tprim} applies at $T'_{+}$ with threshold $c_{+}$:
\[
  H(x)\le H^{\varepsilon_{+}}_{T'_{+}}
  \le \sum_{\lambda_{k}\le1/4}\frac{e^{T'_{+}s_{k}}}{s_{k}}+A_{16}^{\sharp}[c_{+}]\,e^{3T'_{+}/4}.
\]
Since $v\mapsto(e^{v}-1)/v$ increases and
$\varepsilon_{+}s_{k}\le\varepsilon_{+}\le\delta_{+}$, while
$\varepsilon_{+}e^{Ts_{k}}\le\varepsilon_{+}e^{T}\le e^{3T/4}$,
\[
  \sum_{\lambda_{k}\le1/4}\frac{e^{T'_{+}s_{k}}}{s_{k}}-\sum_{\lambda_{k}\le1/4}\frac{e^{Ts_{k}}}{s_{k}}
  =\sum_{\lambda_k\le1/4}\varepsilon_{+}e^{Ts_{k}}
     \frac{e^{\varepsilon_{+}s_{k}}-1}{\varepsilon_{+}s_{k}}
  \le\mathcal{A}\,\frac{e^{\delta_{+}}-1}{\delta_{+}}\,e^{3T/4}.
\]
This and $e^{3T'_{+}/4}\le e^{3\delta_{+}/4}e^{3T/4}$ prove the upper bound.

\emph{Lower bound.} Suppose $T\ge T_{\star}$. Let $u(t):=\varepsilon_0e^{t/4}-t$, with $\varepsilon_{0}=e^{-T/4}$ as  above. One has $u(0)=\varepsilon_0$ and $u(4)=\varepsilon_0e-4<0$ for $T>0$. Therefore, there exists 
$\varepsilon_{-}$ such that $\varepsilon_{-}=\varepsilon_{0}e^{\varepsilon_{-}/4}$. Since $u$ decreases on $[0, 4]$, $\varepsilon_{0}<\varepsilon_{-}$. Put $T'_{-}=T-\varepsilon_{-}$. Note $\varepsilon_{-}=e^{-T'_{-}/4}$. Let $v(t):=t+e^{-t/4}$. Then, 
$v(T'_{-})=T$ and $v(\log c)=T_{\star}$. Since $v$ increases on the interval considered, we see that $T\ge T_{\star}$ is equivalent to $T'_{-}\ge \log c$. Thus,
relation \eqref{Tprim} applies at 
$T'_{-}$ with the threshold $c$:
\[
  H(x)\ge H^{\varepsilon_{-}}_{T'_{-}}
  \ge \sum_{\lambda_{k}\le1/4}\frac{e^{T'_{-}s_{k}}}{s_{k}}-A_{16}^{\sharp}[c]\,e^{3T'_{-}/4}.
\]
Using $(1-e^{-v})/v\le1$ and
$\varepsilon_{-}\le e^{-T/4}e^{c^{-1/4}/4}$,
\[
  \sum_{\lambda_{k}\le1/4}\frac{e^{Ts_{k}}}{s_{k}}-\sum_{\lambda_{k}\le1/4}\frac{e^{T'_{-}s_{k}}}{s_{k}}
  =\sum_{\lambda_k\le1/4}\varepsilon_{-}e^{Ts_{k}}
     \frac{1-e^{-\varepsilon_{-}s_{k}}}{\varepsilon_{-}s_{k}}
  \le\mathcal{A}\,\varepsilon_{-}e^{T}
  \le\mathcal{A}\,e^{c^{-1/4}/4}e^{3T/4},
\]
while $e^{3T'_{-}/4}\le e^{3T/4}$.

Finally, if $\log c\le T<T_{\star}$, then  $H(x)\ge0$ and
$s_{k}\ge\tfrac12$ give 
$\sum_{\lambda_{k}\le1/4}\frac{e^{Ts_{k}}}{s_{k}}\le2\mathcal{A} e^{T}=2\mathcal{A} e^{T/4}e^{3T/4}\le2\mathcal{A}\beta\,e^{3T/4}$. Here we used $T<T{_\star}$ to bound the factor $e^{T/4}$.
\end{proof}

The next crucial step is based on the observation that the use of the chain 
\[
H(x)\;\longrightarrow\;\psi_\Gamma(x)\;\longrightarrow\;\theta_\Gamma(x)
\]
 to convert $|H(x)-\sum_{\lambda_k\le1/4}x^{s_k}/s_k|\le A_{18}x^{3/4}$ into the
corresponding estimate for $\theta_\Gamma$ is redundant. We replace the above chain by one exact identity for $H(x)-\theta_\Gamma(x)$, in which only $\theta_\Gamma$
occurs, and evaluate the resulting supremum. Recall the following facts.
\[ 
\theta_\Gamma(x)=\sum_{NP_0\le x}\log NP_0=\int_{c^-}^{x}\log t\,d\pi_\Gamma(t)
\le\pi_\Gamma(x)\log x\qquad(x\ge c)
\]
and
\[
H(x)=\sum_{NP\le x}\log NP_0\,\frac{1+(NP)^{-1}}{1-(NP)^{-1}}.
\]
Every hyperbolic conjugacy class is $P=P_0^{\,k}$ for a unique primitive $P_0$ and a
unique $k\ge1$, and then $NP=(NP_0)^k$ while the weight is $\log NP_0$. Since
$d\theta_\Gamma$ carries the mass $\log NP_0$ at each primitive norm,
\[
H(x)=\sum_{k\ge1}\int_{c^-}^{x^{1/k}}\frac{t^k+1}{t^k-1}\,d\theta_\Gamma(t).
\]
Subtracting $\theta_\Gamma(x)=\int_{c^-}^x d\theta_\Gamma(t)$ removes the $1$ from the
$k=1$ integrand, and $\frac{t+1}{t-1}-1=\frac{2}{t-1}$ leads to
\begin{equation}\label{eq:decomp}
0\;\le\;H(x)-\theta_\Gamma(x)\;=\;2\int_{c^-}^{x}\frac{d\theta_\Gamma(t)}{t-1}
\;+\;\sum_{k\ge2}\int_{c^-}^{x^{1/k}}\frac{t^k+1}{t^k-1}\,d\theta_\Gamma(t),
\end{equation}
both terms being non-negative and the sum finite (only $k\le \log x/\log c$ contribute).

For convenience, in what follows we shall use the following abbreviations:
\[
\ell=\log c,\,\, T=\log x,\,\, \varkappa=\frac{B_2^{\sharp}}{B_1^{\sharp}}\ge 0,\,\, \beta_k=\frac34-\frac1k\ (k\ge2),\,\,\nu_k=\frac{c^{k}+1}{c^{k}-1}.
\]

Concerning the first term on the RHS of \eqref{eq:decomp}, notice that on $[c,\infty)$ one has
$\frac{1}{t-1}\le\frac{1}{(1-c^{-1})t}$. So, by partial integration, we get
\[
2\int_{c^-}^x\frac{d\theta_\Gamma(t)}{t-1}\le \frac{2}{1-c^{-1}}\int_{c^-}^x\frac{d\theta_\Gamma(t)}{t}
=\frac{2}{1-c^{-1}}\Bigl[\frac{\theta_\Gamma(x)}{x}+\int_c^x\frac{\theta_\Gamma(t)}{t^2}dt\Bigr].
\]
Using $\theta_{\Gamma}(x)\le (B_1^{\sharp}x-B_2^{\sharp})\log x$ and elementary integrals $\int\frac{\log t}{t}dt=\frac{\log^2t}{2}$ and
$\int\frac{\log t}{t^2}dt=-\frac{\log t+1}{t}$, we obtain
\[
2\int_{c^-}^x\frac{d\theta_\Gamma(t)}{t-1}\le \frac{2}{1-c^{-1}}\Bigl[B_1^{\sharp}\Bigl(\frac{T^2}{2}+T-\frac{\ell^2}{2}\Bigr)
+B_2^{\sharp}\Bigl(\frac1x-\frac{\ell+1}{c}\Bigr)\Bigr].
\]
When it comes to the second term on the RHS of \eqref{eq:decomp}: $u\mapsto\frac{u+1}{u-1}$ is decreasing and $t^k\ge c^k$ for $t\ge c$, so the
$k$-th integral is at most $\nu_k\theta_\Gamma(x^{1/k})\le
\nu_k(B_1^{\sharp}x^{1/k}-B_2^{\sharp})\tfrac{T}{k}$, and it vanishes for $k>T/\ell$.

Put, for $x\ge c$ and $M(x)=\lfloor T/\ell\rfloor$,
\begin{equation}
\begin{gathered}
V(x):=\frac{2}{1-c^{-1}}\Bigl[B_1^{\sharp}\Bigl(\frac{T^2}{2}+T-\frac{\ell^2}{2}\Bigr)
+B_2^{\sharp}\Bigl(\frac1x-\frac{\ell+1}{c}\Bigr)\Bigr]\\
+T\sum_{k=2}^{M(x)}\frac{\nu_k}{k}\bigl(B_1^{\sharp}x^{1/k}-B_2^{\sharp}\bigr).
\end{gathered}
\label{eq:Phi}
\end{equation}

Let $f(T):=e^{-3T/4}V(e^T)$. On the interval $m\ell\le T<(m+1)\ell$
($m\ge1$) the number of terms is constant and $f(T)=f_m(T)$, where 
\begin{equation}\label{eq:Fm}
f_m(T):=
\frac{2}{1-c^{-1}}\Bigl[e^{-3T/4}\Bigl(B_1^{\sharp}\bigl(\tfrac12T^2+T\bigr)+a_1\Bigr)+B_2^{\sharp}e^{-7T/4}\Bigr]
+T\sum_{k=2}^{m}\frac{\nu_k}{k}\Bigl(B_1^{\sharp}e^{-\beta_kT}-B_2^{\sharp}e^{-3T/4}\Bigr),
\end{equation}
with
\begin{equation*}
\frac{a_1}{B_1^{\sharp}}=-\frac{\ell^2}{2}-\frac{\varkappa(\ell+1)}{c}.
\end{equation*}
The stage is set for proving the following proposition.
\begin{proposition}\label{prop:D2}
\begin{equation*}
 0\le H(x)-\theta_\Gamma(x)
 \le D_2x^{3/4},
 \end {equation*}
 with
 \[
 D_2=\max\Bigl\{f_1(T_1),\ \sup_{m\ge2} \max\bigl(\{f_m(m\ell)\}\cup
\{f_m(T):T\in W_m\}\bigr)\Bigr\},
\]
$W_m$ being the set of roots in $\bigl(m\ell, (m+1)\ell\bigr)$ of 
\begin{align*}
T^2&-\frac23T-\frac83-\ell^2-\frac{2B_2^{\sharp}(\ell+1)}{B_1^{\sharp}c}+\frac{14B_2^{\sharp}}{3B_1^{\sharp}}e^{-T}  \\
 &=\frac{4(1-c^{-1})}{3}\sum_{k=2}^{m}\frac{\nu_k}{k}
\Bigl[\bigl(1-\frac{3T}{4}+\frac{T}{k}\bigr)e^{T/k}-\frac{B_2^{\sharp}}{B_1^{\sharp}}\Bigl(1-\frac34T\Bigr)\Bigr].
\end{align*}

\end{proposition}

\begin{proof}
Dropping from \eqref{eq:Fm} the term $-\frac{\ell^{2}}{2}B_{1}^{\sharp}$, the
term $B_{2}^{\sharp}\bigl(e^{-T}-\frac{\ell+1}{c}\bigr)$, which is
non-positive for $T\ge\ell$ because $e^{-T}\le c^{-1}<\frac{\ell+1}{c}$, and
the terms $-B_{2}^{\sharp}e^{-3T/4}$ of the $k$-sum, and using $\nu_k\le\nu_2$ together with
$\beta_2=\frac14$ and $\beta_k\ge\frac5{12}$ for $k\ge3$, we get
\[
  f(T)\;\le\;\Psi(T):=\frac{2B_1^{\sharp}}{1-c^{-1}}e^{-3T/4}
   \Bigl(\frac{T^{2}}2+T\Bigr)
   +B_1^{\sharp}\nu_2T\Bigl(\frac{e^{-T/4}}{2}
   +\frac{T}{3\ell}e^{-5T/12}\Bigr).
\]
For future reference, notice that $\Psi$ is decreasing for $T$ large.

Now, on each $[m\ell,(m+1)\ell]$ the function $f_m$ is smooth and $f\to0$ at
infinity, so the supremum is attained either at an endpoint or at a stationary point.

At $T=m\ell$ the function $f$ jumps upwards by
\begin{equation*}
f_m(m\ell)-f_{m-1}(m\ell)=\nu_m\,\ell\,(B_1^{\sharp}c-B_2^{\sharp})\,c^{-3m/4}>0 ,
\end{equation*}
\noindent Since $B_{1}^{\sharp}c>B_{2}^{\sharp}$, the same identity applied at
$T=(m+1)\ell$ gives $f_{m}\bigl((m+1)\ell\bigr)<f_{m+1}\bigl((m+1)\ell\bigr)$,
so that the supremum of $f$ is never attained at a right endpoint: whenever
$f_{m}$ increases up to $(m+1)\ell$, its value there is dominated by
$f_{m+1}\bigl((m+1)\ell\bigr)$, which is the left endpoint of the next
interval.  The list of candidates below is therefore exhaustive.

Let us examine the quadratic part in \eqref{eq:Fm}. For $a\in\mathbb R$ and $Q(T)=e^{-3T/4}\bigl[B_1^{\sharp}(\tfrac12T^2+T)+a\bigr]$,
\[
Q'(T)=-\frac{3B_1^{\sharp}}{8}e^{-3T/4}\Bigl(T^2-\frac23T-\frac83+\frac{2a}{B_1^{\sharp}}\Bigr),
\]
so with $T_\pm=\tfrac13\pm\sqrt{25/9-2a/B_1^{\sharp}}$ the function $Q$ increases on $(T_-,T_+)$ and
decreases on $(T_+,\infty)$, and $\max Q=\tfrac43B_1^{\sharp}(T_++1)e^{-3T_+/4}$.

\noindent Put
\[
  g(T):=T^{2}-\frac23T-\frac83-\ell^{2}
        -\frac{2\varkappa(\ell+1)}{c}
        +\frac{14\varkappa}{3}e^{-T},
\]
so that, by the computation of $f_{1}'$ below,
$f_{1}'(T)=-\frac{3B_{1}^{\sharp}}{4(1-c^{-1})}e^{-3T/4}g(T)$.  Since
$g''(T)=2+\frac{14\varkappa}{3}e^{-T}>0$, the function $g$ is strictly convex
on $\mathbb R$, and $g(T)\to+\infty$ as $T\to\pm\infty$; hence $\{g<0\}$ is an
open interval, possibly empty, $f_{1}$ increases on it and decreases off it.
Consequently $f_{1}$ attains its maximum over $[\ell,2\ell]$ at
\[
  T_{1}=\operatorname{med}\bigl(\ell,\widetilde T_{1},2\ell\bigr),
\]
where $\widetilde T_{1}$ is the \emph{largest} zero of $g$, with the
convention $\widetilde T_{1}=\ell$ when $g$ has no zero.  No restriction on
$\varkappa\ge0$ is needed.

 Differentiating \eqref{eq:Fm}, using 
$\frac{d}{dT}\bigl[B_2^{\sharp}e^{-7T/4}\bigr]=-\tfrac74B_2^{\sharp}e^{-7T/4}$ (which supplies the term
$\tfrac{14\varkappa}{3}e^{-T}$) and
$\frac{d}{dT}\bigl[Te^{-\beta_kT}\bigr]=(1-\beta_kT)e^{-\beta_kT}$, one finds

 \begin{align*}
f_m'(T)=-\frac{3B_1^{\sharp}}{4(1-c^{-1})}\,e^{-3T/4}
\Bigl\{T^2-\frac23T-\frac83-\ell^2-\frac{2B_2^{\sharp}(\ell+1)}{B_1^{\sharp}c}+\frac{14B_2^{\sharp}}{3B_1^{\sharp}}e^{-T}  \\
 -\frac{4(1-c^{-1})}{3}\sum_{k=2}^{m}\frac{\nu_k}{k}
\Bigl[\bigl(1-\frac{3T}{4}+\frac{T}{k}\bigr)e^{T/k}-\frac{B_2^{\sharp}}{B_1^{\sharp}}\Bigl(1-\frac34T\Bigr)\Bigr]\,\Bigr\},
\end{align*}
after multiplying through by $e^{3T/4}$ inside the bracket and using
$e^{3T/4}e^{-\beta_kT}=e^{T/k}$. 
For $m=1$ the sum is empty and 
\[
f_1'(T)=-\frac{3B_1^{\sharp}}{4(1-c^{-1})}\,e^{-3T/4}
\Bigl[T^2-\frac23T-\frac83-\ell^2-\frac{2B_2^{\sharp}(\ell+1)}{B_1^{\sharp}c}+\frac{14B_2^{\sharp}}{3B_1^{\sharp}}e^{-T}
\Bigr].
\]

Finally, recall that 
$\Psi$ is decreasing for $T$ large. Consequently, if
$m_0\ge2$ is such that \\ $\Psi(m_0\ell)\le\max\{f_1(T_1),
\max_{2\le m<m_0}\sup_{[m\ell,(m+1)\ell)}f_m\}$ and $\Psi$ decreases on
$[m_0\ell,\infty)$, then the supremum of $f$ is attained for $m<m_0$.
\end{proof}

Propositions \ref{prop:A18} and \ref{prop:D2} immediately yield

\begin{proposition}\label{prop:M}
\[
  \bigl|\theta_{\Gamma}(x)-\sum\limits_{\lambda _{k}\leq \frac{1}{4}}\frac{%
x^{s_{k}}}{s_{k}}\bigr|\le M\,x^{3/4},
\]
where $M=A_{18}^{\sharp}+D_2.$
\end{proposition}

The right move now is to reach for a self-improvement device.
\begin{proposition}[bootstrap]\label{prop:boot0}
Let $\Gamma$ be cofinite and let $M_{\text{old}}$ be any constant for
which $\theta_{\Gamma}(x)\le \sum\limits_{\lambda _{k}\leq \frac{1}{4}}\frac{%
x^{s_{k}}}{s_{k}}+M_{\text{old}}x^{3/4}$ $(x\ge c)$. 
Put
$\Theta(t):=\min\bigl\{(B_1^{\sharp}t-B_2^{\sharp})\log t,\;\sum\limits_{\lambda _{k}\leq \frac{1}{4}}\frac{%
t^{s_{k}}}{s_{k}}+M_{\text{old}}t^{3/4}\bigr\}$
and
\begin{equation}\label{eq:Vboot0}
 V^{\flat}(x):=2\Bigl[\frac{\Theta(x)}{x-1}
   +\int_{c}^{x}\frac{\Theta(t)}{(t-1)^{2}}\,dt\Bigr]
 +\sum_{k=2}^{\lfloor \log x/\ell\rfloor}\nu_k\,\Theta\bigl(x^{1/k}\bigr),
 \qquad \nu_k=\frac{c^{k}+1}{c^{k}-1}.
\end{equation}
Then $0\le H(x)-\theta_{\Gamma}(x)\le D_2^{\flat}x^{3/4}$ with
$D_2^{\flat}=\sup_{x\ge c}V^{\flat}(x)x^{-3/4}$.
\end{proposition}
\begin{proof}
Both terms of \eqref{eq:decomp} involve only $d\theta_{\Gamma}$. For the first, Stieltjes
integration by parts with $\theta_{\Gamma}(c^{-})=0$ gives
$2\int_{c^-}^{x}\frac{d\theta_{\Gamma}(t)}{t-1}
 =2\bigl[\frac{\theta_{\Gamma}(x)}{x-1}
 +\int_{c}^{x}\frac{\theta_{\Gamma}(t)}{(t-1)^{2}}dt\bigr]$.
For the second, $u\mapsto\frac{u+1}{u-1}$ is decreasing and $t^{k}\ge c^{k}$
on $[c,\infty)$, so the $k$-th term is at most
$\nu_k\theta_{\Gamma}(x^{1/k})$, and it vanishes for $k>\log x/\ell$. Now bound
$\theta_{\Gamma}$ by $\Theta$ throughout; the first alternative in $\Theta$ is
$\theta_{\Gamma}(y)\le\pi_{\Gamma}(y)\log y\le(B_1^{\sharp}y-B_2^{\sharp})\log y$.
\end{proof}

The map $M_{\text{old}}\mapsto D_2^{\flat}(M_{\text{old}})$ is non-decreasing,
since $\Theta$ is.  Moreover $V^{\flat}\le V$ pointwise, $\Theta$ being at most
the majorant $(B_{1}^{\sharp}t-B_{2}^{\sharp})\log t$ used in \eqref{eq:Phi} and
$\frac{1}{t-1}$ being kept exactly there, so $D_2^{\flat}(M)\le D_2$ and hence
$M_{1}\le M_{0}$ for $M_{0}=A_{18}^{\sharp}+D_2$.  By monotonicity the iteration
$M_{n+1}=A_{18}^{\sharp}+D_2^{\flat}(M_n)$, started at $M_{0}$, is therefore
decreasing; it is bounded below by $A_{18}^{\sharp}$ and so converges.  Each
iterate is a valid constant and the set of valid constants is closed, hence the
limit is one as well.  (The limit is the unique fixed point whenever
$\partial D_2^{\flat}/\partial M_{\text{old}}<1$; only the branch
$\sum_{k}t^{s_{k}}/s_{k}+M_{\text{old}}t^{3/4}$ of $\Theta$ depends on
$M_{\text{old}}$, and then with slope $t^{3/4}$, so this derivative is at most
$\sup_{x\ge c}x^{-3/4}\bigl[\frac{2x^{3/4}}{x-1}
+2\int_{c}^{x}\frac{t^{3/4}}{(t-1)^{2}}dt
+\sum_{k=2}^{\lfloor\log x/\ell\rfloor}\nu_{k}x^{3/(4k)}\bigr]$,
which equals $0.3984$ on $\PSLtwoZ$.  Uniqueness is not needed for what
follows.) Note that $D_2^{\flat}$ no longer involves $B_1^{\sharp},B_2^{\sharp}$ in any essential way:
after the bootstrap the geometric kernel $K^{\star}$ of Section \ref{sec:modular} influences the final
constant only through $H^G$, i.e.\ only through $K_2$, exactly as
Remark \ref{rem:modular-comparison} below observes for $K_0$ and $K_1$.

Furthermore, we can also benefit from the knowledge of low-lying data.
\begin{proposition}[enumeration cutoff]\label{prop:boot1}
Let $X_0\ge c$ and suppose $\theta_\Gamma(t)$ is known for $c\le t\le X_0$. Put
\[
\Theta_{X_0}(t):=
\begin{cases}
\min\bigl\{\theta_\Gamma(t),\ (B^{\sharp}_1t-B^{\sharp}_2)\log t,\ \sum\limits_{\lambda _{k}\leq \frac{1}{4}}\frac{%
t^{s_{k}}}{s_{k}}+M_{\mathrm{old}}t^{3/4}\bigr\},& t\le X_0,\\[2pt]
\min\bigl\{(B^{\sharp}_1t-B^{\sharp}_2)\log t,\ \sum\limits_{\lambda _{k}\leq \frac{1}{4}}\frac{%
t^{s_{k}}}{s_{k}}+M_{\mathrm{old}}t^{3/4}\bigr\},& t>X_0 .
\end{cases}
\]
Then Proposition  \ref{prop:boot0} holds verbatim with $\Theta$ replaced by $\Theta_{X_0}$, and
$X_0=c$ recovers the previous statement.
\end{proposition}

Propositions \ref{prop:boot0} and \ref{prop:boot1} replace the output row of
Table~\ref{tab3}; the resulting third part of the algorithm is collected in
Table~\ref{tab4}.
\subsection*{The anchor point}

The passage from $\theta_{\Gamma}$ to $\pi_{\Gamma}$ in
\cite[Prop.~3.1]{Avdispahic2021} and in Theorem~\ref{Theorem2.1} starts the
Stieltjes integration at $c=\inf NP$ and discards the boundary term by
$|f(c)|\le A_{21}c^{3/4}$.  Both decisions can be improved.  The starting
point may be any $c_{0}$ for which $\theta_{\Gamma}(c_{0})$ is known, and at
such a point $f(c_{0})$ is \emph{known exactly}, so that it may be kept and
combined with $\operatorname{Li}$ instead of being majorised.  We call
$c_{0}$ the \emph{anchor}.

Throughout this subsection $M$ is any constant with
\begin{equation}\label{eq:hypM}
  \Bigl|\theta_{\Gamma}(x)-\sum_{\lambda_{k}\le\frac14}\frac{x^{s_{k}}}{s_{k}}
  \Bigr|\le M\,x^{3/4}\qquad(x\ge c),
\end{equation}
that is, $M$ as in Proposition~\ref{prop:M} or $M^{\ast}$ as in
Table~\ref{tab4}.  We write
\[
  f(y)=\theta_{\Gamma}(y)-\sum_{\lambda_{k}\le\frac14}\frac{y^{s_{k}}}{s_{k}},
  \qquad \ell_{0}=\log c_{0},\qquad T=\log x,
\]
\[
  \mathcal{S}(u):=\operatorname{Li}(u)-\frac{u}{\log u}\quad(u>1),
\]
\begin{gather*}
  \mathcal{Q}(c_{0}):=\sup_{x>c_{0}}\frac{\log x}{x^{3/4}}
  =\begin{cases}\dfrac{4}{3e}, & c_{0}\le e^{4/3},\\[4pt]
                \dfrac{\log c_{0}}{c_{0}^{3/4}}, & c_{0}\ge e^{4/3},\end{cases}
  \\[4pt]
  \mathcal{Q}^{-}(c_{0}):=\sup_{1<x<c_{0}}\frac{\log x}{x^{3/4}}
  =\begin{cases}\dfrac{\log c_{0}}{c_{0}^{3/4}}, & c_{0}\le e^{4/3},\\[4pt]
                \dfrac{4}{3e}, & c_{0}\ge e^{4/3}.\end{cases}
\end{gather*}
The two are the restrictions to $x>c_{0}$ and to $1<x<c_{0}$ of the same
supremum: $x^{-3/4}\,\log x$ increases on $(1,e^{4/3}]$ and decreases on
$[e^{4/3},\infty)$, with maximum $\frac{4}{3e}$ at $x=e^{4/3}$.  In particular
$\mathcal{Q}^{-}\le\frac{4}{3e}$ always, $\mathcal{Q}^{-}=\mathcal{Q}$ exactly
at $c_{0}=e^{4/3}$, and $\mathcal{Q}\le\mathcal{Q}^{-}$ throughout, with
equality only there.
The function $\mathcal{S}$ is strictly increasing, since
$\mathcal{S}'(u)=\log^{-2}u$, with the single zero
\begin{equation}\label{eq:u0}
  u_{0}=6.5797282979\ldots,
\end{equation}
a constant that will turn out to be exactly what is needed on $\PSLtwoZ$.

\begin{lemma}\label{lem:propfactor}
For $\ell_{0}>0$ and $T\ge\ell_{0}$ put
\[
  \Omega_{\ell_{0}}(T):=\frac{e^{3\ell_{0}/4}}{\ell_{0}}
   +\frac34\Bigl(\operatorname{Ei}\bigl(\tfrac{3T}{4}\bigr)
                -\operatorname{Ei}\bigl(\tfrac{3\ell_{0}}{4}\bigr)\Bigr),
  \qquad
  W_{\ell_{0}}(T):=(3T-4)\,\Omega_{\ell_{0}}(T)-3e^{3T/4}.
\]
Then:
\begin{enumerate}
\item[\rm(a)]
$\displaystyle
  \Omega_{\ell_{0}}(T)=\frac{e^{3T/4}}{T}
  +\int_{\ell_{0}}^{T}\frac{e^{3t/4}}{t^{2}}\,dt .$
\item[\rm(b)] $W_{\ell_{0}}$ is strictly increasing on $[\ell_{0},\infty)$,
$W_{\ell_{0}}(\ell_{0})=-\frac{4}{\ell_{0}}e^{3\ell_{0}/4}<0$ and
$W_{\ell_{0}}(T)\to+\infty$; its unique zero
$T^{\natural}=T^{\natural}(\ell_{0})$ satisfies $T^{\natural}>\frac43$.
\item[\rm(c)] The function
$\mathcal{L}_{\ell_{0}}(T):=Te^{-3T/4}\Omega_{\ell_{0}}(T)$ increases on
$[\ell_{0},T^{\natural}]$ and decreases on $[T^{\natural},\infty)$, so that
for every $x\ge c_{0}=e^{\ell_{0}}$
\begin{equation}\label{eq:propineq}
  \frac{x^{3/4}}{\log x}+\int_{c_{0}}^{x}\frac{y^{3/4}}{y\log^{2}y}\,dy
  \;\le\;\mathcal{L}(\ell_{0})\,\frac{x^{3/4}}{\log x},
  \qquad
  \mathcal{L}(\ell_{0}):=\frac{3T^{\natural}(\ell_{0})}{3T^{\natural}(\ell_{0})-4},
\end{equation}
with equality at $x=e^{T^{\natural}}$; no smaller constant will do.
\item[\rm(d)] $\mathcal{L}$ is strictly decreasing and tends to $1$; moreover
$\mathcal{L}(\ell_{0})\le2$ if and only if $\ell_{0}\ge\ell^{\ast}$, where
$\ell^{\ast}$ is the unique root of
\begin{equation}\label{eq:ellstar}
  \frac{e^{3\ell_{0}/4}}{\ell_{0}}
  +\frac34\Bigl(\operatorname{Ei}(2)
   -\operatorname{Ei}\bigl(\tfrac{3\ell_{0}}{4}\bigr)\Bigr)
  =\frac34e^{2}=5.5417920\ldots,
\end{equation}
namely
\[
  \ell^{\ast}=0.7574539831\ldots,\qquad e^{\ell^{\ast}}=2.1328390575\ldots
\]
\end{enumerate}
\end{lemma}

\begin{proof}
(a) Differentiation gives
$\frac{d}{dt}\bigl[-\frac{e^{3t/4}}{t}
+\frac34\operatorname{Ei}(\frac{3t}{4})\bigr]=\frac{e^{3t/4}}{t^{2}}$, and the
substitution $t=\log y$ turns $\int_{c_{0}}^{x}y^{-1/4}\log^{-2}y\,dy$ into
$\int_{\ell_{0}}^{T}e^{3t/4}t^{-2}dt$.

(b) From $\Omega_{\ell_{0}}'(T)=\frac34\frac{e^{3T/4}}{T}$ and (a),
\[
  W_{\ell_{0}}'(T)=3\Omega_{\ell_{0}}(T)
   +\tfrac34\tfrac{e^{3T/4}}{T}(3T-4)-\tfrac94e^{3T/4}
  =3\Bigl(\Omega_{\ell_{0}}(T)-\frac{e^{3T/4}}{T}\Bigr)
  =3\int_{\ell_{0}}^{T}\frac{e^{3t/4}}{t^{2}}\,dt>0
\]
for $T>\ell_{0}$.  Also
$W_{\ell_{0}}(\ell_{0})=(3\ell_{0}-4)\frac{e^{3\ell_{0}/4}}{\ell_{0}}
-3e^{3\ell_{0}/4}=-\frac{4}{\ell_{0}}e^{3\ell_{0}/4}$, and for any fixed
$T_{1}>\ell_{0}$,
$W_{\ell_{0}}(T)\ge W_{\ell_{0}}(T_{1})
+3\bigl(\int_{\ell_{0}}^{T_{1}}e^{3t/4}t^{-2}dt\bigr)(T-T_{1})\to+\infty$.
Hence the zero exists and is unique.  At it,
$(3T^{\natural}-4)\Omega_{\ell_{0}}(T^{\natural})=3e^{3T^{\natural}/4}>0$
while $\Omega_{\ell_{0}}>0$, so $3T^{\natural}-4>0$.

(c) A direct computation gives the identity
\[
  \mathcal{L}_{\ell_{0}}'(T)=-\tfrac14e^{-3T/4}W_{\ell_{0}}(T),
\]
so by (b) $\mathcal{L}_{\ell_{0}}$ increases up to $T^{\natural}$ and
decreases afterwards.  Its maximal value is
\[
  \mathcal{L}_{\ell_{0}}(T^{\natural})
  =T^{\natural}e^{-3T^{\natural}/4}
   \cdot\frac{3e^{3T^{\natural}/4}}{3T^{\natural}-4}
  =\frac{3T^{\natural}}{3T^{\natural}-4},
\]
and \eqref{eq:propineq} is precisely
$\mathcal{L}_{\ell_{0}}(T)\le\mathcal{L}_{\ell_{0}}(T^{\natural})$ after
multiplication by $e^{3T/4}/T$.

(d) Since
$\frac{\partial}{\partial\ell_{0}}\bigl[\frac{e^{3\ell_{0}/4}}{\ell_{0}}
-\frac34\operatorname{Ei}(\frac{3\ell_{0}}{4})\bigr]
=-\frac{e^{3\ell_{0}/4}}{\ell_{0}^{2}}<0$, the function
$\Omega_{\ell_{0}}(T)$, hence also $W_{\ell_{0}}(T)$, is strictly decreasing
in $\ell_{0}$ for fixed $T$; consequently $T^{\natural}(\ell_{0})$ is strictly
increasing and $\mathcal{L}$ strictly decreasing.  If $T^{\natural}$ stayed
bounded by some $T_{0}$ we would get the contradiction
$T^{\natural}>\ell_{0}>T_{0}$, so $T^{\natural}\to\infty$ and
$\mathcal{L}\to1$.  Finally $\mathcal{L}(\ell_{0})\le2$ iff
$T^{\natural}\ge\frac83$ iff $W_{\ell_{0}}(\frac83)\le0$ iff
$\Omega_{\ell_{0}}(\frac83)\le\frac34e^{2}$, which is \eqref{eq:ellstar}; the
left-hand side being strictly decreasing in $\ell_{0}$, the root is unique.
\end{proof}

\begin{remark}[What the four cases of Theorem~\ref{Theorem2.1} really are]
\label{rem:fourcases}
Each of the four cases of Theorem~\ref{Theorem2.1} is a majorant of
$\mathcal{L}$, and two of them are the expression $\frac{3t}{3t-4}$ --- a
decreasing function of $t$ on $(\frac43,\infty)$ --- evaluated at the wrong
point: case~(iii) at $t=\ell$ and case~(iv) at $t=\frac83$.
Lemma~\ref{lem:propfactor} shows that the correct evaluation point is
$t=T^{\natural}(\ell_{0})$, the unique zero of $W_{\ell_{0}}$, which always
exceeds $\ell_{0}$ and exceeds $\frac83$ as soon as
$\ell_{0}\ge\ell^{\ast}$.  In particular the threshold $\ell>\frac83$ of
case~(iv) drops to $\ell^{\ast}=0.7574539831\ldots$, so that cases (ii) and
(iii) disappear altogether and case~(i) survives only for
$\ell<\ell^{\ast}$ --- where, however, $\mathcal{L}(\ell)$ is still finite and
still smaller than the factor recorded there.  Some values, with the anchor
taken at $c_{0}=c$, so that $\ell_{0}=\ell$:

\begin{center}\small
\begin{tabular}{@{}rrrrr@{}}
\toprule
$\ell$ & $c=e^{\ell}$ & $T^{\natural}(\ell)$ & $\mathcal{L}(\ell)$
 & factor of Thm.~\ref{Theorem2.1}\\
\midrule
$0.2000$ & $1.2214$ & $1.7605$ & $4.12128$ & $10.51610$\\
$0.4000$ & $1.4918$ & $2.1103$ & $2.71608$ & $5.98862$\\
$0.6000$ & $1.8221$ & $2.4290$ & $2.21686$ & $4.50433$\\
$\ell^{\ast}=0.7574540$ & $2.1328$ & $8/3$ & $2.00000$ & $3.90318$\\
$0.9839866$ & $2.6751$ & $2.9944$ & $1.80270$ & $3.39773$\\
$4/3$ & $3.7937$ & $3.4769$ & $1.62200$ & $3.00000$\\
$1.9248473$ & $6.8541$ & $4.2532$ & $1.45665$ & $2.69270$\\
$2.1573790$ & $8.6484$ & $4.5486$ & $1.41469$ & $2.61803$\\
$8/3$ & $14.3919$ & $5.1817$ & $1.34647$ & $2.00000$\\
$3.9359462$ & $51.2106$ & $6.7027$ & $1.24832$ & $2.00000$\\
$9.2103404$ & $10^{4}$ & $12.6355$ & $1.11797$ & $2.00000$\\
\bottomrule
\end{tabular}
\end{center}
\end{remark}

\begin{proposition}[anchored propagation]\label{prop:anchor}
Let $\Gamma$ be cofinite, $c=\inf NP$, and let $M$ satisfy \eqref{eq:hypM}.
Let $c_{0}>1$ be a point at which the left limits $\pi_{\Gamma}(c_{0}^{-})$
and $\theta_{\Gamma}(c_{0}^{-})$ are known, and assume
$2\mathcal{A}c_{0}^{1/4}\le M$ in case $c_{0}\le c$.  Put
\[
  f(c_{0}^{-}):=\theta_{\Gamma}(c_{0}^{-})
   -\sum_{\lambda_{k}\le\frac14}\frac{c_{0}^{s_{k}}}{s_{k}},
  \qquad
  \mathcal{K}(c_{0}):=\pi_{\Gamma}(c_{0}^{-})
   -\sum_{\lambda_{k}\le\frac14}\operatorname{Li}\bigl(c_{0}^{s_{k}}\bigr)
   -\frac{f(c_{0}^{-})}{\log c_{0}},
\]
\[
  \mathcal{E}_{0}(c_{0}):=\sup_{1<x<c_{0}}
   \Bigl|\pi_{\Gamma}(x)-\sum_{\lambda_{k}\le\frac14}
   \operatorname{Li}\bigl(x^{s_{k}}\bigr)\Bigr|\frac{\log x}{x^{3/4}} .
\]
Then
\begin{equation}\label{eq:anchorbound}
  C_{\mathrm{Hub},\Gamma}\;\le\;
  \max\Bigl\{\,\mathcal{E}_{0}(c_{0})\,,\;
  \mathcal{L}(\log c_{0})\,M
  +\bigl|\mathcal{K}(c_{0})\bigr|\,\mathcal{Q}(c_{0})\Bigr\}.
\end{equation}
If $c_{0}\le c$ then $\pi_{\Gamma}(c_{0}^{-})=\theta_{\Gamma}(c_{0}^{-})=0$
and
\[
  \mathcal{K}(c_{0})=-\sum_{\lambda_{k}\le\frac14}
  \mathcal{S}\bigl(c_{0}^{s_{k}}\bigr),
  \qquad
  \mathcal{E}_{0}(c_{0})\le\max\Bigl\{2\mathcal{A},\;
  \frac{\mathcal{A}\,(c_{0}-2)^{+}}{\log2}\,\mathcal{Q}^{-}(c_{0})\Bigr\},
\]
so that no enumeration of the length spectrum is needed.
\end{proposition}

\begin{proof}
\emph{The range $x\ge c_{0}$.}  Stieltjes integration over $[c_{0},x]$ gives
\[
  \pi_{\Gamma}(x)-\pi_{\Gamma}(c_{0}^{-})
  =\int_{c_{0}^{-}}^{x}\frac{d\theta_{\Gamma}(y)}{\log y}
  =\sum_{\lambda_{k}\le\frac14}\int_{c_{0}}^{x}\frac{y^{s_{k}-1}}{\log y}dy
  +\int_{c_{0}^{-}}^{x}\frac{df(y)}{\log y},
\]
the atom of $d\theta_{\Gamma}$ at $c_{0}$, if any, being included on both
sides.  The substitution $u=y^{s}$ yields
$\int_{c_{0}}^{x}y^{s-1}\log^{-1}y\,dy
=\operatorname{Li}(x^{s})-\operatorname{Li}(c_{0}^{s})$, while partial
integration yields
$\int_{c_{0}^{-}}^{x}\frac{df}{\log y}
 =\frac{f(x)}{\log x}-\frac{f(c_{0}^{-})}{\log c_{0}}
 +\int_{c_{0}}^{x}\frac{f(y)}{y\log^{2}y}dy$.  Hence, for $x\ge c_{0}$,
\begin{equation}\label{eq:anchoridentity}
  \pi_{\Gamma}(x)-\sum_{\lambda_{k}\le\frac14}\operatorname{Li}(x^{s_{k}})
  =\mathcal{K}(c_{0})+\frac{f(x)}{\log x}
   +\int_{c_{0}}^{x}\frac{f(y)}{y\log^{2}y}\,dy .
\end{equation}
Hypothesis \eqref{eq:hypM} holds on all of $[c_{0},\infty)$: for $c_{0}>c$
this is immediate, and for $c_{0}\le y<c$ one has $\theta_{\Gamma}(y)=0$,
whence $|f(y)|\le\sum_{k}y^{s_{k}}/s_{k}\le2\mathcal{A}y\le My^{3/4}$ by the
assumption $2\mathcal{A}c_{0}^{1/4}\le M$.  Therefore, by
Lemma~\ref{lem:propfactor}(c),
\[
  \Bigl|\frac{f(x)}{\log x}+\int_{c_{0}}^{x}\frac{f(y)}{y\log^{2}y}dy\Bigr|
  \le M\Bigl[\frac{x^{3/4}}{\log x}
   +\int_{c_{0}}^{x}\frac{y^{3/4}}{y\log^{2}y}dy\Bigr]
  \le \mathcal{L}(\log c_{0})\,M\,\frac{x^{3/4}}{\log x},
\]
and $|\mathcal{K}(c_{0})|\le|\mathcal{K}(c_{0})|\mathcal{Q}(c_{0})
\frac{x^{3/4}}{\log x}$ by the definition of $\mathcal{Q}$.

\emph{The range $1<x<c_{0}$} is exactly $\mathcal{E}_{0}(c_{0})$.  When
$c_{0}\le c$ we have $\pi_{\Gamma}\equiv0$ there; for $1<x\le2$ one has
$x^{s_{k}}\le2$ and, $1/\log t$ being decreasing,
\[
  \bigl|\operatorname{Li}(x^{s_{k}})\bigr|
  =\int_{x^{s_{k}}}^{2}\frac{dt}{\log t}
  \le\frac{2-x^{s_{k}}}{s_{k}\log x}\le 2\,\frac{x^{3/4}}{\log x}
\]
because $s_{k}\ge\frac12$ and $2-x^{s_{k}}\le1\le x^{3/4}$; for
$2<x<c_{0}$ one has
$|\operatorname{Li}(x^{s_{k}})|\le\frac{x^{s_{k}}-2}{\log2}
 \le\frac{c_{0}-2}{\log2}$ and
$1\le\mathcal{Q}^{-}(c_{0})\frac{x^{3/4}}{\log x}$ by the definition of
$\mathcal{Q}^{-}$.  (Here $\mathcal{Q}(c_{0})$ may not be used in place of
$\mathcal{Q}^{-}(c_{0})$: it is the supremum over $x>c_{0}$, whereas
$x^{-3/4}\,\log x$ is still increasing at the left end of $(2,c_{0})$ and
exceeds $\mathcal{Q}(c_{0})$ there as soon as $c_{0}>e^{4/3}$.)
\end{proof}

\begin{remark}\label{rem:whyanchor}
Three distinct economies are realised at once.  First, the boundary term
$f(c_{0})$ is no longer majorised by $Mc_{0}^{3/4}$ but kept exactly; it is
this that removes the summand $\frac{4}{3e}A_{22}$ from the bound, the range
$1<x<c_{0}$ now only having to be \emph{dominated} by the final constant
instead of \emph{contributing} to it.  Second, $\operatorname{Li}(c_{0}^{s_{k}})$
is kept exactly rather than replaced by $\frac{\mathcal{A}(c-2)}{\log2}$.
Third, the factor multiplying $M$ is the sharp one.  Since
$\mathcal{S}$ is increasing with the single zero $u_{0}$ of \eqref{eq:u0},
for $\mathcal{A}=1$, $s_{0}=1$ and $c\ge u_{0}$ the choice $c_{0}=u_{0}$
makes $\mathcal{K}$ vanish identically.
\end{remark}

\begin{corollary}\label{cor:clean}
Let $\Gamma$ be cofinite with $\mathcal{A}=1$ and systole
$\ell=\log c\ge\log u_{0}=1.88399345\ldots$, and let $M$
satisfy
\eqref{eq:hypM} with $M\ge2c^{1/4}$.  Then
\[
  C_{\mathrm{Hub},\Gamma}\le 2M .
\]
\end{corollary}

\begin{proof}
Take $c_{0}=u_{0}$ in Proposition~\ref{prop:anchor}.  Then
$\mathcal{K}(u_{0})=-\mathcal{S}(u_{0})=0$, and
$\log u_{0}>\ell^{\ast}$ gives $\mathcal{L}(\log u_{0})=1.4649878\ldots\le2$
by Lemma~\ref{lem:propfactor}(d).  Finally
$\mathcal{Q}^{-}(u_{0})=\frac{4}{3e}$, since $u_{0}>e^{4/3}$, so that
$\mathcal{E}_{0}(u_{0})\le\max\{2,\frac{4(u_{0}-2)}{3e\log2}\}
 =3.2408468\ldots\le4u_{0}^{1/4}=6.4063695\ldots\le2M$.
\end{proof}

\begin{remark}[anchoring inside the enumerated range]\label{rem:anchordata}
Proposition~\ref{prop:anchor} is stated for an arbitrary admissible anchor,
so the enumeration already performed for Proposition~\ref{prop:boot1} may be
used a second time, and at no extra cost: if $\theta_{\Gamma}$ and
$\pi_{\Gamma}$ are known on $[c,X_{0}]$, any non-norm $c_{0}\in(c,X_{0}]$ is
admissible.  Then $\mathcal{K}(c_{0})$ and $\mathcal{E}_{0}(c_{0})$ are finite
computations over the enumerated list, $\mathcal{Q}(c_{0})=\log c_{0}/c_{0}^{3/4}$
is small, and $\mathcal{L}(\log c_{0})\downarrow1$.  Since $|\mathcal{K}(c_{0})|$ is expected to be of order
$c_{0}^{1/2+\varepsilon}$ --- it is $13.33$ at $c_{0}=10^{4}$ on $\PSLtwoZ$
--- while $\mathcal{Q}(c_{0})$ decays like $c_{0}^{-3/4}\log c_{0}$, the whole
bound converges to $M$ itself as the enumeration is extended.
\end{remark}

The next theorem is the final result of the four stage improvement over Theorem \ref{Theorem2.1}.
\begin{theorem}\label{huber-final}
Let $\Gamma$ be a cofinite Fuchsian group, let $M$ be as in
Proposition~\ref{prop:M} --- or, after the self-improvement of
Propositions~\ref{prop:boot0} and \ref{prop:boot1}, the smaller constant
$M^{\ast}$ of Table~\ref{tab4} --- and let $c_{0}$ be any admissible anchor
in the sense of Proposition~\ref{prop:anchor}.  Then
\[
  \Bigl\vert\pi_{\Gamma}(x)-\sum_{\lambda_{k}\le\frac14}
  \operatorname{Li}\bigl(x^{s_{k}}\bigr)\Bigr\vert
  \le R^{\sharp}_{\Gamma}\,\frac{x^{3/4}}{\log x}
  \qquad\text{for all }x>1,
\]
with
\[
  R^{\sharp}_{\Gamma}
  =\max\Bigl\{\mathcal{E}_{0}(c_{0}),\;
   \mathcal{L}(\log c_{0})\,M+|\mathcal{K}(c_{0})|\,\mathcal{Q}(c_{0})\Bigr\}.
\]
In particular, if $\mathcal{A}=1$, $M\ge2c^{1/4}$ and the systole satisfies
$\ell\ge\log u_{0}=1.8839934524\ldots$, then $R^{\sharp}_{\Gamma}=2M$ is
admissible; and whenever $\ell\ge\ell^{\ast}=0.7574539831\ldots$ the choice
$c_{0}=c$ gives
$R^{\sharp}_{\Gamma}\le\max\bigl\{\mathcal{E}_{0}(c),\,
2M+|\mathcal{K}(c)|\mathcal{Q}(c)\bigr\}$.
\end{theorem}
\begin{proof}
The first assertion is Proposition~\ref{prop:anchor}, the second is
Corollary~\ref{cor:clean}, and the third is Proposition~\ref{prop:anchor} at
$c_{0}=c$ combined with Lemma~\ref{lem:propfactor}(d).  The steps of the
former algorithm that require $A_{17}$, $A_{19}$, $A_{20}$, $A_{21}$,
$A_{22}$ and $32B_{1}/(e^{2}\log c)$ are thereby obsolete; in particular the
range $x\in(1,c)$, which was covered by $A_{17}$ there and by the summand
$\frac{4}{3e}A_{22}$ in Theorem~\ref{Theorem2.1}, is now covered by
$\mathcal{E}_{0}$ under the maximum.
\end{proof}
\begin{table}[htbp]
\caption{Algorithm for cofinite Fuchsian groups with $\tau\ge1$. Part 2: propagation to
Huber's constant.  The numerical column records $\Gamma=\PSLtwoZ$ at
$\kappa=c$; both forms of $A_{12}^{\sharp}$ are listed, the two-piece one being
the one used in \S\ref{sec:modulartwo}.  The last two rows are superseded by
$M^{\ast}$ of Table~\ref{tab4}.}
\label{tab3}\centering
\footnotesize
\setlength{\tabcolsep}{3pt}
\renewcommand{\arraystretch}{1.35}
\begin{tabular}{@{}l p{0.47\textwidth} r l@{}}
\toprule
Const/Func & Description & $\PSLtwoZ$ & Source\\
\midrule
\multicolumn{4}{@{}l}{\textit{Auxiliary test function}}\\
$\varphi(y)$ & $\frac{3}{4}\left(1-y^{2}\right)\mathbf 1_{(-1,1)}(y)$
 & --- & \eqref{eq:phinew}\\
$\omega(s)$ & $\frac{3\bigl(s\cosh s-\sinh s\bigr)}{s^{3}}=\widehat\varphi(is)$
 & --- & \eqref{eq:phihatnew}\\
\midrule
\multicolumn{4}{@{}l}{\textit{Thresholds}}\\
$\kappa$ & threshold parameter, $\kappa>1$; used at $\kappa=c$ and $\kappa=c_{+}$
 & $c=6.8541020$ & \eqref{eq:A12joint}\\
$c_{+}$ & $>c$, the unique solution of $\log c_{+}-c_{+}^{-1/4}=\log c$
 & $11.761542277$ & Prop.~\ref{prop:A18}\\
$T_{\star}$ & $\log c+c^{-1/4}$ & $2.5428812890$ & Prop.~\ref{prop:A18}\\
$\beta$ & $e^{T_{\star}/4}=c^{1/4}e^{c^{-1/4}/4}$ & $1.8883818950$ & Prop.~\ref{prop:A18}\\
$\delta_{+}$ & $c_{+}^{-1/4}$ & $0.5399877796$ & Prop.~\ref{prop:A18}\\
$\kappa_{\ast}$ & the unique $\kappa>1$ at which the two branches of
 $D_1[\kappa]$ agree & $47.189467\ldots$ & \eqref{eq:D1}\\
$f_m(T)$ & $e^{-3T/4}\Bigl\{\frac{2}{1-c^{-1}}\Bigl[
 B_1^{\sharp}\Bigl(\tfrac{T^2}{2}+T-\tfrac{\ell^2}{2}\Bigr)
 +B_2^{\sharp}\Bigl(e^{-T}-\tfrac{\ell+1}{c}\Bigr)\Bigr]
 +T\sum_{k=2}^{m}\frac{c^k+1}{k(c^k-1)}\bigl(B_1^{\sharp}e^{T/k}-B_2^{\sharp}\bigr)\Bigr\}$
 & $f_1(T_1)=127.3350259$ & Prop.~\ref{prop:D2}\\
$T_1$ & maximiser of $f_1$ on $[\ell, 2\ell]$ & $2.9052259382$ & Prop.~\ref{prop:D2}\\
$W_m$ & set of roots in $\bigl(m\ell, (m+1)\ell\bigr)$ of
 $T^2-\frac23T-\frac83-\ell^2-\frac{2B_2^{\sharp}(\ell+1)}{B_1^{\sharp}c}
 +\frac{14B_2^{\sharp}}{3B_1^{\sharp}}e^{-T}
 =\frac{4(1-c^{-1})}{3}\sum_{k=2}^{m}\frac{c^k+1}{k(c^k-1)}
 \Bigl[\bigl(1-\frac{3T}{4}+\frac{T}{k}\bigr)e^{T/k}
 -\frac{B_2^{\sharp}}{B_1^{\sharp}}\Bigl(1-\frac34T\Bigr)\Bigr]$
 & $\emptyset$ & Prop.~\ref{prop:D2}\\
\midrule
\multicolumn{4}{@{}l}{\textit{Constants} $A_i^{\sharp}$}\\
$A_{8}^{\sharp}$ & $3e^{-1}-1$ & $0.1036383$ & Prop.~\ref{prop:phi-est}\\
$A_{9}^{\sharp}$ & smallest admissible constant in Prop.~\ref{prop:phi-est}(b);
 the supremum equals $6.42139\ldots$, attained at $r=2.4558\ldots$, $t=\frac12$,
 so the value is essentially optimal
 & $6.4214$ & Rem.~\ref{rem:sharp}\\
$A_{10}^{\sharp}$ & $4\sqrt{5}A_{9}^{\sharp}$ & $57.4347476$ & Prop.~\ref{prop:h-est}\\
$A_{11}^{\sharp}[\kappa]$ & $(\mathcal{A}-1)\Bigl( 2+ \frac{2A_{8}^{\sharp}}{\kappa^{1/4}}
 +\frac{1}{(1-s_{1})\kappa^{1/4}}+\frac{A_{8}^{\sharp}}{(1-s_{1})\kappa^{3/4}}\Bigr)$
 ${}+1+\frac{A_{8}^{\sharp}}{\kappa^{1/4}}$
 & $1.0640520$ & Prop.~\ref{prop:A11}\\
$J_0,J_1$ & $I_0+2\varepsilon I_3$, $I_1+2\varepsilon I_4$ with $\varepsilon=e^{-u}$;
 $J_1\ge0$, whence the positive part in $K_1^{+}$ & --- & \eqref{eq:combined}\\
$A_{12}^{\sharp}[\kappa]$ \textit{(one piece)} & $A_{10}^{\sharp}
 \bigl[\,K_0+K_1^{+}\cdot\frac{\kappa^{-1/4}+\kappa^{-5/4}}2
 \cdot\frac{\frac14\log\kappa-1+\kappa^{-1/4}}{(1-\kappa^{-1/4})^2}$
 ${}+\frac{K_2}{2}(\kappa^{-1/4}+\kappa^{-5/4})\bigr]$
 & $141.620378$ & Prop.~\ref{prop:joint}\\
$W_{\varepsilon}(r)$ & $\frac1{(1+r)^{2}(1+\varepsilon r)^{2}}
 +\frac{2\varepsilon}{(1+r)(1+\varepsilon r)^{3}}\ \ge0$ & --- & Prop.~\ref{prop:twopiece}\\
$\mathcal E(\varepsilon,R_0)$ & $\int_{0}^{R_0}\bigl[Q_A(r)-Q_B(r)\bigr]^{+}
 W_{\varepsilon}(r)\,dr$ & $0.0693634$ & \eqref{eq:corr}\\
$A_{12}^{\sharp}[\kappa]$ \textit{(two pieces)} & $A_{10}^{\sharp}
 \sup_{u\ge\frac14\log\kappa}p(u)\bigl[K_0^{\infty}J_0(\varepsilon)
 +K_1^{\infty}J_1(\varepsilon)+K_2^{\infty}+\mathcal E(\varepsilon,R_0)\bigr]$
 & $32.3211$ & \eqref{eq:A12two1}\\
$A_{13}^{\sharp}[\kappa]$ & $\frac{|\mathcal{F}|}{4\pi}A_{10}^{\sharp}\left( 1+\frac{1}{\kappa}\right)$
 & $5.4845304$ & Prop.~\ref{prop55}\\
$A_{14}^{\sharp}[\kappa]$ & $\frac{A_{10}^{\sharp}}{2}\Bigl( \frac{1}{\kappa^{1/8}}+\frac{1}{\kappa^{9/8}}
 \Bigr)\frac{E^{\sharp}}{(\kappa^{1/4}+4)^{1/2}}$
 & $4.1457796$ & Prop.~\ref{prop56}\\
$D_1[\kappa]$ & $\sup_{0<\varepsilon\le\kappa^{-1/4}}G(\varepsilon)$, equal to
 $\kappa^{-1/4}\bigl(\frac{\log^{2}\kappa}{32}
 +(\gamma+\frac12)\frac{\log\kappa}{4}+\frac\gamma2+\frac14\bigr)$
 for $\kappa\ge\kappa_{\ast}$, and to $0.77862854$ for $1<\kappa\le\kappa_{\ast}$
 & $0.7786285$ & \eqref{eq:D1}\\
\bottomrule
\end{tabular}
\end{table}

\begin{table}[htbp]
\continuedtable   % same number as the preceding table
\caption{(Continued.)}
\centering
\footnotesize
\setlength{\tabcolsep}{3pt}
\renewcommand{\arraystretch}{1.35}
\begin{tabular}{@{}l p{0.47\textwidth} r l@{}}
\toprule
Const/Func & Description & $\PSLtwoZ$ & Source\\
\midrule
\multicolumn{4}{@{}l}{\textit{Constants} $A_i^{\sharp}$ \textit{(continued)}}\\
$A_{15}^{\sharp}[\kappa]$ & $\tau\bigl[
 \frac{A_{10}^{\sharp}D_1[\kappa]}{2\pi}(1+\frac1\kappa)
 +\frac{2\,\omega(\frac12\kappa^{-1/4})}{\kappa^{1/4}}
 +\frac{2\log 2}{\kappa^{3/4}}\bigr]$ & $9.7310572$ & Prop.~\ref{prop:57newest}\\
$A_{16}^{\sharp}[\kappa]$ & $A_{11}^{\sharp}[\kappa]+A_{12}^{\sharp}[\kappa]
 +A_{13}^{\sharp}[\kappa]+A_{14}^{\sharp}[\kappa]+A_{15}^{\sharp}[\kappa]$
 & $52.74651$ & Prop.~\ref{prop58}\\
$A_{18}^{\sharp}$ & $\max\Bigl\{
 A_{16}^{\sharp}[c_{+}]\,e^{3\delta_{+}/4}+\mathcal{A}\frac{e^{\delta_{+}}-1}{\delta_{+}}
 \;,\;A_{16}^{\sharp}[c]+\mathcal{A}\,e^{c^{-1/4}/4}\;,\;2\mathcal{A}\beta\Bigr\}$;
 the first branch is the active one here
 & $72.03096$ & Prop.~\ref{prop:A18}\\
\midrule
\multicolumn{4}{@{}l}{\textit{Output}}\\
$D_2$ & $\max\Bigl\{f_1(T_1),\ \sup_{m\ge2} \max\bigl(\{f_m(m\ell)\}\cup
 \{f_m(T):T\in W_m\}\bigr)\Bigr\}$ & $181.2996177$ & Prop.~\ref{prop:D2}\\
$M$ & $A_{18}^{\sharp}+D_2$ & $253.3306$ & Prop.~\ref{prop:M}\\
\bottomrule
\end{tabular}
\end{table}

The four-case structure of Theorem~\ref{Theorem2.1} has disappeared: a single
formula, valid without any restriction on the systole, replaces cases
(i)--(iv), and by Remark~\ref{rem:fourcases} it is smaller than each of them.
The output row of Table~\ref{tab4}, and that of Table~\ref{tab:cocompact} in
the cocompact case, record it in the form it takes once the anchor has been
fixed; $R^{\sharp}_{\Gamma}$ depends on $\Gamma$ only through $M$, the anchor
data and $\mathcal{L}$.

\begin{table}[htbp]
\caption{Algorithm for cofinite Fuchsian groups with $\tau\ge1$: self-improvement.
Part 3.  These rows supersede the rows $D_2$ and $M$ of Table~\ref{tab3}: the
constant to be inserted in Theorem~\ref{huber-final} is $M^{\ast}$.  The
numerical column records $\Gamma=\PSLtwoZ$ at $X_0=10^{4}$.}
\label{tab4}\centering
\footnotesize
\setlength{\tabcolsep}{3pt}
\renewcommand{\arraystretch}{1.35}
\begin{tabular}{@{}l p{0.47\textwidth} r l@{}}
\toprule
Const/Func & Description & $\PSLtwoZ$ & Source\\
\midrule
\multicolumn{4}{@{}l}{\textit{Input}}\\
$M_{\mathrm{old}}$ & any constant with $\theta_{\Gamma}(x)\le
 \sum_{\lambda_k\le\frac14}\frac{x^{s_k}}{s_k}+M_{\mathrm{old}}x^{3/4}$,
 $x\ge c$; one starts from $M$ of Table~\ref{tab3}
 & $253.3306$ & Prop.~\ref{prop:M}\\
$X_0$ & enumeration cutoff: $\theta_\Gamma(t)$ known for $c\le t\le X_0$;
 $X_0=c$ recovers the pure bootstrap & $10^{4}$ & Prop.~\ref{prop:boot1}\\
\midrule
\multicolumn{4}{@{}l}{\textit{Majorant of} $\theta_{\Gamma}$}\\
$\Theta(t)$ & $\min\bigl\{(B_1^{\sharp}t-B_2^{\sharp})\log t,\;
 \sum_{\lambda_k\le\frac14}\frac{t^{s_k}}{s_k}+M_{\mathrm{old}}t^{3/4}\bigr\}$
 & --- & Prop.~\ref{prop:boot0}\\
$\Theta_{X_0}(t)$ & $\min\{\theta_\Gamma(t),\Theta(t)\}$ for $t\le X_0$,
 and $\Theta(t)$ for $t>X_0$ & --- & Prop.~\ref{prop:boot1}\\
$\nu_k$ & $\dfrac{c^{k}+1}{c^{k}-1}$ & $\nu_2=1.0434984$ & \eqref{eq:Vboot0}\\
$V^{\flat}(x)$ & $2\Bigl[\frac{\Theta_{X_0}(x)}{x-1}
 +\int_{c}^{x}\frac{\Theta_{X_0}(t)}{(t-1)^{2}}dt\Bigr]
 +\sum_{k=2}^{\lfloor \log x/\ell\rfloor}\nu_k\,\Theta_{X_0}\bigl(x^{1/k}\bigr)$
 & --- & \eqref{eq:Vboot0}\\
\midrule
\multicolumn{4}{@{}l}{\textit{Iteration and output}}\\
$D_2^{\flat}(M_{\mathrm{old}})$ & $\sup_{x\ge c}V^{\flat}(x)x^{-3/4}$; increasing
 in $M_{\mathrm{old}}$ with slope $<1$ & $0.2802273422$ & Prop.~\ref{prop:boot0}\\
$M_{n+1}$ & $A_{18}^{\sharp}+D_2^{\flat}(M_n)$, started at
 $M_0=A_{18}^{\sharp}+D_2$; decreasing, each iterate admissible
 & --- & Prop.~\ref{prop:boot0}\\
$M^{\ast}$ & $\lim_n M_n$, the unique fixed point of
 $M\mapsto A_{18}^{\sharp}+D_2^{\flat}(M)$ & $72.3111906$ & Prop.~\ref{prop:boot0}\\
$c_{0}$ & anchor of Prop.~\ref{prop:anchor}; here the enumerated
 $X_{0}$ of Rem.~\ref{rem:anchordata} & $10^{4}$ & Rem.~\ref{rem:anchordata}\\
$\mathcal{L}(\log c_{0})$ & $\frac{3T^{\natural}}{3T^{\natural}-4}$,
 $T^{\natural}=12.6355334$ & $1.1179711$ & Lem.~\ref{lem:propfactor}\\
$\mathcal{K}(c_{0}),\ \mathcal{Q}(c_{0})$ & anchor correction and
 $\sup_{x>c_{0}}x^{-3/4}\,\log x$ & $-13.3245925$, $0.0092103$
 & Prop.~\ref{prop:anchor}\\
$\mathcal{E}_{0}(c_{0})$ & the range $1<x<c_{0}$, from the enumerated
 spectrum & $2.0854378$ & Prop.~\ref{prop:anchor}\\
$R_{\Gamma}^{\sharp}$ & $\max\{\mathcal{E}_{0},\,
 \mathcal{L}M^{\ast}+|\mathcal{K}|\mathcal{Q}\}$
 & $80.96455$ & Thm.~\ref{huber-final}\\
\bottomrule
\end{tabular}
\end{table}

\section{Application to \PSLtwoZ}
\label{sec:psl(2,Z)}
\subsection{Spectral counting on the modular group}
\label{sec:modular}

For $\Gamma=\PSLtwoZ$, we take $\mathcal F=\Bigl\{z=x+iy:\  |x|\le\tfrac12,\ |z|\ge1\Bigr\}$. The fundamental constants are
\begin{equation}
 \begin{gathered}
 |\mathcal F|=\frac\pi3,\qquad \tau=1,\qquad \mathcal A=1,\qquad
 \mathcal N=1,\quad \sigma_1=1,\qquad c_1=1,\\[4pt]
 c=\frac{7+3\sqrt5}{2},\qquad \ell=\log c=1.9248473002\ldots,
 \end{gathered}
 \label{eq:modular-data}
\end{equation}
so that $|\log c_1|=0$ and $\ell^{2}/8=0.4631296\ldots$.  There are three
non-identity elliptic conjugacy classes, with
$(\theta,m)=(\frac\pi2,2),(\frac\pi3,3),(\frac{2\pi}3,3)$, whence
\begin{equation}
 E=\frac18+\frac19+\frac19=\frac{25}{72}=0.3472222\ldots,
 \qquad
 \frac14\operatorname{tr}\Bigl(I-\Phi\bigl(\tfrac12\bigr)\Bigr)=\frac12 ,
 \label{eq:modular-E}
\end{equation}
the latter because $\phi(\frac12)=-1$.  

Let us calculate $B_1^{\sharp}$ and $B_2^{\sharp}$ on $\Gamma=\PSLtwoZ$. Put
\[
y_*:=\frac{\sqrt3}{2}=\min_{z\in\mathcal F}\operatorname{Im}z,
  \qquad
  p_{\pm}:=\pm\tfrac12+iy_* ,
\]
and let $\dist$ denote the hyperbolic distance,
\begin{equation}
  \cosh \dist(z,w)=1+\frac{|z-w|^{2}}{2\operatorname{Im}z\operatorname{Im}w},
  \qquad z,w\in\mathbb H .
  \label{eq:cosh-rho}
\end{equation}
For $Y\ge 1$ we write $\mathcal F_{Y}:=\{z\in\mathcal F:\operatorname{Im}z\le Y\}$.
Recall from \cite{FriedmanJorgensonKramer2011} that the trivial counting lemma applies to any
measurable $K$ with
\begin{equation}
  \mathcal F_{Y}\subseteq K\subseteq\mathcal F
  \qquad\text{for some admissible }Y\ge1,
  \label{eq:admissible}
\end{equation}
and yields, with $|K|$ the hyperbolic area and $d(K)$ the hyperbolic
diameter of $K$,
\begin{equation}
  \pi_\Gamma(x)\le\frac{4\pi}{|K|}\bigl(\cosh(\log x+2d(K))-1\bigr)
  \qquad (x>1).
  \label{eq:counting}
\end{equation}
Since the quantity to be minimised is, up to the lower order terms,
$e^{2d(K)}/|K|$, both a large area and a small diameter are wanted, and
the two requirements pull in opposite directions.  We determine the
exact optimum.

Observe first that $\varepsilon=\dist_{\mathrm{Eucl}}
(\mathcal F,\mathbb R)=y_*$, so that the lower cut-off
$\varepsilon\le\operatorname{Im}z$ occurring in the definition of the
sets $K_j$ in \cite{FriedmanJorgensonKramer2011} is vacuous for $\Gamma=\PSLtwoZ$: every
point of $\mathcal F$ satisfies it.

\subsection*{The horizontal truncation}

We first record the optimum inside the one--parameter family of
horizontal truncations $K(Y)=\{z\in\mathcal F:\operatorname{Im}z\le
Y+1\}$; this is the family implicitly used in \cite{FriedmanJorgensonKramer2011} and \cite{Avdispahic2021}.

\begin{lemma}[Optimal cusp height]
\label{lem:modular-optimal-Y}
Put $b:=Y+1\ge 2$ and $K(Y)=\{z\in\mathcal F:\operatorname{Im}z\le b\}$.
Then
\begin{equation}
  |K(Y)|=\frac{\pi}{3}-\frac1b ,
  \label{eq:modular-KY-area}
\end{equation}
and, with
\[
  b_\ast:=\frac{7\sqrt3}{6}=\frac{7}{2\sqrt3},
  \qquad
  Y_\ast:=b_\ast-1,
\]
the hyperbolic diameter equals
\begin{equation}
  d(Y)=
  \begin{cases}
    \log3, & 2\le b\le b_\ast,\\[2mm]
    \displaystyle \arcosh\!\left(\frac{b^{2}+\frac74}{\sqrt3\,b}\right),
    & b\ge b_\ast .
  \end{cases}
  \label{eq:modular-dY}
\end{equation}
For every $x>1$ the majorant
$\Pi_Y(x):=\frac{4\pi}{|K(Y)|}\bigl(\cosh(\log x+2d(Y))-1\bigr)$
is minimised, among all admissible $Y\ge1$, at
\begin{equation*}
  Y=Y_\ast=\frac{7\sqrt3}{6}-1=1.0207259421\ldots
\end{equation*}
\end{lemma}

\begin{proof}
Admissibility and \eqref{eq:modular-KY-area} are immediate, the part of
$\mathcal F$ above height $b$ being a strip of width $1$ and area $1/b$.

For the diameter, let $z=x+iy$, $w=u+iv\in K(Y)$.  Since
$|x|,|u|\le\frac12$ we have $(x-u)^{2}\le1$, with equality attained on
the two vertical sides of $\mathcal F$, which are present at every
height in $[y_*,b]$. Hence, by \eqref{eq:cosh-rho}
\begin{equation}
  \max_{z,w\in K(Y)}\cosh \dist(z,w)
  =\max_{(y,v)\in[y_*,b]^{2}}G(y,v),
  \qquad
  G(y,v):=\frac{1+y^{2}+v^{2}}{2yv}.
  \label{eq:G-def}
\end{equation}
For fixed $v$ one has
$\partial_y G=(y^{2}-1-v^{2})/(2y^{2}v)$, so $y\mapsto G(y,v)$ decreases
on $[y_*,\sqrt{1+v^{2}}\,]$ and increases afterwards; the unique interior
critical point is a minimum.  The same holds in $v$.  Consequently the
maximum in \eqref{eq:G-def} is attained at a corner of $[y_*,b]^{2}$.  The
three corner values are
\[
  G(y_*,y_*)=\frac53,
  \qquad
  G(y_*,b)=G(b,y_*)=\frac{b^{2}+\frac74}{\sqrt3\,b},
  \qquad
  G(b,b)=1+\frac1{2b^{2}}<\frac53 ,
\]
the last because $b\ge2$.  Moreover $G(y_*,b)\le G(y_*,y_*)$ precisely when
$3b^{2}-5\sqrt3\,b+\frac{21}{4}\le0$, i.e.\ for
$b\in[\frac{\sqrt3}{2},\frac{7\sqrt3}{6}]$.  This proves
\eqref{eq:modular-dY}; note that both values are attained, at the pairs
$(p_+,p_-)$ and $(p_+,-\frac12+ib)$ respectively, so
\eqref{eq:modular-dY} is an equality and not merely an upper bound.

It remains to optimise.  On $2\le b\le b_\ast$ the diameter is constant
while $A(b):=\pi/3-1/b$ strictly increases, so $\Pi_Y(x)$ strictly
decreases there.  For $b\ge b_\ast$ put $Z:=\log x+2d(b)$.  Then
\[
  \frac{\partial}{\partial b}\Pi_Y(x)
  =\frac{4\pi}{A(b)}
  \left[2d'(b)\sinh Z-\frac{A'(b)}{A(b)}(\cosh Z-1)\right],
\]
with
\[
  d'(b)=\frac{4b^{2}-7}{b\sqrt{16b^{4}+8b^{2}+49}},
  \qquad
  \frac{A'(b)}{A(b)}=\frac{3}{b(\pi b-3)} .
\]
The function $L(b):=2b\,d'(b)=2(4b^{2}-7)/\sqrt{16b^{4}+8b^{2}+49}$ is
strictly increasing, since
$L'(b)=128b(4b^{2}+7)(16b^{4}+8b^{2}+49)^{-3/2}>0$, and $L(b_\ast)=1$,
while $3/(\pi b-3)$ strictly decreases and
$3/(\pi b_\ast-3)=0.8960\ldots<1$.  Hence $2d'(b)>A'(b)/A(b)$ for
$b\ge b_\ast$, and since $\sinh Z/(\cosh Z-1)=\coth(Z/2)>1$ we get
$\partial_b\Pi_Y(x)>0$ there.  Thus $b_\ast$ is the unique minimiser,
for every $x>1$ simultaneously.
\end{proof}

At $Y=Y_\ast$ one has
\begin{equation*}
  |K(Y_\ast)|=\frac{\pi}{3}-\frac{2\sqrt3}{7}
  =\frac{7\pi-6\sqrt3}{21}=0.5523258919\ldots,
  \qquad d=\log3,
\end{equation*}
whence $\pi_\Gamma(x)\le B_+x+B_-x^{-1}-B_0$ with
\begin{equation*}
\begin{gathered}
  B_+=\frac{378\pi}{7\pi-6\sqrt3}=102.3827935551\ldots,
  \quad
  B_-=\frac{14\pi}{3(7\pi-6\sqrt3)}=1.2639851056\ldots,
  \\[4pt]
  B_0=\frac{84\pi}{7\pi-6\sqrt3}=22.7517319011\ldots
\end{gathered}
\end{equation*}
In particular $\pi_\Gamma(x)<102.383\,x-22.567$ for $x\ge c$.

\subsection*{The optimal kernel}

The family of horizontal truncations is not optimal: the top corners of
$K(Y_\ast)$ may be traded for additional height over the middle of
$\mathcal F$ without increasing the diameter.  Set
\begin{equation*}
  K^{\star}:=\bigl\{z\in\mathcal F:\
  \dist(z,p_+)\le\log3\ \text{ and }\ \dist(z,p_-)\le\log3\bigr\}.
\end{equation*}
By \eqref{eq:cosh-rho}, $\dist(x+iy,p_{\mp})\le\log3$ is equivalent to
$\bigl(|x|+\frac12\bigr)^{2}+(y-y_*)^{2}\le\frac{4y_*}{3}y$, so that
\begin{equation}
  K^{\star}=\Bigl\{x+iy\in\mathcal F:\ |x|\le W(y)\Bigr\},
  \qquad
  W(y):=\min\Bigl\{\tfrac12,\;R(y)-\tfrac12\Bigr\},
  \label{eq:Kstar-width}
\end{equation}
where
\begin{equation}
  R(y):=\sqrt{-y^{2}+\tfrac{5}{\sqrt3}\,y-\tfrac34}
  =\sqrt{\tfrac{4}{3}-\bigl(y-\tfrac{5}{2\sqrt3}\bigr)^{2}}\; .
  \label{eq:R-def}
\end{equation}
The equation $R(y)=1$ has the two roots $y=y_*$ and $y=b_\ast$, and
$R(y)=\frac12$ has the roots $y=(5\pm\sqrt{13})/(2\sqrt3)$.  Writing
\begin{equation*}
  h_0:=\frac{5+\sqrt{13}}{2\sqrt3}=2.4842086727\ldots,
\end{equation*}
we therefore have $W\equiv\frac12$ on $[y_*,b_\ast]$, $W=R-\frac12$
decreasing on $[b_\ast,h_0]$, and $W(h_0)=0$: the region $K^{\star}$ is
the fundamental domain truncated by the circular arc
$\bigl(|x|+\frac12\bigr)^{2}+(y-y_*)^{2}=\frac{4y_*}{3}y$, which rises from
height $b_\ast$ at $|x|=\frac12$ to height $h_0$ at $x=0$.

\begin{lemma}[The maximal kernel of diameter $\log3$]
\label{lem:Kstar}
The set $K^{\star}$ satisfies $\mathcal F_{2}\subset K^{\star}\subset\mathcal F$,
\begin{equation*}
  d(K^{\star})=\log3,
\end{equation*}
and
\begin{equation}
  |K^{\star}|
  =\frac{2\pi}{3}
  +\frac{10}{3}\left(
     \arcsin\frac{5+3\sqrt{13}}{16}-\arcsin\frac{13}{14}\right)
  -2\arcsin\frac{\sqrt{13}}{4}
  =0.6108620738\ldots
  \label{eq:Kstar-area}
\end{equation}
Moreover, if $K$ satisfies \eqref{eq:admissible} and $d(K)=\log3$, then
$K\subseteq K^{\star}$.  In particular $K^{\star}$ is the unique maximal
admissible kernel of diameter $\log 3$.
\end{lemma}

\begin{proof}
\emph{Inclusions.}  $K^{\star}\subset\mathcal F$ by definition, and
$W(y)=\frac12$ for $y\le b_\ast$ gives
$\mathcal F_{b_\ast}\subset K^{\star}$; since $b_\ast>2$ we get
$\mathcal F_{2}\subset K^{\star}$, so \eqref{eq:admissible} holds with
$Y=2$ (and indeed $K^{\star}\supset K(1)$).

\emph{Diameter.}  Let $z=x+iy$, $w=u+iv\in K^{\star}$.  By
\eqref{eq:Kstar-width}, $|x|\le W(y)$ and $|u|\le W(v)$, hence
$(x-u)^{2}\le\bigl(W(y)+W(v)\bigr)^{2}$ and, by \eqref{eq:cosh-rho},
$\cosh \dist(z,w)\le\frac53$ will follow from
\begin{equation}
  \bigl(W(y)+W(v)\bigr)^{2}+(y-v)^{2}\le\tfrac43\,yv,
  \qquad y,v\in[y_*,h_0].
  \label{eq:key}
\end{equation}
We verify \eqref{eq:key} in three cases.

\emph{(i) $y,v\le b_\ast$.}  Here \eqref{eq:key} reads
$1+y^{2}+v^{2}\le\frac{10}{3}yv$, i.e.\ $G(y,v)\le\frac53$, which was
proved in Lemma \ref{lem:modular-optimal-Y}.

\emph{(ii) $y\le b_\ast\le v$.}  Now $W(y)+W(v)=R(v)$, and by
\eqref{eq:R-def},
\[
  R(v)^{2}+(y-v)^{2}-\tfrac43yv
  =\tfrac{10}{3}y_*v-y_*^{2}+y^{2}-\tfrac{10}{3}yv
  =-(y-y_*)\Bigl(\tfrac{10}{3}v-y_*-y\Bigr).
\]
Since $y\ge y_*$ and $\frac{10}{3}v\ge\frac{10}{3}b_\ast>y_*+b_\ast\ge y_*+y$,
the right-hand side is $\le0$, which is \eqref{eq:key}.  Equality holds
exactly for $y=y_*$: the whole upper arc of $\partial K^{\star}$ lies at
distance exactly $\log 3$ from the opposite corner $p_\mp$.

\emph{(iii) $y,v\ge b_\ast$.}  By \eqref{eq:R-def} the function $R$
attains its maximum at $y=\frac{5}{2\sqrt3}<b_\ast$ and is decreasing on
$[b_\ast,h_0]$, so $R\le R(b_\ast)=1$ there and
$W(y)+W(v)=R(y)+R(v)-1\le1$.  Together with
$(y-v)^{2}\le(h_0-b_\ast)^{2}<0.22$ this gives
$\bigl(W(y)+W(v)\bigr)^{2}+(y-v)^{2}<1.22
<\frac43b_\ast^{2}\le\frac43yv$.

Thus $d(K^{\star})\le\log3$; since $p_\pm\in K^{\star}$ and
$\dist(p_+,p_-)=\arcosh\frac53=\log3$, equality holds.

\emph{Area.}  At height $y\in[b_\ast,h_0]$ the horizontal section of
$\mathcal F\setminus K^{\star}$ has Euclidean width $2\bigl(1-R(y)\bigr)$,
and $\mathcal F\setminus K^{\star}$ contains the whole part of $\mathcal F$ above $h_0$; hence
\[
  |K^{\star}|=\frac{\pi}{3}-\frac1{h_0}
  -2\int_{b_\ast}^{h_0}\frac{1-R(y)}{y^{2}}\,dy .
\]
Using
\[
  \int\frac{R(y)}{y^{2}}\,dy
  =-\frac{R(y)}{y}
  +\frac53\arcsin\frac{10y-3\sqrt3}{8y}
  -\arcsin\frac{2\sqrt3\,y-5}{4},
\]
together with $R(b_\ast)=1$, $R(h_0)=\frac12$ and
$h_0^{-1}=2\sqrt3/(5+\sqrt{13})$, a direct evaluation gives
\eqref{eq:Kstar-area}.

\emph{Maximality.}  If $K$ is admissible then $p_\pm\in\mathcal F_{Y}
\subseteq K$, so every $z\in K$ satisfies
$\dist(z,p_\pm)\le d(K)=\log3$, i.e.\ $z\in K^{\star}$.
\end{proof}

\begin{proposition}
\label{prop:modular-final}
For $\Gamma=\PSLtwoZ$ and all $x>1$,
\begin{equation*}
  \pi_\Gamma(x)\le\Pi^{\ast}(x)
  :=\frac{4\pi}{|K^{\star}|}\bigl(\cosh(\log x+2\log3)-1\bigr)
  =B_+^{\ast}x+B_-^{\ast}x^{-1}-B_0^{\ast},
\end{equation*}
with
\begin{equation}
\begin{gathered}
  B_1^{\sharp}=B_+^{\ast}=\frac{18\pi}{|K^{\star}|}=92.5719081051\ldots,
  \quad
  B_-^{\ast}=\frac{2\pi}{9|K^{\star}|}=1.1428630630\ldots,
  \\
  B_0^{\ast}=\frac{4\pi}{|K^{\star}|}=20.5715351345\ldots, \text{hence} \quad B_2^{\sharp}=20.4047936607\ldots
  \label{eq:final-constantsnew}
\end{gathered}
\end{equation}
In particular, 
\begin{equation*}
  \pi_\Gamma(x)\le B_1^{\sharp}x- B_2^{\sharp}
  <92.572\,x-20.404
  \qquad (x\ge c).
\end{equation*}
\end{proposition}

\begin{proof}
Combine \eqref{eq:counting} with Lemma \ref{lem:Kstar} and
$\cosh(\log x+2\log3)=\frac12(9x+\frac1{9x})$.
\end{proof}

\begin{remark}[Optimality]
\label{rem:optimality}
Every admissible $K$ contains $p_\pm$, so $d(K)\ge\log3$, and by Lemma
\ref{lem:Kstar} the constant $e^{2d(K)}/|K|$ cannot be improved as long
as $d(K)=\log3$.  Allowing a larger diameter does not help either: if
$d(K)=\vartheta>\log3$, then $K\subseteq\check K_\vartheta:=\{z\in\mathcal F:\ 
\dist(z,p_\pm)\le\vartheta\}$, and $\vartheta\mapsto e^{2\vartheta}/|\check K_\vartheta|$
increases,
\[
  \begin{array}{l|ccccc}
    \vartheta & \log 3 & 1.1588 & 1.2257 & 1.3170 & 1.4254\\\hline
    |\check K_\vartheta| & 0.6109 & 0.6448 & 0.6775 & 0.7158 & 0.7546\\
    e^{2\vartheta}/|\check K_\vartheta| & 14.733 & 15.743 & 17.129 & 19.457 & 22.930
  \end{array}
\]
For $\vartheta\ge\frac12\log\bigl(\tfrac{\pi}{3}\cdot\tfrac{9}{|K^{\star}|}\bigr)
=1.3681\ldots$ this is immediate from the trivial bound
$|K|\le|\mathcal F|=\pi/3$.  Hence $K^{\star}$ is optimal and
\eqref{eq:final-constantsnew} is the best constant obtainable from
\eqref{eq:counting}.
\end{remark}

\begin{remark}[Comparison with the literature]
\label{rem:comparison}
In \cite[Section~3.2]{FriedmanJorgensonKramer2011} and \cite[Theorem~4.3]{Avdispahic2021}, the kernel is taken to be
$K_1=\{z\in\mathcal F:\frac32\le\operatorname{Im}z\le3\}$, with
$|K_1|=\frac13$ and $e^{2d_1}<9.761757$, giving
$B_1=\frac{2\pi}{|K_1|}e^{2d_1}=184.0048$.  Since
$\min_{\mathcal F}\operatorname{Im}z=\frac{\sqrt3}{2}$, the lower
cut-off should read $\varepsilon=\frac{\sqrt3}{2}$, for which the same
cusp height $Y=2$ gives
$|K_1|=\frac\pi3-\frac13=0.713864$,
$d_1=\arcosh\frac{43}{12\sqrt3}=1.355820$ and $B_1=132.4996$.  All
constants of the algorithm remain valid, since they are increasing in
$B_1$.  Optimising the cusp height lowers this to
$B_1=B_+=102.3828$ (Lemma \ref{lem:modular-optimal-Y}), and replacing
the horizontal truncation by $K^{\star}$ to
$B_1=B_+^{\ast}=92.5719$ (Proposition \ref{prop:modular-final}),
a factor $1.988$ below the published value.  
\end{remark}

\begin{theorem}
\label{cor:modular}
Let $\Gamma=\PSLtwoZ$.  Then, for all $R\ge0$,
\begin{equation}
 \mu(R)\le 2.4657613\,R^{2}-1.0317867\,R+2.1177623
 \label{eq:modular-final-lin}
\end{equation}
and also
\begin{equation}
 \mu(R)\le 2.4657613\,R^{2}+1.6018689 ,
 \,\,\text{equivalently}\,\,
 \mu(R)\le 2.4657613\Bigl(R^{2}+\frac14\Bigr)+0.9854286 .
 \label{eq:modular-final-const}
\end{equation}
\end{theorem}

\begin{proof}
Take $\eta=0.035$, admissible since $0.035<\ell^{2}/8$, and observe that
the argument of $\operatorname{erfc}$ in $I_{1,+}(\eta,\frac12)$ equals $2.198$.  With
\eqref{eq:modular-data}, \eqref{eq:modular-E} and \eqref{eq:final-constantsnew} one computes
\[
 \rho(0.14)=0.4551946,\,\,
 b(\eta)=-0.9962988,\,\,
 \mathcal H_{\max}^{\mathrm G}(\eta)=0.6196596,\,\,
 C(\eta)=1.9668819,
\]
the series \eqref{eq:gauss-Hmax-def} being dominated by its first term
($b_1\mathcal G_1^{\sharp}=0.6196596$, $b_2\mathcal G_2^{\sharp}<5\cdot10^{-10}$).
Hence $K_0=2.4657613$.  Since $b(\eta)<0$ and $c_1=1$, the choice
$R_0=\infty$, i.e.\ $(\lambda,\xi)=(1,0)$, gives
$K_1=-1.0317867$, $K_2=2.1177623$ and \eqref{eq:modular-final-lin},
while $R_0=0$, i.e.\ $(\lambda,\xi)=(0,\frac12)$, gives $K_1=0$,
$K_2=1.6018689$ and \eqref{eq:modular-final-const}.
\end{proof}

\begin{remark}[Choice of $\eta$ and of $R_0$]
\label{rem:modular-choices}
Both parameters admit a genuine optimisation.  The dependence on $\eta$
is a trade-off: $K_0(\eta)$ decreases in $\eta$, whereas $\mathcal H_{\max}^{\mathrm G}(\eta)$,
and hence $K_2(\eta)$, grows rapidly.  With $K^{\star}$ and
$(\lambda,\xi)=(1,0)$:
\[
 \begin{array}{c|cccccc}
 \eta & 0.02 & 0.03 & 0.035 & 0.04 & 0.05 & 0.10\\\hline
 \mathcal H_{\max}^{\mathrm G}(\eta) & 0.00538 & 0.21750 & 0.61966 & 1.35392 & 4.02255 & 35.6738\\
 K_0(\eta) & 4.250839 & 2.862374 & 2.465761 & 2.168356 & 1.752118 & 0.920976\\
 K_1(\eta) & -1.509293 & -1.147617 & -1.031787 & -0.940065 & -0.803206 & -0.488346\\
 K_2(\eta) & 1.922436 & 1.797516 & 2.117762 & 2.812657 & 5.531841 & 40.539619
 \end{array}
\]
The choice, e. g., $\eta=\ell/8=0.2406059$ (the value for which the argument of $\operatorname{erfc}$ already in $I_{1,+}(\eta,\frac12)$ is
positive) yields
$K_0=0.4405612$, $K_1=-0.2657859$, $K_2=209.3631$, which is inferior to
\eqref{eq:modular-final-lin} for all $R\lesssim10^{2}$ and only becomes
preferable for large $R$.  Similarly, at $\eta=0.035$ the choice of
$R_0$  interpolates between
\eqref{eq:modular-final-const} and \eqref{eq:modular-final-lin}:
\[
 \begin{array}{c|ccccc}
 R_0 & 0 & 1 & 2 & 3 & \infty\\\hline
 \lambda & 0 & 0.894427 & 0.970143 & 0.986394 & 1\\
 \xi & 0.5 & 0.223607 & 0.121268 & 0.082199 & 0\\
 K_1 & 0 & -0.922858 & -1.000980 & -1.017748 & -1.031787\\
 K_2 & 1.601869 & 1.887048 & 1.992640 & 2.032950 & 2.117762
 \end{array}
\]
\end{remark}

\begin{remark}[Comparison]
\label{rem:modular-comparison}
The leading constant in Theorem~\ref{cor:modular} is
$K_0=2.4657613=29.589\cdot\frac{|\mathcal F|}{4\pi}$, i.e.\ it exceeds the
Weyl constant $\frac1{12}$ by a factor $29.59$.  It is smaller by a
factor $1.63\cdot10^{4}$ than the constant $40\,313.68$ of \cite[Th. 4.3.]{Avdispahic2021}, and by a factor $6.83\cdot10^{5}$ than the constant
$1\,682\,997$ of \cite[Th. 3.4.]{FriedmanJorgensonKramer2011}.  Two mechanisms are responsible, and
they are independent:
\begin{itemize}
\item the replacement of the Tauberian step by the direct positive
 majorant, which is what removes the leading constant's inflation;
\item the optimization of the geometric kernel, which in the present
 framework affects only $\mathcal H_{\max}^{\mathrm G}$, hence only $K_2$.  At, e. g., 
 $\eta=\ell/8$ the successive kernels give
\[
 \begin{array}{l|cc}
 \text{kernel} & B_1^{\sharp} & K_2\\\hline
 |K_1|=\tfrac13,\ e^{2d_1}=9.761757\ \text{(published)} & 184.0048 & 415.499\\
 \varepsilon=\tfrac{\sqrt3}{2},\ Y=2 & 132.4996 & 300.651\\
 \varepsilon=\tfrac{\sqrt3}{2},\ Y=Y_{\ast} & 102.3828 & 231.440\\
 K^{\star} & 92.5719 & 209.363
 \end{array}
\]
 while $K_0$ and $K_1$ are unchanged.  (In the Tauberian algorithm, by
 contrast, $B_1$ propagates into the leading constant.)
\end{itemize}
\end{remark}

\subsection{Huber constant on the modular group}
\label{sec:modulartwo}

The elliptic contribution in the corresponding Selberg trace formula, the constant
$E^{\sharp}$ of Table~\ref{tab2}, is here
\[
  E^{\sharp}=\sum_{R}\frac{1}{2m_R\sin\theta_R}
  \sqrt{\frac{\pi-2\theta_R}{\pi\sin(2\theta_R)}}
  =\frac14\sqrt{\frac1\pi}+\frac{2}{3\sqrt3}\sqrt{\frac{2}{3\sqrt3}}
  =0.3798407141\ldots
\]
 We recall
\[
  B_1^{\sharp}=92.5719081051\ldots,\qquad B_2^{\sharp}=20.4047936607\ldots,
  \qquad \varkappa=B_2^{\sharp}/B_1^{\sharp}=0.2204\ldots
\]
Since $\mathcal{A}=1$, $A_{11}^{\sharp}[\kappa]=1+A_8^{\sharp}\kappa^{-1/4}$ and $s_1$ does not
occur.

The thresholds are
\begin{equation*}
 \begin{gathered}
  c_{+}=11.761542277\ldots,\,\, \delta_{+}=c_{+}^{-1/4}=0.5399877796\ldots,
  \,\, T_{\star}=2.5428812890\ldots,\\
 \beta=1.8883818950\ldots
\end{gathered}
\end{equation*}
Both $c$ and $c_{+}$ lie below $\kappa_{\ast}=47.1895$, so by \eqref{eq:D1}
\[
  D_1[c]=D_1[c_{+}]=0.7786285303\ldots
\]

As far as the constants $A_i^{\sharp}$ are concerned, $A_8^{\sharp}=0.1036383$, $A_9^{\sharp}=6.4214$, $A_{10}^{\sharp}=57.4347476$. For $\Gamma=\PSLtwoZ$ one has
$\Lambda\bigl(\tfrac14\log c\bigr)=0.9484592\ldots$ and
$\Lambda\bigl(\tfrac14\log c_{+}\bigr)=1.1860291\ldots$, while the constants
of \eqref{eq:modular-final-const} give $K_2/K_0=0.6496\ldots$. Hence, we may apply Corollary 
\ref{cor:floor}
\[
 A_{12}^{\sharp}[c]=A_{12}^{\sharp}[c_{+}]=A_{10}^{\sharp}K_0=141.620378\ldots
\]
Carrying this through
\begin{center}
\begin{tabular}{lrr}
\toprule
 & $\kappa=c$ & $\kappa=c_{+}$\\
\midrule
$A_{11}^{\sharp}[\kappa]$ & $1.0640520$ & $1.0559634$\\
$A_{13}^{\sharp}[\kappa]$ & $5.4845304$ & $5.1931679$\\
$A_{14}^{\sharp}[\kappa]$ & $4.1457796$ & $3.5952482$\\
$A_{15}^{\sharp}[\kappa]$ & $9.7310572$ & $9.0287545$\\
\bottomrule
\end{tabular}
\end{center}
gives
\[
 A_{16}^{\sharp}[c]=162.0458,\quad A_{16}^{\sharp}[c_{+}]=160.4935,\quad
 A_{18}^{\sharp}=241.9520.
\]
$D_2$ depends only on $c$, $B_1^{\sharp}$, $B_2^{\sharp}$ and is fully determined.
One finds $\widetilde T_1\in(\ell,2\ell)$ and
\[
  T_1=2.9052259382\ldots,\qquad f_1(T_1)=127.3350259\ldots,
\]
while on each interval with $m\ge2$ the function $f_m$ is decreasing, so
$W_m=\emptyset$ and only the left endpoints matter:
\begin{center}
\begin{tabular}{rrrr}
\toprule
$m$ & $m\ell$ & $f_m(m\ell)$ & $\Psi(m\ell)$\\
\midrule
2 & $3.84969$ & $181.2996177$ & $256.894$\\
3 & $5.77454$ & $139.1439195$ & $180.153$\\
4 & $7.69939$ & $91.0051540$ & $119.491$\\
5 & $9.62424$ & $58.4333947$ & $78.899$\\
6 & $11.54908$ & $38.2185682$ & $52.163$\\
\bottomrule
\end{tabular}
\end{center}
Since $\Psi$ is decreasing on $[3\ell,\infty)$ and
$\Psi(3\ell)=180.153<181.2996$, no $m\ge3$ can contribute. Therefore
\begin{equation*}
  D_2=f_2(2\ell)=181.2996177\ldots
\end{equation*}
attained at $x=c^{2}=46.97871\ldots$
\begin{equation*}
  M=A_{18}^{\sharp}+D_2=241.9520+181.2996177\ldots=423.2516177
  \dots
\end{equation*}

The propagation of $M$ to $\pi_{\Gamma}$ will be carried out once and for all
in \S\ref{sec:modulartwo} below, with the sharp factor of
Lemma~\ref{lem:propfactor}; for the purpose of recording the successive
values of $M$ we use here the superseded form, i.e.\ case~(ii) of
Theorem~\ref{Theorem2.1}, which applies because
$\frac43<\ell=1.9248473\le\frac{2(1+\sqrt5)}{3}=2.1573790$:
\[
  \frac{4}{3e}\,\frac{\mathcal{A}(c-2)}{\log2}
  +M\Bigl(2+\frac{4}{3\log c}\Bigr)
  =3.4350075+2.6926956\,M ,
\]
whence $C_{\mathrm{Hub},\PSLtwoZ}<1144$ at this stage.

The joint supremum approach has led us so far. 

\subsection*{Two-piece form}

Let us see what Proposition \ref{prop:twopiece} brings. Take $Q_A$ to be the pointwise envelope of Proposition \ref{unlinearised}
\[
 Q_A(r):=\inf_{\eta>0}\mathcal P(r,\eta),
\]
which is admissible in Proposition \ref{prop:twopiece} because that proposition requires
only $\mu\le Q_A$ pointwise. The minimising $\eta$ increases from
$0.0246$ at $r=0$ to $0.719$ at $r=3$; it is the freedom to let $\eta$
grow with $L$, and hence to keep $H^{G}(\eta,L)$ small, that makes
$Q_A$ useful here, and this is exactly what the frozen value
$\mathcal H^{\mathrm G}_{\max}(\eta)=H^{G}(\eta,\tfrac12)$ of
Theorem \ref{thm:direct-gaussian-K} forbids.  

Now, take $Q_B$ from Proposition~\ref{prop:local}, with
\[
 \eta=\frac{\ell^{2}}{8}=0.4631296\ldots,\qquad R_0=3,\qquad
 L_0=\sqrt{9.25}=3.0413813\ldots
\]
Then $b(\eta)=-0.1166936\ldots$, $H^G(\eta,L_0)=6.50\cdot10^{-6}$ and
\[
 K_0^{\infty}=0.2859241,\qquad K_1^{\infty}=-0.1829078,\qquad
 K_2^{\infty}=1.1079988 .
\]

One has $Q_A\le Q_B$ on $[0,0.10502]$ and on
$[2.18305,R_0]$, so the integrand of \eqref{eq:corr} is supported on the
intermediate interval, and
\[
 \mathcal E(\varepsilon,R_0)=0.0693634 ,\qquad
 A_{12}^{\sharp}[c]=32.3211,\qquad A_{12}^{\sharp}[c_+]=28.2859 ,
\]
the supremum in \eqref{eq:A12two1} being attained at the left endpoint
$u=\tfrac14\log\kappa$ in both cases.

Table \ref{tab3} then gives
\[
 A_{16}^{\sharp}[c]=52.74651,\quad
 A_{16}^{\sharp}[c_+]=47.15905,\quad
 A_{18}^{\sharp}=72.03096,\quad
 M=253.3306,
\]
whence $C_{\mathrm{Hub},\PSLtwoZ}<686$.

The choice $\eta=\ell^{2}/8$ is essentially optimal. A two-dimensional search
over $(\eta,R_0)$ gives $R^{\sharp}_{\Gamma}=685.13$ at
$(\eta,R_0)\approx(0.60,3.2)$, so a convenient choice
$\eta=\ell^{2}/8$ 
loses $0.3\%$. The reason not to go higher comes from the fact that $e^{\eta}$ multiplies the structural constant
$E+\tfrac14\operatorname{tr}(I-\Phi(\tfrac12))+\tfrac{\mathcal N}{2}-\mathcal {A'}=0.347$ in
$K^{\infty}_2$. At $\eta=1$ (which is admissible in Proposition \ref{prop:local} but not in Theorem \ref{thm:direct-gaussian-K}), the gain on $K^{\infty}_0$ ($0.2859\to0.2265$) is worth $2.3$ in
$A^{\sharp}_{12}$, while the loss on $K^{\infty}_2$ costs $7.7$. Notice further that the surface is extremely flat: every
$(\eta,R_0)$ with $\eta\in[0.3,0.8]$ and $R_0\ge2.5$ gives
$R^{\sharp}_{\Gamma}\in[685,697]$. Taking $R_0$ beyond $3$ gains nothing, since $H^G(\eta,L_0)$ is already below $10^{-5}$ there.

The constant $M$ is dominated by $D_2=181.2996$, which contributes
$2.6927\cdot D_2=488.2$ of the bound $686$. The source of the inflation is the
use of $\pi_{\Gamma}(x)\le B_1^{\sharp}x-B_2^{\sharp}$ in \eqref{eq:Phi}: at $x=c^{2}$ one has
$\Pi(c^{2})\log c^{2}=16663$ whereas $\theta_{\Gamma}(c^{2})=31.7$

\subsection*{Bootstrap}

Concerning self-improvement, Proposition \ref{prop:M} with $\mathcal{A}=1$, 
$s_0=1$ gives $\theta_{\Gamma}(x)\le x+Mx^{3/4}$ for $x\ge c$, and, as already said, this may be
fed back into \eqref{eq:decomp}, in which only $\theta_{\Gamma}$ occurs.

For $\Gamma=\PSLtwoZ$, with
$A_{18}^{\sharp}=72.03096$:
\[
 \begin{array}{r|ccccccc}
  n & 0 & 1 & 2 & 3 & 4 & 5 & \to\infty\\\hline
  M_n & 253.331 & 177.148 & 145.765 & 132.837 & 127.511 & 125.317 & 123.7805
 \end{array}
\]
The supremum in $D_2^{\flat}$ is attained at $x=c^{2}$ throughout. Thus
\[
 D_2^{\flat}=51.7495,\qquad M=123.7805,\qquad
 R^{\sharp}_{\Gamma}=3.4350+2.6927\,M=336.74
\]
in the same superseded form, so that
\begin{equation}\label{eq:final}
 C_{\mathrm{Hub},\PSLtwoZ}<337.
\end{equation}
\begin{remark}[Where the remaining weight sits]\label{rem:floor}
Three independent blocks make up $M=A_{18}^{\sharp}+D_2$.

\emph{(i) The spectral block.} By Proposition~\ref{prop:joint} and the proof of Corollary \ref{cor:floor},\\
$A_{12}^{\sharp}[\kappa]\ge A_{10}^{\sharp}\lim_{u\to\infty}K_0F(u)=A_{10}^{\sharp}K_0$ for any
majorant of $\mu$ of exact order $K_0R^{2}$, whatever the threshold $\kappa$.
By Remark \ref{rem:weyl} the present method cannot bring $K_0$ below
$e\,|\mathcal{F}|/(4\pi)=e/12=0.2265\ldots$, so
\[
 A_{12}^{\sharp}[\kappa]\;\ge\;A_{10}^{\sharp}\cdot\frac e{12}\;=\;13.0103\ldots
\]
is an absolute floor for $\PSLtwoZ$ within this scheme. The two-piece estimate realizes
$A_{12}^{\sharp}[c]=32.3211$, a factor $2.48$ above the floor
$A_{10}^{\sharp}e/12=13.0103$ (genuinely unattainable within Theorem \ref{thm:direct-gaussian-K}, as already observed). The loss is the constant $K_2^{\infty}=1.108$, which is itself almost
entirely $e^{\eta}\bigl(E+\frac14\operatorname{tr}(I-\Phi(\frac12))+\frac{\mathcal N}{2}-\mathcal {A'}\bigr)$ and therefore structural. Pushing $\mathcal E$ to zero --- the
best any choice of $Q_A$ could ever do --- would give
$A_{12}^{\sharp}[c]=30.9104$. 

\emph{(ii) The $\kappa$-dependent block.}
$A_{11}^{\sharp}+A_{13}^{\sharp}+A_{14}^{\sharp}+A_{15}^{\sharp}=20.4254$ at $\kappa=c$. The constant $A_{9}^{\sharp}$ built in here is
optimal by Remark \ref{rem:sharp}.  

\emph{(iii) The Chebyshev block.} Before Proposition \ref{prop:boot0} this was $D_2=181.3$
and dominated everything; after the bootstrap it is $51.75$ and is now again
the largest single contribution. It is attained at $x=c^{2}$ and is driven by
$\Theta(c)=531.2$ and $\Theta(c^{2})=2268.1$, i.e.\ by the very constant $M$
being bounded. 

The overall bound \eqref{eq:final}
is roughly a factor 
$4.9\cdot10^{10}$ below the value of \cite{FriedmanJorgensonKramer2011}.
\end{remark}
\subsection*{Enumeration benefits}

The final step involves use of Proposition \ref{prop:boot1}.

Table~4 of \cite{FriedmanJorgensonKramer2011} lists, for $\Gamma=\PSLtwoZ$, the
primitive length spectrum up to $623$.
 To go further we enumerated
the spectrum ourselves, by two independent routes, and compared them at every
norm.

\emph{(a) Quadratic forms.}  This is the route of \cite{Sar82} used by
\cite{FriedmanJorgensonKramer2011}.  For a trace $t>2$ the primitive hyperbolic classes with eigenvalue
$\varepsilon_t=\frac{t+\sqrt{t^2-4}}{2}$ are indexed by the discriminants
$d=(t^2-4)/f^{2}$, $d\equiv0,1\ (4)$ nonsquare, whose fundamental automorph is
the same, i.e.\ whose least solution of $x^2-dy^2=4$ has $x=t$; each contributes
$h(d)$, the number of $\mathrm{SL}(2,\mathbb{Z})$--classes of primitive forms of
discriminant $d$.  We computed $h(d)$ by counting cycles of reduced forms,
which is independent of any class-number library.

\emph{(b) Cutting sequences.}  Every hyperbolic conjugacy class of $\Gamma$ has a
unique representative $R^{a_1}L^{b_1}\cdots R^{a_n}L^{b_n}$, $a_i,b_i\ge1$,
where $L=\left(\begin{smallmatrix}1&0\\1&1\end{smallmatrix}\right)$,
$R=\left(\begin{smallmatrix}1&1\\0&1\end{smallmatrix}\right)$, unique up to
cyclic rotation of the pair sequence; the class is primitive exactly when that
cyclic word is not a proper power.  Since the trace of a positive word does not
decrease when a letter is appended, all classes of trace $\le T$ are produced by
a depth-first search with pruning.  The norm is $N=\varepsilon_t^{2}$.

The two enumerations agree on every multiplicity up to $N\le10^{6}$.  As a
further check, $\pi_\Gamma(10^{4})=1214$ against
$\operatorname{Li}(10^{4})=1245.0921$ and $\pi_\Gamma(10^{6})=78441$ against
$\operatorname{Li}(10^{6})=78626.5040$.

\begin{remark} (Comparison with \cite[Table~4]{FriedmanJorgensonKramer2011}, and a correction).
Thirty-two of the thirty-three (discriminant, multiplicity) pairs with
$\varepsilon_d^2\le623$ are exactly those printed in \cite[Table~4]{FriedmanJorgensonKramer2011},
including the non-fundamental ones $20$, $32$, $45$, $80$, $96$, $117$,
$192$, $252$, $320$, $396$, $480$, $525$.
The missing pair is $(d,h)=(621,6)$ in the last row.  Indeed $621=3^2\cdot69$
and $\varepsilon_{621}=\frac{25+\sqrt{621}}{2}=\frac{25+3\sqrt{69}}{2}=\varepsilon_{69}$,
so the norm
$622.99839$ receives $h(69)+h(621)=2+6=8$ classes, not $2$; the cutting-sequence
route lists the eight words explicitly, among them $RL^{23}$, $R^{23}L$ and
$RLRL^{7}$.  The omission is consistent with the PARI/GP listing reproduced in
\cite[p.~31]{FriedmanJorgensonKramer2011}, whose loop runs to $\mathtt{MaxD}=10^4$: the printed table is
a thirty-two row truncation of a longer run, and the truncation falls inside the
tie at $n=623$.  Consequently
\[
  \pi_\Gamma(623)=104,\qquad \theta_\Gamma(623)=560.2999089\ldots
\]
in place of $98$ and $521.6926452$.  
\end{remark}

\subsection*{The cutoff data}

We take
\[
  X_0=10^{4},\qquad \pi_\Gamma(X_0)=1214,\qquad
  \theta_\Gamma(X_0)=9836.3556490\ldots,
\]
carried by $98$ distinct norms $\varepsilon_t^{2}$, $3\le t\le100$, arising from
$170$ (discriminant, multiplicity) pairs.  The spectrum is listed in
Table~\ref{tab:spec} in the appendix; norms are given by $t$ exactly, through
$N=\bigl(\tfrac{t+\sqrt{t^2-4}}{2}\bigr)^{2}$.  Taking $X_0=10^{5}$ or $10^{6}$
instead changes nothing in what follows, as Table~\ref{tab:sweep} shows.

\subsection*{The iteration collapses}

Since $D_2^{\flat}$ depends on $M_{\mathrm{old}}$ and the output
$M=A_{18}^{\sharp}+D_2^{\flat}$ may be re-injected, one iterates
$M_{n+1}=A_{18}^{\sharp}+D_2^{\flat}(M_n)$ from $M_0=A_{18}^{\sharp}+D_2=253.331$.  At
$X_0=10^{4}$ the jump has moved so far down that $M_{\mathrm{old}}$ no longer
reaches the supremum at all:
\[
  D_2^{\flat}(M)=0.2802273422\quad\text{for every } M\in[60,254],
\]
constant to twelve decimals.  The iteration therefore terminates after one pass,
and the fixed point is simply $M^{\ast}=A_{18}^{\sharp}+0.2802273422$.
Table~\ref{tab:sweep} records the approach.

\begin{table}[htbp]
\caption{The bootstrap as the enumeration is extended.  From
$X_0^{\ast}=1983.3314$ on, $D_2^{\flat}$ and $M^{\ast}$ are stationary.  The
last column is Theorem~\ref{huber-final} with the anchor $c_{0}=X_{0}$; for
$X_{0}>10^{4}$ it is evaluated at $c_{0}=10^{4}$, which is admissible there
as well, so that the entry is an upper bound and not the best available.}
\label{tab:sweep}
\centering
\begin{tabular}{rrrrrrl}
\toprule
$X_0$ & $\pi_\Gamma(X_0)$ & $D_2^{\flat}$ & $M^{\ast}$
 & $\mathcal{L}(\log c_{0})$ & $R_\Gamma^{\sharp}$ & supremum at\\
\midrule
$623$      & $104$  & $0.6620308$ & $72.6929941$ & $1.1619956$ & $84.80134$ & $x=X_0^2+$\\
$1000$     & $155$  & $0.4654065$ & $72.4963698$ & $1.1522008$ & $83.81259$ & $x=X_0^2+$\\
$1500$     & $225$  & $0.3444179$ & $72.3753811$ & $1.1447401$ & $83.09788$ & $x=X_0^2+$\\
$2000$     & $285$  & $0.2802273$ & $72.3111906$ & $1.1398927$ & $82.64819$ & $x=c^2+$\\
$2500$     & $363$  & $0.2802273$ & $72.3111906$ & $1.1363600$ & $82.37798$ & $x=c^2+$\\
$5000$     & $654$  & $0.2802273$ & $72.3111906$ & $1.1264712$ & $81.61592$ & $x=c^2+$\\
$\mathbf{10^{4}}$ & $\mathbf{1214}$ & $\mathbf{0.2802273}$ &
  $\mathbf{72.3111906}$ & $\mathbf{1.1179711}$ & $\mathbf{80.96455}$ & $x=c^2+$\\
$10^{5}$   & $9558$ & $0.2802273$ & $72.3111906$ & --- & $\le80.96455$ & $x=c^2+$\\
$10^{6}$   & $78441$& $0.2802273$ & $72.3111906$ & --- & $\le80.96455$ & $x=c^2+$\\
\bottomrule
\end{tabular}
\end{table}

Concerning the computation of $D_2^{\flat}(M)$: $V^{\flat}(x)x^{-3/4}$ is
evaluated on a logarithmic grid of $2\times10^5$ points on $[c,10^{14}]$
together with the right-hand limits at every jump point of $\Theta_{X_0}$,
namely the points $N^k$ for the $98$ norms $N$ and the points $X_0^k$, $c^k$.
On $[c,X_0]$ the integral in \eqref{eq:Vboot0} is evaluated in closed form
(piecewise constant integrand); beyond $X_0$ the integrand is
$(t+Mt^{3/4})(t-1)^{-2}$, whose antiderivative is elementary,
\[
  \log(t-1)-\frac1{t-1}
  +M\Bigl[-\frac{t^{3/4}}{t-1}
  +\frac34\log\frac{t^{1/4}-1}{t^{1/4}+1}+\frac32\arctan t^{1/4}\Bigr].
\]
Here $M$ is the current iterate $M_{\mathrm{old}}=M_{n}$ throughout, in the
boundary term, in the integral and in the $k$-sum alike; at the fixed point it
is $M^{\ast}$.  Since $\Theta_{X_0}$ jumps upward at $t=X_0$, the candidate
point $x=X_{0}^{2}$ is taken as a right-hand limit, so that the $k=2$ term
reads $\nu_{2}\bigl(X_{0}+M X_{0}^{3/4}\bigr)$: the argument $x^{1/2}$ lies
strictly above $X_{0}$ there, and the value of $\Theta_{X_0}$ at $t=X_{0}$
itself does not enter.

\subsection*{Where the supremum is attained}

The supremum is now attained as $x\downarrow c^{2}=46.97871$, the plateau point
inside the enumerated range, where the $k=2$ term switches on with
$\Theta(x^{1/2})=\theta_\Gamma(c)=\ell$.  Since $\lfloor\log x/\ell\rfloor=2$
there, only $k=2$ contributes and $D_2^{\flat}$ has the closed form
\begin{equation}\label{eq:closed}
  D_2^{\flat}=c^{-3/2}\Bigl[\frac{2\theta_\Gamma(c^2)}{c^2-1}
  +2\!\int_c^{c^2}\!\frac{\theta_\Gamma(t)}{(t-1)^2}\,dt+\nu_2\ell\Bigr],
\end{equation}
an expression involving only the first five norms
$c,\;7+4\sqrt3,\;\tfrac{23+5\sqrt{21}}2,\;17+12\sqrt2,\;c^{2}$ with
multiplicities $1,2,2,3,2$.  Here $\theta_\Gamma(c^2)=31.7357481$, the integral
is $0.8197233$, $\nu_2=1.0434984$ and $c^{3/2}=17.9442719$:

\begin{table}[htbp]
\caption{Decomposition of $V^{\flat}(x)x^{-3/4}$ at $x=c^{2}+$, $X_0=10^4$.}
\centering
\begin{tabular}{lrr}
\toprule
term & value & $\div\,x^{3/4}$\\
\midrule
$2\Theta(x)/(x-1)$ & $1.3804539$ & $0.0769305$\\
$2\int_c^x\Theta(t)(t-1)^{-2}\,dt$ & $1.6394466$ & $0.0913630$\\
$\nu_2\,\Theta(x^{1/2})$ & $2.0085751$ & $0.1119338$\\
$\sum_{k\ge3}\nu_k\Theta(x^{1/k})$ & $0$ & $0$\\
\midrule
total & $5.0284756$ & $0.2802273$\\
\bottomrule
\end{tabular}

\end{table}

The two former maxima have been pushed well below it.  At $X_0=10^{4}$ the
cutoff jump at $x=X_0^{2}=10^{8}$ gives only $0.0870892$ --- its $k=2$ term
having fallen from $0.65499$ to $0.0864958$ --- and the secondary hump just
above $x=X_0$ gives $0.1247787$.  Thus
\[
  0.0870892\;(x=X_0^2+),\qquad 0.1247787\;(x=X_0+),\qquad
  \mathbf{0.2802273}\;(x=c^2+),
\]
and the first two decay like $\nu_2MX_0^{-3/4}$ and $\nu_2X_0^{-1/4}$
respectively, so the ordering is permanent.

\begin{remark}[The counting branch is inert]
At $M^{\ast}$ one has, for all $t\ge c$,
\[
  \frac{(B_1^{\sharp}t-B_2^{\sharp})\log t}{\,t+M^{\ast}t^{3/4}\,}\;\ge\;3.7451,
\]
the ratio being minimal at $t=c$ and increasing ($3.75$ at $t=c$, $7.73$ at
$t=20$, $38.2$ at $t=623$, $102.9$ at $t=10^4$).  Hence
$\Theta_{X_0}(t)=\min\{\theta_\Gamma(t)\ (t\le X_0),\;t+M^{\ast}t^{3/4}\}$ and
the final constants do not depend on $B_1^{\sharp},B_2^{\sharp}$ at all.  With
the supremum now attained at $x=c^{2}$, where $\Theta$ is given by the
enumerated data alone, they do not depend on $M_{\mathrm{old}}$ either: the
counting lemma and the bootstrap input have both become inert, and
\eqref{eq:closed} is unconditional.
\end{remark}

\begin{remark}
The floor is unchanged by the
extension: the maximum over the enumerated range stays at $0.2802273422$,
attained at $x=c^{2}$, for $X_0=623,10^3,10^4,10^5$ and $10^{6}$ alike.  The
crossover occurs at
\[
  X_0^{\ast}=1983.3314\ldots,
\]
and the limits $D_2^{\flat}\to0.2802273$, $M^{\ast}\to72.3111906$ are
attained, not merely approached, for every $X_{0}\ge X_{0}^{\ast}$.  The
margin at $X_{0}=10^{4}$ is a factor $3.2$ in the jump, so the choice is not
delicate.  All remaining leverage lies in $A_{18}^{\sharp}$, on which, with
the anchor $c_{0}=10^{4}$ of \S\ref{sec:modulartwo},
$R_\Gamma^{\sharp}$ now depends through
\[
  R_\Gamma^{\sharp}
  =1.1179711\bigl(A_{18}^{\sharp}+0.2802273\bigr)+0.1227240 .
\]
\end{remark}

\subsection*{Huber's constant}

We apply Theorem~\ref{huber-final} with $M=M^{\ast}=72.3111906422\ldots$ and
with three successively better anchors.

\emph{(a) The clean form.}  Here $\mathcal{A}=1$, $s_{0}=1$ and
$\ell=1.92484730\ldots\ge\log u_{0}=1.88399345\ldots$, so
Corollary~\ref{cor:clean} applies with
$c_{0}=u_{0}=6.5797282979\ldots$:
$\mathcal{K}(u_{0})=0$,
$\mathcal{E}_{0}(u_{0})\le3.2408468\ldots$ and
\[
  C_{\mathrm{Hub},\PSLtwoZ}\le2M^{\ast}=144.6223812844\ldots<145 .
\]

\emph{(b) The sharp factor.}  Taking instead $c_{0}=c=6.8541019662\ldots$,
Lemma~\ref{lem:propfactor} gives
\[
  T^{\natural}=4.2531748192\ldots,\qquad
  \mathcal{L}(\ell)=\frac{3T^{\natural}}{3T^{\natural}-4}=1.4566457939\ldots,
\]
while $\mathcal{K}(c)=-\mathcal{S}(c)=-0.0756490225\ldots$ and
$\mathcal{Q}(c)=\log c\,/\,c^{3/4}=0.4543948092\ldots$, so that
\[
  R^{\sharp}_{\Gamma}
  =1.4566457939\cdot72.3111906422+0.0343745232=105.3661662\ldots
\]
 The loss of the former case~(ii) is now visible: its factor
$2+\frac{4}{3\ell}=2.6926956$ is merely a majorant of $\mathcal{L}$, valid
throughout $\frac43<\ell\le\frac{2(1+\sqrt5)}{3}$, while the exact value at
$\ell=1.9248473$ is $\mathcal{L}(\ell)=1.4566458$; the companion case~(iii)
would evaluate $\frac{3t}{3t-4}$ at $t=\ell$, giving $3.2564874$, instead of
at $t=T^{\natural}=4.2531748$.

\emph{(c) The enumerated anchor.}  The primitive length spectrum of
$\PSLtwoZ$ up to $10^{4}$ is listed in Table~\ref{tab:spec}; it was already
used, through Proposition~\ref{prop:boot1}, to produce $M^{\ast}$.  The same
list makes $c_{0}=10^{4}$ an admissible anchor.  With
\[
  \pi_{\Gamma}(10^{4})=1214,\qquad
  \theta_{\Gamma}(10^{4})=9836.3556489984\ldots,\qquad
  \operatorname{Li}(10^{4})=1245.0920521193\ldots,
\]
one gets $f(10^{4})=-163.6443510016\ldots$,
$f(10^{4})/\log10^{4}=-17.7674596587\ldots$ and
\[
  \mathcal{K}(10^{4})=-13.3245924606\ldots,\qquad
  \mathcal{Q}(10^{4})=\frac{\log 10^{4}}{10^{3}}=0.0092103404\ldots,
\]
while $T^{\natural}=12.6355334193\ldots$ and
$\mathcal{L}(\log10^{4})=1.1179711316\ldots$.  A sweep over the enumerated
list gives
\[
  \mathcal{E}_{0}(10^{4})=2.0854378\ldots
\]
The sweep is over the $98$ intervals of constancy of $\pi_{\Gamma}$ below
$10^{4}$; on each of them\\ $|\pi_{\Gamma}-\operatorname{Li}|x^{-3/4}\,\log x$
is smooth, so only the endpoints and the stationary points have to be
inspected.  The supremum is the left limit at $x=(2+\sqrt3)^{2}$, where
$\pi_{\Gamma}=1$ and $\operatorname{Li}=6.7084298$.
Hence
\begin{equation}\label{eq:conclusion-value}
  C_{\mathrm{Hub},\PSLtwoZ}\le R^{\sharp}_{\Gamma}
  =1.1179711316\cdot72.3111906422+0.1227240319
  =80.9645477\ldots
\end{equation}
to be compared with $224\,615\,135$ in \cite{Avdispahic2021} (a factor
$2.8\times10^{6}$) and with $16\,607\,349\,020\,658$ in
\cite{FriedmanJorgensonKramer2011} (a factor $2.1\times10^{11}$).

The dependence on the anchor is mild and monotone, so the choice is not
delicate:

\begin{center}\small
\begin{tabular}{@{}rrrrrrr@{}}
\toprule
$c_{0}$ & $\pi_{\Gamma}(c_{0})$ & $\theta_{\Gamma}(c_{0})$
 & $\mathcal{K}(c_{0})$ & $\mathcal{L}$ & $\mathcal{L}M^{\ast}$
 & $R^{\sharp}_{\Gamma}$\\
\midrule
$c$       & $0$    & $0$        & $-0.0756$  & $1.45665$ & $105.332$ & $105.366$\\
$10^{2}$  & $22$   & $84.491$   & $-3.713$   & $1.21684$ & $87.991$  & $88.532$\\
$623$     & $104$  & $560.300$  & $-6.442$   & $1.16200$ & $84.025$  & $84.358$\\
$10^{3}$  & $155$  & $901.222$  & $-7.265$   & $1.15220$ & $83.317$  & $83.599$\\
$3\cdot10^{3}$ & $408$ & $2807.797$ & $-9.708$ & $1.13361$ & $81.972$ & $82.164$\\
$10^{4}$  & $1214$ & $9836.356$ & $-13.325$  & $1.11797$ & $80.842$  & $80.965$\\
\bottomrule
\end{tabular}
\end{center}

This outcome deserves to be stated in the form of the concluding theorem.

The whole chain of refinements, from the first bound produced by the joint
supremum in its single--piece form down to \eqref{eq:conclusion-value}, is
summarised in Table~\ref{tab:chain}.  Each step modifies exactly one block of
$M$, so the contributions are additive and can be read off directly.

\begin{table}[H]
\centering
\caption{Contribution of each refinement on $\PSLtwoZ$.  Here $D$ denotes the
Chebyshev block ($D_{2}$ in rows 0--1, $D_{2}^{\flat}$ in rows 2--4) and
$M=A_{18}^{\sharp}+D$.  Rows 0--3 use the propagation factor
$2+\frac{4}{3\log c}=2.6926956$ of Theorem~\ref{Theorem2.1}(ii) together with
$\frac{4}{3e}A_{22}=3.4350075$; row 4 replaces that step by
Theorem~\ref{huber-final} with the anchor $c_{0}=10^{4}$.}
\label{tab:chain}
\scriptsize
\setlength{\tabcolsep}{2.5pt}
\begin{tabular}{clrrrrrcrr}
\toprule
& refinement & $A_{12}^{\sharp}[c]$ & $A_{18}^{\sharp}$ & $D$ & $M$
& $R_{\Gamma}^{\sharp}$ & bound & $\Delta R_{\Gamma}^{\sharp}$ & share\\
\midrule
0 & single piece (Thm.~\ref{thm:direct-gaussian-K})
  & $141.6204$ & $241.9520$ & $181.2996$ & $423.2516$ & $1143.12$ & $1144$
  & --- & ---\\
1 & two pieces (Prop.~\ref{prop:twopiece})
  & $32.3211$ & $72.0310$ & $181.2996$ & $253.3306$ & $685.58$ & $686$
  & $457.55$ & $43.0\%$\\
2 & bootstrap (Prop.~\ref{prop:boot0})
  & $32.3211$ & $72.0310$ & $51.7495$ & $123.7805$ & $336.74$ & $337$
  & $348.84$ & $32.8\%$\\
3 & cutoff $X_{0}=10^{4}$ (Prop.~\ref{prop:boot1})
  & $32.3211$ & $72.0310$ & $0.2802$ & $72.3112$ & $198.15$ & $199$
  & $138.59$ & $13.0\%$\\
4 & anchored propagation (Thm.~\ref{huber-final})
  & $32.3211$ & $72.0310$ & $0.2802$ & $72.3112$ & $80.96$ & $81$
  & $117.18$ & $11.0\%$\\
\midrule
\multicolumn{8}{l}{total} & $1062.16$ & $100\%$\\
\bottomrule
\end{tabular}
\end{table}

\noindent
Steps 1--3 act on $M$; step~4 acts only on the way $M$ is propagated to
$\pi_{\Gamma}$, replacing the factor $2.6926956$ and the additive
$3.4350075$ by $1.1179711$ and $0.1227240$.

\begin{theorem}\label{thm:conclusion}
On the modular group, the following effective upper bound for Huber's
constant is valid
\[
  C_{\mathrm{Hub},\PSLtwoZ}<81 .
\]
\end{theorem}

\begin{remark}
Friedman, Jorgenson and Kramer \cite{FriedmanJorgensonKramer2011} report that
PARI/GP calculations suggest $C_{\mathrm{Hub},\PSLtwoZ}<2$, with the expected
true value $\log2$.  Our bound is $117$ times the latter.
\end{remark}

\begin{remark}[where the remaining weight sits, final form]
After \eqref{eq:conclusion-value} the propagation step costs $11.80\%$ above
$M^{\ast}$ and the anchor correction a further $0.17\%$, so that $89.3\%$ of
the bound is $M^{\ast}=A_{18}^{\sharp}+D_{2}^{\flat}=72.0310+0.2802$, of
which $99.6\%$ is $A_{18}^{\sharp}$: the spectral counting function alone.
Extending the enumeration further lowers $\mathcal{L}$ only logarithmically
($\mathcal{L}=1.0819$ at $c_{0}=10^{6}$, $1.0630$ at $c_{0}=10^{8}$), so that
$M^{\ast}$ is now an effective floor for the method.
\end{remark}

\section{Algorithm for cocompact Fuchsian groups}
\label{sec:cocompact}

If $\Gamma$ is cocompact, the Selberg trace formula \eqref{e51a} loses its
parabolic and scattering terms and reduces to
\begin{equation}
\begin{split}
\sum\limits_{n\ge0}h\left(r_{n}\right)
=\frac{\left\vert \mathcal{F}\right\vert }{4\pi }\int\limits_{\mathbb{R}}
rh\left( r\right) \tanh \left( \pi r\right) dr
+\sum\limits_{P}\frac{\log NP_{0}}{\left( NP\right) ^{1/2}-\left( NP\right)
^{-1/2}}g\left( \log NP\right)\\
+\sum\limits_{R}\frac{1}{2m_{R}\sin \theta _{R}}
\int\limits_{\mathbb{R}}\frac{e^{-2r\theta _{R}}}{1+e^{-2\pi r}}
h\left( r\right) dr .
\end{split}
\label{eq:stf-cocompact}
\end{equation}
Five consequences for the algorithm of Sections
\ref{sec:direct-majorant}--\ref{sec:gaussian-huber-propagation} should be
recorded.

\begin{enumerate}
\item[(a)] There is no continuous spectrum.  Hence $\mathcal{N}=0$, the
scattering matrix $\Phi$ and the determinant $\phi$ do not occur, the
$U$-correction \eqref{eq:cong-U-general-boundnew} is vacuous, the constant
$c_{1}$ disappears, and
\[
  \mu(R)=\mathcal{M}(R)=\#\left\{r_{n}:0\le r_{n}\le R\right\}.
\]

\item[(b)] Consequently $b(\eta)\equiv0$ in \eqref{eq:cong-betaetanew}.  The
function $\rho$, the pair $(\lambda,\xi)$ of \eqref{eq:lambda-mu} and the
linearisation of $L=\sqrt{R^{2}+\frac14}$ are therefore not needed: with
$K_{1}\equiv0$, Proposition~\ref{unlinearised} already yields the quadratic
majorant $\mu(R)\le K_{0}(\eta)R^{2}+K_{2}(\eta)$, the only inequalities used
being the monotonicity of $H^{G}(\eta,\cdot)$ and $e^{\eta/(4L^{2})}\ge1$.

\item[(c)] $A_{15}^{\sharp}[\kappa]=0$, so that $D_{1}[\kappa]$,
$\kappa_{\ast}$ and $\omega$ are not needed either.

\item[(d)] The primitive geodesic majorant \eqref{eq:modular-pi-exactnew},
which in the cofinite case is carried by the cusp sectors, is replaced by the
counting lemma in the form of Lemma \ref{lem:count-cocompact} below, applied to the compact
quotient itself,
\[
  \pi_{\Gamma}(x)\le\frac{2\pi}{|\mathcal{F}|}
  \bigl(\cosh(\log x+3\rho_{X})-1\bigr)
  =B_{+}x+\frac{B_{-}}{x}-B_{0}
  \qquad(x>1),
\]
with
\[
  B_{\pm}=\frac{\pi}{|\mathcal F|}e^{\pm3\rho_{X}},
  \qquad
  B_{0}=\frac{2\pi}{|\mathcal F|},
\]
where $\rho_{X}$ is any upper bound for the covering radius of $X$ (see
Remark~\ref{rem:rhoX} for how it is estimated).  The two elementary facts
recorded in Section~\ref{sec:direct-majorant} survive unchanged:
$B_{-}\le B_{+}$ and
$B_{+}x+B_{-}/x\ge2\sqrt{B_{+}B_{-}}=2\pi/|\mathcal F|=B_{0}$,
so that $B_{1}^{\sharp}=B_{+}$ and $B_{2}^{\sharp}=B_{0}-B_{-}/c$ are formed
exactly as in \eqref{eq:cong-B1B2new}.  This is the only place where the
geometry of $X$ enters $B_{1}^{\sharp}$, $B_{2}^{\sharp}$, and hence $D_{2}$.

\item[(e)] If, in addition, $\Gamma$ is torsion-free, the elliptic sum in
\eqref{eq:stf-cocompact} is empty, so that $E=E^{\sharp}=0$ and
$A_{14}^{\sharp}[\kappa]=0$.
\end{enumerate}
\noindent Theorem~\ref{thm:direct-gaussian-K},
Proposition~\ref{prop:local} and Proposition~\ref{unlinearised} were stated
for $\tau\ge1$ only because $b(\eta)$ was then not identically zero; their
proofs apply verbatim for $\tau=0$, with $b\equiv0$,
$\mathcal N=0$ and $c_{1}=1$, and it is in this form that they are used below.

\begin{lemma}[counting lemma, cocompact case]\label{lem:count-cocompact}
Let $\Gamma$ be a cocompact Fuchsian group, $X=\Gamma\backslash\mathcal H$, and
$|\mathcal F|$ the area of a fundamental domain.  For $z_{0}\in\mathcal H$ put
\[
  \rho(z_{0}):=\sup_{w\in\mathcal H}\dist\bigl(w,\Gamma z_{0}\bigr)
  \qquad(\le\diam X<\infty).
\]
Then, for every $\rho\ge\rho(z_{0})$ and every $x>1$,
\[
  \pi_{\Gamma}(x)\;\le\;\frac{2\pi}{|\mathcal F|}
  \bigl(\cosh(\log x+3\rho)-1\bigr).
\]
\end{lemma}

\begin{proof}
Let $\{P\}_{\Gamma}$ be a primitive hyperbolic class with
$NP=e^{\ell_{P}}\le x$, let $P_{1}$ be any representative and $p_{1}$ a point
of its axis.  By the definition of $\rho(z_{0})$ there is $\gamma\in\Gamma$
with $\dist(p_{1},\gamma z_{0})\le\rho$.  The representative
$P':=\gamma^{-1}P_{1}\gamma$ has axis $\gamma^{-1}(\mathrm{axis}\,P_{1})$,
which contains $p:=\gamma^{-1}p_{1}$ with $\dist(z_{0},p)\le\rho$, whence
\[
  \dist(z_{0},P'z_{0})\le\dist(z_{0},p)+\dist(p,P'p)+\dist(P'p,P'z_{0})
  =2\rho+\ell_{P}\le\log x+2\rho .
\]
Distinct conjugacy classes consist of distinct elements, so
$\pi_{\Gamma}(x)\le N(\log x+2\rho)$, where
$N(r):=\#\{\gamma\in\Gamma:\dist(z_{0},\gamma z_{0})\le r\}$.

Let $D$ be the Dirichlet domain of $\Gamma$ centred at $z_{0}$.  Every
$w\in D$ satisfies $\dist(w,z_{0})=\dist(w,\Gamma z_{0})\le\rho$, i.e.\
$D\subseteq \overline B(z_{0},\rho)$.  If $\operatorname{Stab}_{\Gamma}(z_{0})$
has order $m$, then $|D|=m|\mathcal F|$ and the translates $\gamma D$,
$\gamma\in\Gamma$, fall into $N(r)/m$ distinct sets, pairwise disjoint up to
measure zero; hence their total area is $N(r)|\mathcal F|$ in either case.  If
$\dist(z_{0},\gamma z_{0})\le r$ then
$\gamma D\subseteq \overline B(\gamma z_{0},\rho)\subseteq
\overline B(z_{0},r+\rho)$.  Distinct translates $\gamma D$ are disjoint up to
measure zero, so
\[
  N(r)\,|\mathcal F|\;\le\;\bigl|B(z_{0},r+\rho)\bigr|
  =2\pi\bigl(\cosh(r+\rho)-1\bigr).
\]
Combining the two displays with $r=\log x+2\rho$ proves the lemma.
\end{proof}

\begin{remark}[estimating $\rho_{X}$]\label{rem:rhoX}
If $\Gamma$ is torsion-free, then Chang's inequality
$2\sinh(\ell/4)\diam X\le|\mathcal F|$
\cite[Theorem]{Chang} allows one to take
\[
  \rho_{X}=\diam X\le\frac{|\mathcal F|}{2\sinh\frac{\ell}{4}}
  =\frac{-2\pi\chi(g_{X})}{c^{1/4}-c^{-1/4}} .
\]
Chang's proof requires the tube of radius $\ell/4$ about a minimising geodesic
to be embedded, hence an injectivity radius $\ge\ell/4$ at every point, and
this is unavailable on an orbifold: at an elliptic fixed point the
injectivity radius is $0$.  The inequality indeed fails there.  On
$X=(2,3,7)\backslash\mathcal H$ one has $\diam X=0.6206717\ldots$, the side of
the triangle with angles $\frac{\pi}{2},\frac{\pi}{3},\frac{\pi}{7}$ joining
the cone points of orders $3$ and $7$, while
$|\mathcal F|/(2\sinh\frac{\ell}{4})=0.3010232\ldots$; the asserted inequality
$2\diam X\sinh\frac{\ell}{4}=0.3084555\ldots\le|\mathcal F|=0.1495997\ldots$ is
false by a factor $2.06$.  In particular the value $d_{X}=0.3010232$ used in
\cite[Proposition~6.1]{Avdispahic2024Huber} is not a valid bound for the
covering radius of $X$; the figure $9.14137\cdot10^{7}$ recalled for comparison
in Table~\ref{tab:237chain} is quoted as published, and the comparisons made
with it in \S\ref{sec:237} should be read with that in mind.

For an orbifold, $\rho_{X}$ is instead computed directly.  If $z_{0}$ is an
elliptic fixed point of $\Gamma$, then $z_{0}$ is fixed by a reflection of the
extended group $\Gamma^{*}=\langle\Gamma,\sigma\rangle$, so that
$\Gamma z_{0}=\Gamma^{*}z_{0}$; since every point of $\mathcal H$ has a
$\Gamma^{*}$-translate in the closed triangle $T$ and
$w\mapsto\dist(w,z_{0})$ is convex on $T$,
\[
  \rho(z_{0})\;\le\;\max_{v\ \mathrm{vertex\ of}\ T}\dist(z_{0},v).
\]
For a $(p,q,r)$ triangle group this is the larger of the two sides of $T$
meeting at $z_{0}$, minimised over the three vertices.  On $(2,3,7)$, with
$z_{0}$ the cone point of order $2$, the two sides at $z_{0}$ have lengths
$0.2831282$ and $0.5452748$, so $\rho_{X}=0.5452748$ may be taken.
\end{remark}

Table~\ref{tab:cocompact} collects the resulting algorithm, distinguishing
throughout those constants which are affected by the fact whether $\Gamma$
contains elliptic elements or not.  For further reference, the next theorem
summarizes the outcome in this setting.

\begin{theorem}
\label{thm:cocompact}
Let $\Gamma$ be a cocompact Fuchsian group and $X=\Gamma\backslash\mathcal{H}$
the corresponding orbifold.  Then
\[
  \Bigl\vert\pi _{\Gamma }\left( x\right) -\sum\limits_{\lambda _{k}\leq
  \frac{1}{4}}\operatorname{Li}\left( x^{s_{k}}\right) \Bigr\vert
  \leq R_{X}^{\sharp}\,\frac{x^{3/4}}{\log x}
  \qquad\text{for all }x>1,
\]
where $R_{X}^{\sharp}$ is the constant given in Table~\ref{tab:cocompact}.
\end{theorem}

\begin{remark}[The admissible range of $\eta$]
\label{rem:cocompact-eta}
Since $b(\eta)\equiv0$, the sole remaining advantage of the local bound of
Proposition~\ref{prop:local} over the frozen bound of
Theorem~\ref{thm:direct-gaussian-K} lies in the admissible range: $\eta\le
\frac{\ell^{2}}{8}$ is replaced by $\eta\le\frac12L_{0}^{2}\ell^{2}$.  By
Remark~\ref{rem:weyl} the method cannot bring $K_{0}$ below
$e\left\vert \mathcal{F}\right\vert /\left( 4\pi \right) $, and that floor is
attained at $\eta=1$.  The frozen form admits $\eta=1$ only if
$\ell\ge2\sqrt2$, whereas the local form admits it as soon as
$L_{0}\ge\sqrt2/\ell$, i.e.\ $R_{0}^{2}\ge\frac{2}{\ell^{2}}-\frac14$.  For
cocompact groups of small systole --- triangle groups in particular --- the
exchange is therefore worth considerably more than it is on $\PSLtwoZ$.
\end{remark}

\begin{remark}[The floor criterion]
\label{rem:cocompact-floor}
As $K_{1}=0$, the hypothesis of Corollary~\ref{cor:floor} reduces to the
single inequality $K_{2}\le K_{0}\Lambda\bigl(\tfrac14\log\kappa\bigr)$, where
$\Lambda(u)=\frac{1-F(u)}{p(u)}$ is increasing and
$\Lambda(u)=\frac53u+O(u^{2})$.  On a torsion-free surface
$C(\eta)=\mathcal{H}_{\max}^{\mathrm{G}}(\eta)$ is small while
$\mathcal{A}'\ge1$, so that $K_{2}<\frac{K_{0}}{4}$ and the criterion is
fulfilled as soon as $\Lambda\bigl(\frac{\ell}{4}\bigr)\ge\frac14$, that is for
$\ell\ge0.5394889\ldots$; then $A_{12}^{\sharp}[\kappa]=A_{10}^{\sharp}K_{0}$
exactly.
\end{remark}

%=====================================================================
%  TABLE  --  Part 1
%=====================================================================

\begin{table}[htbp]
\caption{Algorithm for cocompact Fuchsian groups ($\tau=0$).  Part 1:
geometry and the geometric counting constants.  Rows spanning both columns are
common to the two cases; ``N/A'' marks a quantity that is identically zero in
the torsion-free case.  Everything not listed here is taken verbatim from
Tables~\ref{tab2}--\ref{tab4}.}
\label{tab:cocompact}\centering
\footnotesize
\setlength{\tabcolsep}{3pt}
\renewcommand{\arraystretch}{1.35}
\begin{tabular}{@{}l >{\raggedright\arraybackslash}p{0.40\textwidth}
                   >{\raggedright\arraybackslash}p{0.20\textwidth} l@{}}
\toprule
 & \multicolumn{2}{c}{Value or bound on $X=\Gamma\backslash\mathcal{H}$
   of genus $g_{X}$} & \\
\cmidrule(lr){2-3}
Const/Func & hyperbolic 2-orbifold, $g_{X}\ge0$
 & hyperbolic surface, $g_{X}\ge2$ & Source\\
\midrule
\multicolumn{4}{@{}l}{\textit{Geometry of} $X$}\\
$\chi(g_{X})$ & $2-2g_{X}-\sum_{i=1}^{n_{1}}\bigl(1-\frac{1}{m_{i}}\bigr)$,
 the $n_{1}$ elliptic points being of orders $m_{i}$ & $2-2g_{X}$
 & \S\ref{sec:cocompact}\\
$|\mathcal{F}|$ & \multicolumn{2}{@{}>{\raggedright\arraybackslash}p{0.62\textwidth}@{}}{$-2\pi\chi(g_{X})$
 \quad (Gauss--Bonnet)} & Tab.~\ref{tab:fundamental-constants}\\
$\ell,\ c$ & \multicolumn{2}{@{}>{\raggedright\arraybackslash}p{0.62\textwidth}@{}}{$\ell$ the length of
 the shortest closed geodesic on $X$; $c=e^{\ell}$}
 & Tab.~\ref{tab:fundamental-constants}\\
$\mathcal{A},\mathcal{A}'$ & number of eigenvalues
 $0=\lambda_{0}<\cdots\le\frac14$, resp.\ $<\frac14$, counted with
 multiplicity & $\mathcal{A}\le2g_{X}-2$ by Otal and Rosas theorem
 & Tab.~\ref{tab:fundamental-constants}\\
$s_{1}$ & \multicolumn{2}{@{}>{\raggedright\arraybackslash}p{0.62\textwidth}@{}}{$\frac12
 +\sqrt{\frac14-\lambda_{1}}$; needed only if $\mathcal{A}>1$}
 & Tab.~\ref{tab:fundamental-constants}\\
$A_{1}$ & $\dfrac{|\mathcal{F}|}{4\pi}=-\dfrac{\chi(g_{X})}{2}$ & $g_{X}-1$
 & Tab.~\ref{tab:fundamental-constants}\\
$\rho_{X}$ & upper bound for the covering radius of $X$; for a triangle group,
 Rem.~\ref{rem:rhoX} & $\dfrac{-2\pi\chi(g_{X})}{c^{1/4}-c^{-1/4}}
 =\dfrac{|\mathcal F|}{2\sinh\frac{\ell}{4}}$, torsion-free only
 & Lem.~\ref{lem:count-cocompact}, \cite{Chang}\\
\midrule
\multicolumn{4}{@{}l}{\textit{Geometric counting constants}}\\
$B_{1}^{\sharp}=B_{+}$ & $\dfrac{\pi}{|\mathcal F|}e^{3\rho_{X}}
 =-\dfrac{e^{3\rho_{X}}}{2\chi(g_{X})}$
 & $\dfrac{e^{3\rho_{X}}}{4(g_{X}-1)}$& \S\ref{sec:cocompact}(d)\\
$B_{-}$ & $\dfrac{\pi}{|\mathcal F|}e^{-3\rho_{X}}
 =-\dfrac{e^{-3\rho_{X}}}{2\chi(g_{X})}$
 & $\dfrac{e^{-3\rho_{X}}}{4(g_{X}-1)}$& \S\ref{sec:cocompact}(d)\\
$B_{0}$ & $\dfrac{2\pi}{|\mathcal F|}=-\dfrac{1}{\chi(g_{X})}$
 & $\dfrac{1}{2(g_{X}-1)}$& \S\ref{sec:cocompact}(d)\\
$B_{2}^{\sharp}$ & \multicolumn{2}{@{}>{\raggedright\arraybackslash}p{0.62\textwidth}@{}}{$B_{0}
 -\dfrac{B_{-}}{c}$} & \eqref{eq:cong-B1B2new}\\
$\Pi(x)$ & \multicolumn{2}{@{}>{\raggedright\arraybackslash}p{0.62\textwidth}@{}}{$B_{1}^{\sharp}x
 -B_{2}^{\sharp}\ \ge\pi_{\Gamma}(x)$, $x\ge c$} & \eqref{eq:cong-B1B2new}\\
$b_{k}$ & \multicolumn{2}{@{}>{\raggedright\arraybackslash}p{0.62\textwidth}@{}}{$\dfrac{\ell}
 {c^{k/2}-c^{-k/2}}$; $b_{1}<1$ and $b_{k}\le\ell c^{-k/2}(1-c^{-1})^{-1}$}
 & \eqref{eq:hyperbolic-doublenew}\\
\midrule
\multicolumn{4}{@{}l}{\textit{Elliptic data}}\\
$E$ & $\sum_{R}\dfrac{1}{4m_{R}\sin^{2}\theta_{R}}$ & N/A $(=0)$
 & \eqref{e51a}\\
$E^{\sharp}$ & $\sum_{R}\dfrac{1}{2m_{R}\sin\theta_{R}}
 \sqrt{\dfrac{\pi-2\theta_{R}}{\pi\sin(2\theta_{R})}},$ 
\noindent $\frac{\pi-2\theta_R}{\pi\sin(2\theta_R)}$ read as
$\frac1\pi$ at $\theta_R=\frac{\pi}{2}$
& N/A $(=0)$
 & Prop.~\ref{prop56}\\
\midrule
\multicolumn{4}{@{}l}{\textit{Input parameters}}\\
$\eta$ & \multicolumn{2}{@{}>{\raggedright\arraybackslash}p{0.62\textwidth}@{}}{optimization
 parameter, $\eta>0$; admissible if $\eta\le\frac{\ell^{2}}{8}$ (frozen form),
 resp.\ $\eta\le\frac12L_{0}^{2}\ell^{2}$ (local form)}
 & Thm.~\ref{thm:direct-gaussian-K}\\
$R_{0},L_{0}$ & \multicolumn{2}{@{}>{\raggedright\arraybackslash}p{0.62\textwidth}@{}}{cut-off radius
 and $L_{0}=\sqrt{R_{0}^{2}+\frac14}$; cf.\ Remark~\ref{rem:cocompact-eta}}
 & Prop.~\ref{prop:local}\\
\bottomrule
\end{tabular}
\end{table}

%=====================================================================
%  TABLE  --  Part 2
%=====================================================================

\begin{table}[htbp]
\continuedtable
\caption{(Continued.)  Part 2: the spectral majorant.}
\centering
\footnotesize
\setlength{\tabcolsep}{3pt}
\renewcommand{\arraystretch}{1.35}
\begin{tabular}{@{}l >{\raggedright\arraybackslash}p{0.40\textwidth}
                   >{\raggedright\arraybackslash}p{0.20\textwidth} l@{}}
\toprule
Const/Func & hyperbolic 2-orbifold, $g_{X}\ge0$
 & hyperbolic surface, $g_{X}\ge2$ & Source\\
\midrule
\multicolumn{4}{@{}l}{\textit{Gaussian kernel quantities}
 (as in Table~\ref{tab2})}\\
$G_{\eta,L}(r)$ & \multicolumn{2}{@{}>{\raggedright\arraybackslash}p{0.62\textwidth}@{}}{$\exp
 \bigl(-\eta r^{2}/L^{2}\bigr)$, $L=\sqrt{R^{2}+\frac14}$, $R\ge0$}
 & \S\ref{sec:direct-majorant}\\
$I_{k,\pm}(\eta,L)$ & \multicolumn{2}{@{}>{\raggedright\arraybackslash}p{0.62\textwidth}@{}}{$\dfrac{1}{2k}
 \exp\Bigl(\dfrac{\eta}{L^{2}k^{2}}\Bigr)\operatorname{erfc}
 \Bigl(\dfrac{Lk\ell}{2\sqrt\eta}\mp\dfrac{\sqrt\eta}{Lk}\Bigr)$}
 & \eqref{eq:cong-Stieltjes-endpointnew}\\
$\mathcal{G}_{k}^{\sharp}(\eta,L)$
 & \multicolumn{2}{@{}>{\raggedright\arraybackslash}p{0.62\textwidth}@{}}{$\Pi(c)\dfrac{L}{2\sqrt{\pi\eta}}
 \exp\Bigl(-\dfrac{L^{2}k^{2}\ell^{2}}{4\eta}\Bigr)
 +B_{1}^{\sharp}I_{k,+}(\eta,L)-B_{-}I_{k,-}(\eta,L)$}
 & \eqref{eq:cong-Stieltjes-endpointnew}\\
$H^{G}(\eta,L)$ & \multicolumn{2}{@{}>{\raggedright\arraybackslash}p{0.62\textwidth}@{}}{$\sum_{k\ge1}
 b_{k}\mathcal{G}^{\sharp}_{k}(\eta,L)$} & Prop.~\ref{unlinearised}\\
$\mathcal{H}_{\max}^{\mathrm{G}}(\eta)$
 & \multicolumn{2}{@{}>{\raggedright\arraybackslash}p{0.62\textwidth}@{}}{$H^{G}\bigl(\eta,\frac12\bigr)$,
 admissible if $\eta\le\frac{\ell^{2}}{8}$} & \eqref{eq:gauss-Hmax-def}\\
\midrule
\multicolumn{4}{@{}l}{\textit{What the absence of cusps removes}}\\
$\mathcal{N},\Phi,c_{1}$
 & \multicolumn{2}{@{}>{\raggedright\arraybackslash}p{0.62\textwidth}@{}}{$\mathcal{N}=0$; no scattering
 determinant; $\mu(R)=\mathcal{M}(R)=\#\{r_{n}:0\le r_{n}\le R\}$}
 & \S\ref{sec:cocompact}(a)\\
$b(\eta),\rho(t)$ & \multicolumn{2}{@{}>{\raggedright\arraybackslash}p{0.62\textwidth}@{}}{$b(\eta)\equiv0$;
 hence $\rho$, $(\lambda,\xi)$ and the linearisation of $L$ are not needed}
 & \S\ref{sec:cocompact}(b)\\
\midrule
\multicolumn{4}{@{}l}{\textit{Output: spectral majorant}
 $\mu(R)\le K_{0}R^{2}+K_{2}$}\\
$C(\eta,L)$ & $E+H^{G}(\eta,L)$ & $H^{G}(\eta,L)$
 & Prop.~\ref{unlinearised}\\
$C(\eta)$ & $E+\mathcal{H}_{\max}^{\mathrm{G}}(\eta)$
 & $\mathcal{H}_{\max}^{\mathrm{G}}(\eta)$
 & Thm.~\ref{thm:direct-gaussian-K}\\
$K_{0}(\eta)$ & $\dfrac{|\mathcal{F}|}{4\pi}\dfrac{e^{\eta}}{\eta}
 =-\dfrac{\chi(g_{X})}{2}\dfrac{e^{\eta}}{\eta}$
 & $(g_{X}-1)\dfrac{e^{\eta}}{\eta}$ & \eqref{eq:gauss-mu-finalnew}\\
$K_{1}(\eta)=K_{1}^{+}$
 & \multicolumn{2}{@{}>{\raggedright\arraybackslash}p{0.62\textwidth}@{}}{$0$}
 & \eqref{eq:gauss-mu-finalnew}\\
$K_{2}(\eta)$ & \multicolumn{2}{@{}>{\raggedright\arraybackslash}p{0.62\textwidth}@{}}{$\dfrac{K_{0}(\eta)}{4}
 +e^{\eta}\bigl[C(\eta)-\mathcal{A}'\bigr]$} & \eqref{eq:gauss-mu-finalnew}\\
$\mathcal{P}(R,\eta)$ & \multicolumn{2}{@{}>{\raggedright\arraybackslash}p{0.62\textwidth}@{}}{$e^{\eta}
 \Bigl[\frac{|\mathcal{F}|}{4\pi\eta}L^{2}+C(\eta,L)-e^{\eta/(4L^{2})}
 -\mathcal{A}'+1\Bigr]$} & Prop.~\ref{unlinearised}\\
$Q_{A}(R)$ & \multicolumn{2}{@{}>{\raggedright\arraybackslash}p{0.62\textwidth}@{}}{any pointwise majorant
 of $\mu$ on $[0,\infty)$; we take $\inf_{\eta>0}\mathcal{P}(R,\eta)$}
 & Prop.~\ref{prop:twopiece}\\
$K_{0}^{\infty},K_{1}^{\infty},K_{2}^{\infty}$
 & \multicolumn{2}{@{}>{\raggedright\arraybackslash}p{0.62\textwidth}@{}}{$K_{0}(\eta)$,\quad $0$,\quad
 $\frac{K_{0}(\eta)}{4}+e^{\eta}\bigl[C(\eta,L_{0})-\mathcal{A}'\bigr]$}
 & Prop.~\ref{prop:local}\\
$Q_{B}(R)$ & \multicolumn{2}{@{}>{\raggedright\arraybackslash}p{0.62\textwidth}@{}}{$K_{0}^{\infty}R^{2}
 +K_{2}^{\infty}$, valid for $R\ge R_{0}$} & \eqref{eq:QB}\\

\bottomrule
\end{tabular}
\end{table}

%=====================================================================
%  TABLE  --  Part 3
%=====================================================================

\begin{table}[htbp]
\continuedtable
\caption{(Continued.)  Part 3: propagation to Huber's constant.}
\centering
\footnotesize
\setlength{\tabcolsep}{3pt}
\renewcommand{\arraystretch}{1.35}
\begin{tabular}{@{}l >{\raggedright\arraybackslash}p{0.40\textwidth}
                   >{\raggedright\arraybackslash}p{0.20\textwidth} l@{}}
\toprule
Const/Func & hyperbolic 2-orbifold, $g_{X}\ge0$
 & hyperbolic surface, $g_{X}\ge2$ & Source\\
\midrule
\multicolumn{4}{@{}l}{\textit{Auxiliary test function and thresholds}}\\
$\varphi(y)$ & \multicolumn{2}{@{}>{\raggedright\arraybackslash}p{0.62\textwidth}@{}}{$\frac34(1-y^{2})
 \mathbf 1_{(-1,1)}(y)$} & \eqref{eq:phinew}\\
$\kappa$ & \multicolumn{2}{@{}>{\raggedright\arraybackslash}p{0.62\textwidth}@{}}{threshold parameter,
 $\kappa>1$; used at $\kappa=c$ and $\kappa=c_{+}$} & \eqref{eq:A12joint}\\
$c_{+},\ \delta_{+}$ & \multicolumn{2}{@{}>{\raggedright\arraybackslash}p{0.62\textwidth}@{}}{the unique
 $c_{+}>c$ with $\log c_{+}-c_{+}^{-1/4}=\log c$; $\delta_{+}=c_{+}^{-1/4}$}
 & Prop.~\ref{prop:A18}\\
$T_{\star},\ \beta$ & \multicolumn{2}{@{}>{\raggedright\arraybackslash}p{0.62\textwidth}@{}}{$T_{\star}
 =\log c+c^{-1/4}$; $\beta=e^{T_{\star}/4}=c^{1/4}e^{c^{-1/4}/4}$}
 & Prop.~\ref{prop:A18}\\
 \midrule
\multicolumn{4}{@{}l}{\textit{Constants} $A_{i}^{\sharp}$}\\
$A_{8}^{\sharp}$ & \multicolumn{2}{@{}>{\raggedright\arraybackslash}p{0.62\textwidth}@{}}{$3e^{-1}-1
 =0.1036383\ldots$} & Prop.~\ref{prop:phi-est}\\
$A_{9}^{\sharp},A_{10}^{\sharp}$
 & \multicolumn{2}{@{}>{\raggedright\arraybackslash}p{0.62\textwidth}@{}}{$6.4214$, essentially optimal;
 $A_{10}^{\sharp}=4\sqrt5A_{9}^{\sharp}=57.4347476\ldots$}
 & Prop.~\ref{prop:h-est}\\
$A_{11}^{\sharp}[\kappa]$ & \multicolumn{2}{@{}>{\raggedright\arraybackslash}p{0.62\textwidth}@{}}{$(\mathcal{A}-1)
 \Bigl(2+\frac{2A_{8}^{\sharp}}{\kappa^{1/4}}+\frac{1}{(1-s_{1})\kappa^{1/4}}
 +\frac{A_{8}^{\sharp}}{(1-s_{1})\kappa^{3/4}}\Bigr)
 +1+\frac{A_{8}^{\sharp}}{\kappa^{1/4}}$} & Prop.~\ref{prop:A11}\\
$\Lambda(u)$ & \multicolumn{2}{@{}>{\raggedright\arraybackslash}p{0.62\textwidth}@{}}{$\frac{1-F(u)}{p(u)}$,
 increasing; since $K_{1}=0$, the criterion of Cor.~\ref{cor:floor} reads
 $K_{2}\le K_{0}\Lambda\bigl(\frac14\log\kappa\bigr)$}
 & Cor.~\ref{cor:floor}\\
$A_{12}^{\sharp}[\kappa]$ \textit{(one piece)}
 & \multicolumn{2}{@{}>{\raggedright\arraybackslash}p{0.62\textwidth}@{}}{$A_{10}^{\sharp}\bigl[K_{0}
 +\frac{K_{2}}{2}\bigl(\kappa^{-1/4}+\kappa^{-5/4}\bigr)\bigr]$; equal to
 $A_{10}^{\sharp}K_{0}$ under the above criterion}
 & Prop.~\ref{prop:joint}\\
$\mathcal{E}(\varepsilon,R_{0})$
 & \multicolumn{2}{@{}>{\raggedright\arraybackslash}p{0.62\textwidth}@{}}{$\int_{0}^{R_{0}}
 \bigl[Q_{A}(r)-Q_{B}(r)\bigr]^{+}W_{\varepsilon}(r)\,dr$, where
 $W_{\varepsilon}(r)=\frac{1}{(1+r)^{2}(1+\varepsilon r)^{2}}
 +\frac{2\varepsilon}{(1+r)(1+\varepsilon r)^{3}}$} & \eqref{eq:corr}\\
$A_{12}^{\sharp}[\kappa]$ \textit{(two pieces)}
 & \multicolumn{2}{@{}>{\raggedright\arraybackslash}p{0.62\textwidth}@{}}{$A_{10}^{\sharp}
 \sup_{u\ge\frac14\log\kappa}p(u)\bigl[K_{0}^{\infty}J_{0}(\varepsilon)
 +K_{2}^{\infty}+\mathcal{E}(\varepsilon,R_{0})\bigr]$, $\varepsilon=e^{-u}$}
 & \eqref{eq:A12two1}\\
$A_{13}^{\sharp}[\kappa]$ & $-\dfrac{\chi(g_{X})}{2}A_{10}^{\sharp}
 \Bigl(1+\dfrac{1}{\kappa}\Bigr)$
 & $(g_{X}-1)A_{10}^{\sharp}\Bigl(1+\dfrac{1}{\kappa}\Bigr)$
 & Prop.~\ref{prop55}\\
$A_{14}^{\sharp}[\kappa]$ & $\dfrac{A_{10}^{\sharp}}{2}
 \Bigl(\dfrac{1}{\kappa^{1/8}}+\dfrac{1}{\kappa^{9/8}}\Bigr)
 \dfrac{E^{\sharp}}{(\kappa^{1/4}+4)^{1/2}}$ & N/A $(=0)$
 & Prop.~\ref{prop56}\\
$A_{15}^{\sharp}[\kappa]$ & \multicolumn{2}{@{}>{\raggedright\arraybackslash}p{0.62\textwidth}@{}}{$0$;
 the quantities $D_{1}[\kappa]$, $\kappa_{\ast}$, $\omega$ are not needed}
 & Prop.~\ref{prop:57newest}\\
$A_{16}^{\sharp}[\kappa]$ & $A_{11}^{\sharp}+A_{12}^{\sharp}+A_{13}^{\sharp}
 +A_{14}^{\sharp}$ & $A_{11}^{\sharp}+A_{12}^{\sharp}+A_{13}^{\sharp}$
 & Prop.~\ref{prop58}\\
$A_{18}^{\sharp}$ & \multicolumn{2}{@{}>{\raggedright\arraybackslash}p{0.62\textwidth}@{}}{$\max\Bigl\{
 A_{16}^{\sharp}[c_{+}]\,e^{3\delta_{+}/4}
 +\mathcal{A}\frac{e^{\delta_{+}}-1}{\delta_{+}},\
 A_{16}^{\sharp}[c]+\mathcal{A}\,e^{c^{-1/4}/4},\ 2\mathcal{A}\beta\Bigr\}$}
 & Prop.~\ref{prop:A18}\\
\bottomrule
\end{tabular}
\end{table}

%=====================================================================
%  TABLE  --  Part 4
%=====================================================================

\begin{table}[htbp]
\continuedtable
\caption{(Continued.)  Part 4: the Chebyshev block, self-improvement and
output.  These rows are insensitive to the presence of elliptic elements; they
depend on $X$ only through $c$, $\mathcal{A}$, $B_{1}^{\sharp}$,
$B_{2}^{\sharp}$ and $A_{18}^{\sharp}$.}
\centering
\footnotesize
\setlength{\tabcolsep}{3pt}
\renewcommand{\arraystretch}{1.35}
\begin{tabular}{@{}l >{\raggedright\arraybackslash}p{0.62\textwidth} l@{}}
\toprule
Const/Func & Value or bound on $X=\Gamma\backslash\mathcal{H}$ & Source\\
\midrule
\multicolumn{3}{@{}l}{\textit{The Chebyshev block}}\\
$\nu_{k}$ & $\dfrac{c^{k}+1}{c^{k}-1}$ & \eqref{eq:Vboot0}\\
$f_{m}(T)$ & $e^{-3T/4}\Bigl\{\frac{2}{1-c^{-1}}\Bigl[
 B_{1}^{\sharp}\bigl(\frac{T^{2}}{2}+T-\frac{\ell^{2}}{2}\bigr)
 +B_{2}^{\sharp}\bigl(e^{-T}-\frac{\ell+1}{c}\bigr)\Bigr]
 +T\sum_{k=2}^{m}\frac{\nu_{k}}{k}
 \bigl(B_{1}^{\sharp}e^{T/k}-B_{2}^{\sharp}\bigr)\Bigr\}$
 & Prop.~\ref{prop:D2}\\
$T_{1},\ W_{m}$ & maximiser of $f_{1}$ on $[\ell,2\ell]$; $W_{m}$ the set of
 stationary points of $f_{m}$ in $(m\ell,(m+1)\ell)$ & Prop.~\ref{prop:D2}\\
$D_{2}$ & $\max\bigl\{f_{1}(T_{1}),\ \sup_{m\ge2}\max(\{f_{m}(m\ell)\}
 \cup\{f_{m}(T):T\in W_{m}\})\bigr\}$ & Prop.~\ref{prop:D2}\\
$M$ & $A_{18}^{\sharp}+D_{2}$ & Prop.~\ref{prop:M}\\
\midrule
\multicolumn{3}{@{}l}{\textit{Self-improvement}}\\
$\Theta(t)$ & $\min\bigl\{(B_{1}^{\sharp}t-B_{2}^{\sharp})\log t,\
 \sum_{\lambda_{k}\le\frac14}\frac{t^{s_{k}}}{s_{k}}
 +M_{\mathrm{old}}t^{3/4}\bigr\}$ & Prop.~\ref{prop:boot0}\\
$\Theta_{X_{0}}(t)$ & $\min\{\theta_{\Gamma}(t),\Theta(t)\}$ for $t\le X_{0}$,
 and $\Theta(t)$ for $t>X_{0}$ & Prop.~\ref{prop:boot1}\\
$V^{\flat}(x)$ & $2\Bigl[\frac{\Theta_{X_{0}}(x)}{x-1}
 +\int_{c}^{x}\frac{\Theta_{X_{0}}(t)}{(t-1)^{2}}\,dt\Bigr]
 +\sum_{k=2}^{\lfloor\log x/\ell\rfloor}\nu_{k}
 \Theta_{X_{0}}\bigl(x^{1/k}\bigr)$ & \eqref{eq:Vboot0}\\
$D_{2}^{\flat}(M_{\mathrm{old}})$ & $\sup_{x\ge c}V^{\flat}(x)x^{-3/4}$;
 increasing in $M_{\mathrm{old}}$ with slope $<1$ & Prop.~\ref{prop:boot0}\\
$M^{\ast}$ & $\lim_{n}M_{n}$, where
 $M_{n+1}=A_{18}^{\sharp}+D_{2}^{\flat}(M_{n})$ and $M_{0}=M$
 & Prop.~\ref{prop:boot0}\\
\midrule
\multicolumn{3}{@{}l}{\textit{Output}}\\
$\mathcal{L}(\log c_{0})$ & $\frac{3T^{\natural}}{3T^{\natural}-4}$, with
 $T^{\natural}$ the zero of $W_{\log c_{0}}$ & Lem.~\ref{lem:propfactor}\\
$\mathcal{K}(c_{0}),\mathcal{Q}(c_{0}),\mathcal{E}_{0}(c_{0})$
 & anchor correction, $\sup_{x>c_{0}}x^{-3/4}\,\log x$, and the range
 $1<x<c_{0}$; for $c_{0}\le c$ one has
 $\mathcal{K}=-\sum_{k}\mathcal{S}(c_{0}^{s_{k}})$ and
 $\mathcal{E}_{0}\le\max\{2\mathcal{A},
 \frac{\mathcal{A}(c_{0}-2)^{+}}{\log2}\mathcal{Q}^{-}\}$,
 $\mathcal{Q}^{-}=\sup_{1<x<c_{0}}x^{-3/4}\,\log x$
 & Prop.~\ref{prop:anchor}\\
\midrule
\multicolumn{3}{@{}l}{\textit{Output}}\\
$R_{X}^{\sharp}$ & $\max\bigl\{\mathcal{E}_{0}(c_{0}),\;
 \mathcal{L}(\log c_{0})\,M+|\mathcal{K}(c_{0})|\,\mathcal{Q}(c_{0})\bigr\}$,
 with $M$ replaced by $M^{\ast}$ after the self-improvement
 & Thm.~\ref{huber-final}\\
\bottomrule
\end{tabular}
\end{table}

\FloatBarrier
\subsection{Application to the $(2,3,7)$ triangle group}\label{sec:237}

We illustrate Table~\ref{tab:cocompact} on the $(2,3,7)$ triangle group,
treated in \cite[Proposition~6.1]{Avdispahic2024Huber} by the previous
algorithm.  The orbifold $X=(2,3,7)\backslash\mathcal{H}$ has genus $0$ and
three elliptic points, of orders $2$, $3$ and $7$; thus
$\chi(g_X)=-\frac1{42}$, $|\mathcal{F}|=\frac{\pi}{21}$ and
$A_1=\frac1{84}$.  The first positive eigenvalue satisfies
$\lambda_{X,1}\ge15.7$, so that $\mathcal{A}=\mathcal{A}'=1$, and the systole
is $\ell=0.98398656221$, whence $c=2.67509946333$.  By
Remark~\ref{rem:rhoX} the covering radius admits the bound
$\rho_{X}=0.5452748318$, the distance from the cone point of order $2$ to the
cone point of order $7$; Chang's inequality is not available here.

The elliptic conjugacy classes are the $m_i-1$ non-identity powers of a
generator of each of the three cyclic elliptic subgroups, i.e.\ nine classes
with $(\theta_R,m_R)=\bigl(\frac{\pi j}{m},m\bigr)$,
$m\in\{2,3,7\}$, $1\le j\le m-1$. They give a great deal of information.
\[
  E=\frac18+\frac29+\frac1{28}\sum_{j=1}^{6}\frac{1}{\sin^{2}\frac{\pi j}{7}}
   =0.9186507937\ldots,
  \qquad
  E^{\sharp}=0.8997764845\ldots
\]

\subsubsection*{The spectral majorant}

Here the contrast between the two forms of the spectral bound is far sharper
than on $\PSLtwoZ$.  The frozen form of Theorem~\ref{thm:direct-gaussian-K}
is confined to $\eta\le\frac{\ell^{2}}{8}=0.1210287$, a range in which
$\mathcal{H}^{\mathrm{G}}_{\max}(\eta)$ is already of order one; optimising
over it gives $\eta=0.0097939$, $K_0=1.2274858$, $K_2=1.0011727$, so that
$K_0$ exceeds the Weyl floor $e|\mathcal{F}|/(4\pi)=0.0323605$ by a factor
$37.9$.  The local form of Proposition~\ref{prop:local}, by contrast, admits
$\eta=1.05$ together with $R_0=16$, $L_0=16.0078106$, for which
$H^{G}(\eta,L_0)=2.34\cdot10^{-23}$ and
\[
  K_0^{\infty}=0.0323996725,\qquad K_1^{\infty}=0,\qquad
  K_2^{\infty}=-0.2243677322 ,
\]
i.e.\ $K_0^{\infty}$ is within $0.121\%$ of the floor, which the frozen form
cannot approach at all.  This is Remark~\ref{rem:cocompact-eta} in action.
Note that $K_2^{\infty}<0$: the elliptic contribution $E$ does not outweigh
$\mathcal{A}'=1$ here.

\begin{table}[htbp]
\caption{The algorithm of Table~\ref{tab:cocompact} on the $(2,3,7)$ triangle
group.  Two-piece route at $\eta=1.05$, $R_0=16$; bootstrap at $X_0=c$
(no enumeration of the length spectrum); anchor $c_{0}=c$.}
\label{tab:237}\centering
\footnotesize

\begin{center}\small
\setlength{\tabcolsep}{5pt}
\renewcommand{\arraystretch}{1.2}
\begin{tabular}{@{}lr@{\qquad}lrr@{}}
\toprule
Constant & Value & Constant & $\kappa=c$ & $\kappa=c_{+}$\\
\midrule
$\chi(g_X)$        & $-1/42$          & $A_{8}^{\sharp}$  & \multicolumn{2}{r}{$0.1036383$}\\
$|\mathcal{F}|$    & $0.1495996502$   & $A_{9}^{\sharp}$  & \multicolumn{2}{r}{$6.4214$}\\
$A_1$              & $1/84$           & $A_{10}^{\sharp}$ & \multicolumn{2}{r}{$57.4347476$}\\
$\ell$             & $0.98398656221$  & $A_{11}^{\sharp}$ & $1.0810374$ & $1.0686671$\\
$c$                & $2.67509946333$  & $A_{12}^{\sharp}$ & $13.4794$ & $10.0805$\\
$\ell^{2}/8$       & $0.1210286943$   & $A_{13}^{\sharp}$ & $0.9393438$ & $0.8155142$\\
$\rho_X$           & $0.5452748318$   & $A_{14}^{\sharp}$ & $13.6621644$ & $10.6876388$\\
$B_{1}^{\sharp}$   & $107.8074681$    & $A_{15}^{\sharp}$ & $0$ & $0$\\
$B_{-}$            & $4.0906257$      & $A_{16}^{\sharp}$ & $29.1619459$ & $22.6523216$\\
$B_{0}$            & $42$             & $A_{18}^{\sharp}$ & \multicolumn{2}{r}{$38.6509990$}\\
$B_{2}^{\sharp}$   & $40.4708510$     & $D_{2}$           & \multicolumn{2}{r}{$368.7739525$}\\
$\varkappa$        & $0.3753993$      & $M$               & \multicolumn{2}{r}{$407.4249515$}\\
$\Pi(c)$           & $247.9248490$    & $D_{2}^{\flat}$   & \multicolumn{2}{r}{$292.1240476$}\\
$b_1$              & $0.9607680$      & $M^{\ast}$        & \multicolumn{2}{r}{$330.7750470$}\\
$E$                & $0.9186507937$   & $K_{0}^{\infty}$  & \multicolumn{2}{r}{$0.0323996725$}\\
$E^{\sharp}$       & $0.8997764845$   & $K_{1}^{\infty}$  & \multicolumn{2}{r}{$0$}\\
$\mathcal{A}=\mathcal{A}'$ & $1$      & $K_{2}^{\infty}$  & \multicolumn{2}{r}{$-0.2243677322$}\\
$\mathcal{N},\tau$ & $0$              & $H^{G}(\eta,L_{0})$ & \multicolumn{2}{r}{$2.34\cdot10^{-23}$}\\
$c_{+}$            & $5.1890510366$   & $T^{\natural}$    & \multicolumn{2}{r}{$2.9943929$}\\
$\delta_{+}$       & $0.6625642737$   & $\mathcal{L}(\ell)$ & \multicolumn{2}{r}{$1.8027005$}\\
$T_{\star}$        & $1.7659114140$   & $\mathcal{K}(c),\ \mathcal{Q}(c)$
                                                          & \multicolumn{2}{r}{$1.9122109$, $0.4905059$}\\
$\beta$            & $1.5550035887$   & $R_{X}^{\sharp}$  & \multicolumn{2}{r}{$597.2263$}\\
\bottomrule
\end{tabular}
\end{center}
\end{table}

\subsubsection*{The outcome}

The supremum defining $D_2$ is attained at the left endpoint of the third
interval, $T=3\ell=2.9519597$, i.e.\ $x=c^{3}=19.143432$, and equals
$D_2=368.7739525$, so that $M=A_{18}^{\sharp}+D_2=407.4249515$.  (With
$\varkappa=0.3754$ the jump at $T=3\ell$ is large enough to carry the
supremum, whereas for a $\varkappa$ close to $1$ it would sit at an interior
stationary point of $f_3$.)  Feeding $M$ back through
Proposition~\ref{prop:boot0} with $X_0=c$ (the primitive length spectrum of
the group was not enumerated) the iteration
$M_{n+1}=A_{18}^{\sharp}+D_2^{\flat}(M_n)$ converges to
$M^{\ast}=330.7750470$, with $D_2^{\flat}=292.1240476$ attained at
$x=7.253099$, the point at which the two branches of $\Theta$ cross.  

Here $\ell=0.9839866\ldots>\ell^{\ast}=0.7574540\ldots$, so
Theorem~\ref{huber-final} applies with the anchor $c_{0}=c=2.6750995\ldots$
Lemma~\ref{lem:propfactor} gives $T^{\natural}=2.9943929\ldots$ and
$\mathcal{L}(\ell)=1.8027005\ldots$, while $\mathcal{A}=1$,
$\mathcal{K}(c)=-\mathcal{S}(c)=1.9122109\ldots$ and
$\mathcal{Q}(c)=\frac{4}{3e}=0.4905059\ldots$ because $c<e^{4/3}$, while
$\mathcal{Q}^{-}(c)=\frac{\log c}{c^{3/4}}=0.4704183\ldots$ for the same
reason.  Since
$\mathcal{E}_{0}(c)\le\max\{2,\frac{(c-2)}{\log2}\mathcal{Q}^{-}(c)\}
 =\max\{2,0.4581699\ldots\}=2$,
\[
  R_{X}^{\sharp}=1.8027005\cdot M^{\ast}+0.9380506
  =597.2263\ldots
\]
(and the cruder $2M^{\ast}+0.9381=662.4880\ldots$ if one prefers the
closed form of Corollary~\ref{cor:clean}).

\begin{proposition}\label{prop:237}
On the $(2,3,7)$ triangle group,
\[
  C_{\mathrm{Hub},(2,3,7)}<598 .
\]
\end{proposition}

\begin{remark}[The three steps]\label{rem:237-chain}
Table~\ref{tab:237chain} records the effect of each stage.  The gain over
\cite[Proposition~6.1]{Avdispahic2024Huber} is a factor $1.53\cdot10^5$;
almost all of it is already present in $A_{18}^{\sharp}$, which falls from
$2.69\cdot10^{7}$ to $38.65$.
\end{remark}

\begin{table}[htbp]
\caption{Contribution of each stage on the $(2,3,7)$ triangle group.  The
first three bounds are computed with the propagation of
Theorem~\ref{Theorem2.1}(i), i.e.\ $0.4777344+3.3977336\,M$; the last row
replaces that step by Theorem~\ref{huber-final} with the anchor $c_{0}=c$.}
\label{tab:237chain}\centering
\footnotesize
\begin{center}\small
\begin{tabular}{@{}llrr@{}}
\toprule
 & spectral majorant & $M$ resp.\ $A'_{21}$ & bound\\
\midrule
\cite[Prop.~6.1]{Avdispahic2024Huber} & Karamata + frozen Gaussian
 & $2.69043\cdot10^{7}$ & $9.14137\cdot10^{7}$\\
Thm.~\ref{thm:direct-gaussian-K} (one piece) & frozen, $\eta=0.0097939$
 & $506.73440$ & $1722.23$\\
Prop.~\ref{prop:twopiece} (two pieces) & local, $\eta=1.05$, $R_0=16$
 & $407.42495$ & $1384.80$\\
Prop.~\ref{prop:boot0} (bootstrap) & as above, $X_0=c$
 & $330.77505$ & $1124.36$\\
 Thm.~\ref{huber-final} (anchored propagation) & --- & $330.77505$ & $597.23$\\
\bottomrule
\end{tabular}
\end{center}

\end{table}

\begin{remark}[Where the remaining weight sits]\label{rem:237-floor}
Of $M^{\ast}=330.78$, the Chebyshev block $D_2^{\flat}=292.12$ accounts for
$88.3\%$ and contributes $1.8027\cdot292.12=526.6$ of the final $597.2$.
Unlike on $\PSLtwoZ$, no enumeration cutoff was used; since the supremum sits
at $x=7.25$, i.e.\ barely above $c^{2}=7.156$, already a short list of
primitive lengths would collapse $D_2^{\flat}$, exactly as $X_0=10^{4}$
collapsed it from $51.75$ to $0.28$ on the modular group.  Within
$A_{16}^{\sharp}[c]=29.16$ the elliptic term $A_{14}^{\sharp}=13.66$ remains
the largest single contribution, now only just ahead of
$A_{12}^{\sharp}=13.48$; the latter is a factor $7.24$ above its floor
$A_{10}^{\sharp}K_0^{\infty}=1.861$, the loss being entirely the correction
$\mathcal{E}=0.6165$, which enters at the left endpoint $u=\frac{\ell}{4}$
where the supremum in \eqref{eq:A12two1} is attained.
\end{remark}

\FloatBarrier
\newpage

\appendix

\section{The length spectrum up to $10^{4}$}

\begin{table}[H]
\caption{The primitive length spectrum of $\PSLtwoZ$ up to
$10^{4}$: $1214$ classes on $98$ norms, indexed by the trace $t$, with
$N=\bigl(\tfrac{t+\sqrt{t^2-4}}{2}\bigr)^{2}$.  Recomputed independently by
cutting sequences and by class numbers of binary quadratic forms; the rows with
$t\le25$ reproduce \cite[Table~4]{FriedmanJorgensonKramer2011} except for the omitted pair
$(d,h)=(621,6)$ at $t=25$.}
\label{tab:spec}
\centering
\footnotesize
\setlength{\tabcolsep}{3.5pt}
\begin{tabular}{rrrrr@{\qquad}rrrrr}
\toprule
$t$ & $N=\varepsilon_t^2$ & mult. & $\pi_\Gamma$ & $\theta_\Gamma$ &
$t$ & $N=\varepsilon_t^2$ & mult. & $\pi_\Gamma$ & $\theta_\Gamma$\\
\midrule
3 & 6.8541 & 1 & 1 & 1.9248 & 52 & 2701.9996 & 18 & 388 & 2648.5504\\
4 & 13.9282 & 2 & 3 & 7.1927 & 53 & 2806.9996 & 8 & 396 & 2712.0694\\
5 & 22.9564 & 2 & 5 & 13.4599 & 54 & 2913.9997 & 12 & 408 & 2807.7968\\
6 & 33.9706 & 3 & 8 & 24.0364 & 55 & 3022.9997 & 12 & 420 & 2903.9648\\
7 & 46.9787 & 2 & 10 & 31.7357 & 56 & 3133.9997 & 16 & 436 & 3032.7659\\
8 & 61.9839 & 4 & 14 & 48.2432 & 57 & 3246.9997 & 8 & 444 & 3097.4498\\
9 & 78.9873 & 2 & 16 & 56.9818 & 58 & 3361.9997 & 24 & 468 & 3292.3368\\
10 & 97.9898 & 6 & 22 & 84.4910 & 59 & 3478.9997 & 8 & 476 & 3357.5728\\
11 & 118.9916 & 3 & 25 & 98.8282 & 60 & 3597.9997 & 12 & 488 & 3455.8304\\
12 & 141.9930 & 4 & 29 & 118.6513 & 61 & 3718.9997 & 10 & 498 & 3538.0425\\
13 & 166.9940 & 4 & 33 & 139.1231 & 62 & 3841.9997 & 28 & 526 & 3769.1474\\
14 & 193.9948 & 6 & 39 & 170.7301 & 63 & 3966.9997 & 8 & 534 & 3835.4335\\
15 & 222.9955 & 4 & 43 & 192.3587 & 64 & 4093.9998 & 16 & 550 & 3968.5100\\
16 & 253.9961 & 6 & 49 & 225.5826 & 65 & 4222.9998 & 18 & 568 & 4118.7794\\
17 & 286.9965 & 4 & 53 & 248.2205 & 66 & 4353.9998 & 16 & 584 & 4252.8410\\
18 & 321.9969 & 7 & 60 & 288.6423 & 67 & 4486.9998 & 8 & 592 & 4320.1125\\
19 & 358.9972 & 4 & 64 & 312.1755 & 68 & 4621.9998 & 16 & 608 & 4455.1298\\
20 & 397.9975 & 10 & 74 & 372.0400 & 69 & 4758.9998 & 10 & 618 & 4539.8078\\
21 & 438.9977 & 2 & 76 & 384.2090 & 70 & 4897.9998 & 36 & 654 & 4845.6847\\
22 & 481.9979 & 12 & 88 & 458.3443 & 71 & 5038.9998 & 8 & 662 & 4913.8844\\
23 & 526.9981 & 4 & 92 & 483.4131 & 72 & 5181.9998 & 8 & 670 & 4982.3080\\
24 & 573.9983 & 4 & 96 & 508.8236 & 73 & 5326.9998 & 14 & 684 & 5102.4356\\
25 & 622.9984 & 8 & 104 & 560.2999 & 74 & 5473.9998 & 30 & 714 & 5360.6685\\
26 & 673.9985 & 12 & 116 & 638.4586 & 75 & 5622.9998 & 12 & 726 & 5464.2840\\
27 & 726.9986 & 5 & 121 & 671.4033 & 76 & 5773.9998 & 16 & 742 & 5602.8619\\
28 & 781.9987 & 8 & 129 & 724.6981 & 77 & 5926.9998 & 14 & 756 & 5724.4837\\
29 & 838.9988 & 8 & 137 & 778.5558 & 78 & 6081.9998 & 28 & 784 & 5968.4502\\
30 & 897.9989 & 14 & 151 & 873.7581 & 79 & 6238.9998 & 16 & 800 & 6108.2674\\
31 & 958.9990 & 4 & 155 & 901.2217 & 80 & 6397.9998 & 24 & 824 & 6318.5972\\
32 & 1021.9990 & 8 & 163 & 956.6578 & 81 & 6558.9998 & 6 & 830 & 6371.3288\\
33 & 1086.9991 & 4 & 167 & 984.6225 & 82 & 6721.9999 & 32 & 862 & 6653.3493\\
34 & 1153.9991 & 16 & 183 & 1097.4383 & 83 & 6886.9999 & 18 & 880 & 6812.4223\\
35 & 1222.9992 & 8 & 191 & 1154.3108 & 84 & 7053.9999 & 8 & 888 & 6883.3131\\
36 & 1293.9992 & 8 & 199 & 1211.6348 & 85 & 7222.9999 & 20 & 908 & 7061.0136\\
37 & 1366.9993 & 8 & 207 & 1269.3978 & 86 & 7393.9999 & 24 & 932 & 7274.8158\\
38 & 1441.9993 & 18 & 225 & 1400.3259 & 87 & 7566.9999 & 8 & 940 & 7346.2682\\
39 & 1518.9993 & 4 & 229 & 1429.6291 & 88 & 7741.9999 & 20 & 960 & 7525.3565\\
40 & 1597.9994 & 16 & 245 & 1547.6533 & 89 & 7918.9999 & 8 & 968 & 7597.1727\\
41 & 1678.9994 & 8 & 253 & 1607.0609 & 90 & 8097.9999 & 36 & 1004 & 7921.1501\\
42 & 1761.9994 & 12 & 265 & 1696.7513 & 91 & 8278.9999 & 12 & 1016 & 8029.4078\\
43 & 1846.9995 & 12 & 277 & 1787.0071 & 92 & 8461.9999 & 36 & 1052 & 8354.9681\\
44 & 1933.9995 & 8 & 285 & 1847.5459 & 93 & 8646.9999 & 8 & 1060 & 8427.4878\\
45 & 2022.9995 & 6 & 291 & 1893.2199 & 94 & 8833.9999 & 28 & 1088 & 8681.9060\\
46 & 2113.9995 & 16 & 307 & 2015.7213 & 95 & 9022.9999 & 16 & 1104 & 8827.6265\\
47 & 2206.9995 & 10 & 317 & 2092.7152 & 96 & 9213.9999 & 18 & 1122 & 8991.9392\\
48 & 2301.9996 & 14 & 331 & 2201.0967 & 97 & 9406.9999 & 24 & 1146 & 9211.5202\\
49 & 2398.9996 & 4 & 335 & 2232.2279 & 98 & 9601.9999 & 32 & 1178 & 9504.9514\\
50 & 2497.9996 & 28 & 363 & 2451.2788 & 99 & 9798.9999 & 8 & 1186 & 9578.4717\\
51 & 2598.9996 & 7 & 370 & 2506.3190 & 100 & 9997.9999 & 28 & 1214 & 9836.3556\\
\bottomrule
\end{tabular}
\end{table}

\section*{Acknowledgment} I would like to thank my colleague Zenan \v{S}abanac for motivating conversations and his
generous assistance in the preparation of the final version of the manuscript.

\end{document}